\documentclass[10pt]{amsart}
\usepackage{pdfpages}

\usepackage{graphicx,cancel,xcolor}
\usepackage{amsmath,amsthm,amssymb}
\usepackage{a4wide}
\usepackage[utf8]{inputenc}
\usepackage{graphicx}
\usepackage{comment}
\usepackage{dsfont}
\usepackage{float}
\usepackage{caption} 
\usepackage{subcaption}
\makeatletter
\@namedef{r@tocindent4}{0pt}
\@namedef{r@tocindent5}{0pt}
\makeatother
\makeatletter
\renewcommand{\tocsection}[3]{%
	\indentlabel{\@ifnotempty{#2}{\bfseries\ignorespaces#1 #2\quad}}\bfseries#3}
\renewcommand{\tocsubsection}[3]{%
	\indentlabel{\@ifnotempty{#2}{\ignorespaces#1 #2\quad}}#3}
\newcommand\@dotsep{4.5}
\def\@tocline#1#2#3#4#5#6#7{\relax
	\ifnum #1>\c@tocdepth 
	\else
	\par \addpenalty\@secpenalty\addvspace{#2}%
	\begingroup \hyphenpenalty\@M
	\@ifempty{#4}{%
		\@tempdima\csname r@tocindent\number#1\endcsname\relax
	}{%
		\@tempdima#4\relax
	}%
	\parindent\z@ \leftskip#3\relax \advance\leftskip\@tempdima\relax
	\rightskip\@pnumwidth plus1em \parfillskip-\@pnumwidth
	#5\leavevmode\hskip-\@tempdima{#6}\nobreak
	\leaders\hbox{$\m@th\mkern \@dotsep mu\hbox{.}\mkern \@dotsep mu$}\hfill
	\nobreak
	\hbox to\@pnumwidth{\@tocpagenum{\ifnum#1=1\bfseries\fi#7}}\par
	\nobreak
	\endgroup
	\fi}
\AtBeginDocument{%
	\expandafter\renewcommand\csname r@tocindent0\endcsname{0pt}
}
\def\l@subsection{\@tocline{2}{0pt}{2.5pc}{5pc}{}}
\def\l@subsubsection{\@tocline{2}{0pt}{4.5pc}{5pc}{}}
\def\l@paragraph{\@tocline{2}{0pt}{5.5pc}{5pc}{}}
\makeatother

\def\rr{\mathbb R}

\providecommand{\abs}[1]{\lvert#1\rvert}

\theoremstyle{definition}

\theoremstyle{remark}
\newtheorem{remark}{Remark}

\numberwithin{equation}{section}

\numberwithin{equation}{section}

\begin{document}
\thispagestyle{empty}

\title[Stabilization of locally coupled wave equations: a numerical study in 1D]{Numerical study of the stabilization of \\ 1D locally coupled wave equations.}
\author{St\'ephane Gerbi}
\address{Laboratoire de Math\'ematiques UMR 5127 CNRS \& Universit\'e de Savoie Mont Blanc, Campus scientifique, 73376 Le Bourget du Lac Cedex, France}\email{stephane.gerbi@univ-smb.fr}
\author{Chiraz Kassem}
\address{Universit\'e Libanaise,
			EDST, Equipe EDP-AN,
			Hadath, Beirut, Lebanon}
			\email{shiraz.kassem@hotmail.com}
			\author{Amina Mortada}
			\address{ Universit\'e Libanaise,
			EDST, Equipe EDP-AN,
			Hadath, Beirut, Lebanon}
			\email{amina$\_$mortada2010@hotmail.com}
			\author{Ali Wehbe}
			\address{Universit\'e Libanaise,
			Facult\'e des Sciences 1,
			EDST, Equipe EDP-AN,
			Hadath, Beirut, Lebanon}
				\email{ali.wehbe@ul.edu.lb}
	 
	\date{}
	
	\begin{abstract}
	In this paper, we study the numerical stabilization of a 1D system of two wave equations coupled by velocities with an internal, local control acting on only one equation. In the theoretical part of this study \cite{ZAA1}, we distinguished two cases. 
	In the first one, the two waves assumed propagate at the same speed. Under appropriate geometric conditions, we had proved that the energy decays exponentially. While in the second case, when the waves propagate at different speeds, under appropriate geometric conditions, we had proved that the energy decays only at a polynomial rate. 
	In this paper, we confirmed these two results in a 1D numerical approximation. However, when the coupling region does not intersect the damping region, the stabilization of the system  is still theoretically an open problem. But, here in both cases, we observed an unpredicted behavior : the energy decays at an exponential rate when the propagation speeds are the same or at a polynomial rate when they are different.
\end{abstract}
	
	
 \maketitle
\tableofcontents
\section{Introduction}
In \cite{ZAA1,Wehbe-Amina-Chiraz}, the authors considered the stabilization of locally coupled wave equations. The system is described by 
\begin{equation} \label{coupledND}
	\left\{
	\begin{array}{lll}
		u_{tt}-a\Delta u+ c(x) u_t + b(x)y_t &=& 0 \hskip 1.4 cm \mbox{in} \,\,\, \Omega \times \rr_+^{*}\\
		y_{tt}- \Delta y - b(x)u_t &=&  0  \hskip 1.4 cm \mbox{in}\,\,\, \Omega \times \rr_+^{*}\\
		u = y &=& 0 \hskip 1.4 cm \mbox{on} \,\, \,  \Gamma\times \rr_+^{*}.
	\end{array}
	\right.
\end{equation}
where $\Omega$ is a nonempty connected open subset of $\rr^N$ having a boundary $\Gamma$ of class $C^2$, $a>0$ constant, $b\in C^0(\Omega, \mathbb{R})$ and $c \in C^0(\Omega, \mathbb{R}^+)$. In \cite{Wehbe-Amina-Chiraz}, the authors established an exponential energy decay rate of System \eqref{coupledND} provided that the coupling and the damping regions have non empty intersection satisfying the Piecewise Multiplier Geometric Condition (introduced in \cite{Liu1997}, and recalled in Definition 2 in \cite{ZAA1} and denoted by PMGC in short) and that the waves propagate at the same speed (i.e. $a=1$). This result generalize, in the linear case, that of \cite{Alabau2017} in the sense that the coupling coefficient function $b$ is not necessarily assumed to be positive and small enough. This result has been generalized in \cite{ZAA1} to the case when the coupling region is a subset of the damping region and satisfies a weaker geometric condition namely  Geometric Control Condition (introduced in \cite{Bardos-Lebeau-Rauch}, recalled  in Definition 1 in \cite{ZAA1} and, denoted by GCC in short). Moreover, the stabilization of System \eqref{coupledND} when the waves are not necessarily propagate at same speed (i.e. $a\not=1$ ) has been left as an open problem in \cite{Alabau2017}. 
However, in this case (i.e. $a\not=1$ ), the lack of exponential stability was proved and the optimal polynomial energy decay rate of type $\dfrac{1}{t}$ was established under different type of geometric conditions in \cite{ZAA1, Wehbe-Amina-Chiraz}. Finally, a particular and important case when the coupling region does not intersect the damping region, the stabilization of System \eqref{coupledND} is still theoretically an open problem.\\

The purpose of the present work is to focus to confirm numerically these two facts in the 1D model, where geometric conditions are automatically fulfilled, and to numerically study the case when the coupling region and the damping region does not intersect. For this sake, we firstly construct a finite difference numerical approximation of \eqref{coupledND} in a 1D model. We will construct a suitable discrete energy having the same properties of the continuous energy :
\begin{equation*}
	E(t)=\frac{1}{2} \int_\Omega \left(\abs{u_{t}}^2+a\abs{\nabla u}^2+\abs{y_{t}}^2+\abs{\nabla y}^2\right) dx. \label{energy}
\end{equation*}
This will allows us to conclude by the numerical study  the stabilization of  system \eqref{coupledND}.
\section{Finite difference scheme in one dimensional space}
This section is devoted to the numerical approximation of the problem that we considered by a finite difference discretization and to the validation of the theoretical results stated in \cite{ZAA1,Wehbe-Amina-Chiraz}.
We will firstly construct in detail a discretization in the 1D case and we will define its corresponding discrete energy. Numerical experiments are performed to validate the theoretical results. 
In fact, the numerical results in 1D show an exponential stabilization in any case when $a=1$ and a polynomial stabilization in any case in the case $a\neq 1$. 
They are better than expected.

We firstly introduce the finite difference scheme we will work on. Then we will construct the corresponding energy and finally we will perform numerical experiments.
Let us firstly recall the problem we are considered.
\medskip
\\Consider $\Omega = [0,1]$. We are interested to study the stabilization of the following coupled wave equations by velocities: 
\begin{equation} \label{coupled1D}
	\left\{
	\begin{array}{lll}
		u_{tt}-a u_{xx} + b(x)y_t  +c(x)u_t&=&  0 \quad x \in (0,1), t >0 \\
		y_{tt}- y_{xx} - b(x)u_t &=&  0   \quad x \in (0,1), t >0\\
		u(0,t)= u(1,t) = y(0,t) = y(1,t)&=& 0 \quad t > 0,
	\end{array}
	\right.
\end{equation}
with the following initial data
\begin{equation} \label{init1D0}
	u(x,0)=u_0(x), \mbox{ and } y(x,0)=y_0(x) \quad x \in (0,1)
\end{equation}
and 
\begin{equation} \label{init1D1}
	u_t(x,0)=u_1(x) \mbox{ and } y_t(x,0)=y_1(x), \quad x \in (0,1)
\end{equation}
where $a>0$ constant, $b\in C^0([0,1], \mathbb{R})$ and $c \in C^0([0,1], \mathbb{R}^+)$. We will study the two cases $a =1 $ and $a\neq 1$.
%
%
\subsection{Construction of the numerical scheme}
Let $N$ be a non negative integer. Consider the subdivision of $[0,1]$ given by 
$$0=x_0 <x_1< ...<x_N<x_{N+1}=1, \quad \mbox{ i.e. }  x_j =j \Delta x \,,\, j = 0,\ldots,N+1 \ . $$
Set $t^{n+1}-t^{n}=\Delta t$ for all $n \in \mathbb N$.
For $ j = 0,\ldots,N+1$, we denote $b_{j} = b(x_{j})$, $c_{j} = c(x_{j})$. 
The explicit finite-difference discretization of system \eqref{coupled1D} is thus, for $n \in \mathbb{N}$ and $j= 1,\ldots,N$: 
\begin{equation} \label{discrete-1D}
	\hspace*{-0.05\textwidth}	\left\{
	\begin{array}{ll}
		\dfrac{u_j^{n+1}-2u_j^n +u_j^{n-1}}{\Delta t^2} - a\dfrac{u_{j+1}^{n} -2 u_j^n +u_{j-1}^{n}}{\Delta x^2}+b_j \dfrac{y_j^{n+1}-y_j^{n-1}}{2\Delta t} + c_j \dfrac{u_j^{n+1}-u_j^{n-1}}{2 \Delta t} &= 0\\[0.5cm]
		\dfrac{y_j^{n+1}-2y_j^n +y_j^{n-1}}{\Delta t^2} - \dfrac{y_{j+1}^{n} -2 y_j^n +y_{j-1}^{n}}{\Delta x^2}-b_j \dfrac{u_j^{n+1}-u_j^{n-1}}{2\Delta t}&= 0. \\[0.5cm]
		u^{n}_{0} = u^{n}_{N+1} = 0 \\[0.5cm]
		y^{n}_{0} = y^{n}_{N+1} = 0 
	\end{array}
	\right.
\end{equation}
According to the initial conditions given by equations \eqref{init1D0}, we have firstly: for $j =1,\ldots, N$,
\begin{equation}
	u^{0}_{j} = u_{0}(x_{j})\label{u0}
\end{equation}
\begin{equation}
	y^{0}_{j} = y_{0}(x_{j})  \ .\label{y0}  
\end{equation}
We can use the second initial conditions \eqref{init1D1} to find the values of $u$ and $y$ at time $t^{1} = \Delta t$, by employing a ``ghost'' time-boundary (i.e. $t^{-1}= - \Delta t)$ and the second-order central difference formula $\mbox{for } j =1,\ldots, N$:
\begin{equation} \label{u-y-deltat}
	u_{1}(x_{j}) =\left.\dfrac{\partial u}{\partial t}\right|_{x_{j},0} =  \dfrac{ u_{j}^{1} - u_{j}^{-1} } {2 \Delta t} + O(\Delta t^{2}) .
\end{equation}
Thus we have $\mbox{for } j =1,\ldots, N$:
\begin{equation}\label{u-1}
	u_{j}^{-1}=  u_{j}^{1}  - 2 \Delta t \ u_{1}(x_{j}) \ .
\end{equation} 
We use the same discrete form of the initial conditions for $y$, $\mbox{for } j =1,\ldots, N$:
\begin{equation}\label{y-1}
	y_{j}^{-1}=  y_{j}^{1}  - 2 \Delta t \ y_{1}(x_{j}) \ .
\end{equation}  
Setting $n=0$, in the numerical scheme \eqref{discrete-1D}, the two preceding equalities permit us to compute $\left(u^{1}_{j},y^{1}_{j}\right)_{j=0,N}$. 
Finally, the solution $(u,y)$ can be computed at any time $t^n$.
\subsection{Practical implementation and CFL condition} 
Let us denote $\lambda= \dfrac{\Delta t^2}{\Delta x^2}$. We easily remark that the discrete scheme \eqref{discrete-1D} is composed of $N$ linear systems of two equations which can be written under the form:
\begin{equation}\label{system-1D2}
	\mbox{ for }  j =1,\ldots, N \,,\, M_{j} \cdot 
	\begin{pmatrix}
		u_j^{n+1} \\ \\
		y_j^{n+1}
	\end{pmatrix}
	= 
	\begin{pmatrix}
		A_j \\ \\
		B_j
	\end{pmatrix}
\end{equation}
where
$$
M_{j} = \begin{pmatrix}
	1  + \dfrac{c_j\Delta t}{2}&  \dfrac{b_j\Delta t}{2} \\ \\ 
	\dfrac{-b_j\Delta t}{2} & 1
\end{pmatrix}
$$
$$A_j= 2(1-a \lambda) u_j^n + (\dfrac{c_j}{2}\Delta t -1) u_j^{n-1} +a \lambda (u_{j+1}^n + u_{j-1}^n) + \dfrac{b_j}{2} \Delta t y_j^{n-1}$$ 
and
$$B_j = 2(1-\lambda) y_j^n + \lambda(y_{j+1}^n + y_{j-1}^n) -y_j^{n-1} -\dfrac{b_j}{2}\Delta t u_j^{n-1}.$$
Thanks to the hypothesis $\forall x \in (0,1) \,,\ c(x) \geq 0$, for $j =1,\ldots, N$ the determinant of $M_j$ given by 
$$ |M_j|=1  + \dfrac{c_j\Delta t}{2}+  \bigg(\dfrac{b_j\Delta t}{2}\bigg)^2,$$
is a strictly positive quantity.\\
Consequently, system \eqref{system-1D2} admits a unique solution  given by: for $j =1,\ldots, N$,
\begin{align}
	u_j^{n+1}=(1-a \lambda)\alpha_{j} u_j^n &+  \lambda \beta_{j} (u_{j+1}^n+u_{j-1}^n) + \gamma_{j} u_j^{n-1} - (1-\lambda)\varrho_{j} y_j^n \nonumber\\
	&- \lambda\xi_{j} (y_{j+1}^n + y_{j-1}^n) + \kappa_{j} y_j^{n-1} \label{discrete-u-1D}
\end{align}
\begin{align}
	y_j^{n+1}=(1-\lambda)\widetilde \alpha_{j} y_j^n &+ \lambda \widetilde \beta_{j} (y_{j+1}^{n} + y_{j-1}^{n})+ \widetilde \gamma_{j} y_j^{n-1}+ (1-a \lambda)\widetilde\varrho_{j} u_j^n  \nonumber \\
	&+\lambda \widetilde\xi_{j}(u_{j+1}^n+ u_{j-1}^n) + \widetilde \kappa_{j} u_j^{n-1}
	\label{discrete-y-1D} 
\end{align}
where we have set:
$$ 
\begin{array}{ll}
	\alpha_{j} = \dfrac{2}{1+\dfrac{c_j}{2}\Delta t +\left(\dfrac{b_j \Delta t}{2}\right)^{2}} \quad,
	&  \beta_{j} = \dfrac{a}{1+\dfrac{c_j}{2}\Delta t +\left(\dfrac{b_j \Delta t}{2}\right)^{2}}\quad,\\[30pt]
	\gamma_{j} = \dfrac{\dfrac{c_j}{2}\Delta t +\left(\dfrac{b_j}{2}\Delta t\right)^{2} - 1}{1+\dfrac{c_j}{2}\Delta t +\left(\dfrac{b_j \Delta t}{2}\right)^{2}}\quad,
	&\varrho_{j} = \dfrac{b_j \Delta t}{1+\dfrac{c_j}{2}\Delta t +\left(\dfrac{b_j \Delta t}{2}\right)^{2}}\quad,\\[20pt]
	\xi_{j} = \dfrac{b_j \Delta t}{2 \left(1+\dfrac{c_j}{2}\Delta t +\left(\dfrac{b_j \Delta t}{2}\right)^{2}\right)}\quad,
	&\kappa_{j} = \dfrac{b_j \Delta t}{1+\dfrac{c_j}{2}\Delta t +\left(\dfrac{b_j \Delta t}{2}\right)^{2}}\quad,
\end{array}
$$
$$ 
\begin{array}{ll}
	\widetilde \alpha_{j}= 2- \dfrac{(b_j \Delta t)^2}{2 \left(1+\dfrac{c_j}{2}\Delta t +\left(\dfrac{b_j \Delta t}{2}\right)^{2}\right)}\quad,
	&\widetilde \beta_{j} = 1 - \dfrac{(b_j \Delta t)^2}{4 \left(1+\dfrac{c_j}{2}\Delta t +\left(\dfrac{b_j \Delta t}{2}\right)^{2}\right)}\quad,\\[20pt]
	\widetilde \gamma_{j} = \dfrac{\left(b_j \Delta t\right)^2}{2 \left(1+\dfrac{c_j}{2}\Delta t +\left(\dfrac{b_j \Delta t}{2}\right)^{2}\right)} -1\quad, 
	&\widetilde \varrho_{j} = \dfrac{b_j \Delta t}{1+\dfrac{c_j}{2}\Delta t +\left(\dfrac{b_j \Delta t}{2}\right)^{2}}\quad, \\
	\widetilde \xi_{j} = \dfrac{a b_j \Delta t}{2 \left(1+\dfrac{c_j}{2}\Delta t +\left(\dfrac{b_j \Delta t}{2}\right)^{2}\right)} \quad,
	&\widetilde \kappa_{j} = \left[\dfrac{\dfrac{c_j}{2}\Delta t +\left(\dfrac{b_j \Delta t}{2}\right)^2-1} {1+\dfrac{c_j}{2}\Delta t +\left(\dfrac{b_j \Delta t}{2}\right)^{2}} -1\right] \dfrac{b_j \Delta t}{2} \quad .
\end{array}
$$
The implementation of the numerical discretization of the problem \eqref{coupled1D} consists finally of equations \eqref{u0}, \eqref{y0}, \eqref{discrete-u-1D}, \eqref{discrete-y-1D} where $(u^{-1},y^{-1})$ used for $n=0$, are defined by
\eqref{u-1}, \eqref{y-1}.

By a standard von Neumann stability analysis (that is a discrete Fourier analysis, see for instance \cite{1992-Ames}), the numerical scheme is stable if and only if, the following Courant-Friedrichs-Lewy, CFL, condition holds:
$$ \Delta t ^{2}\leq \Delta x ^{2}\mbox{ and } a \; \Delta t ^{2} \leq \Delta x ^{2}$$
which is equivalent to
\begin{equation}\label{CFL-1D}
	\Delta t \leq \min\left(1,\dfrac{1}{\sqrt{a}}\right) \Delta x \quad .
\end{equation}
The number $\min\left(1,\dfrac{1}{\sqrt{a}}\right)$ is called the CFL number and is denoted in the following by $CFL$.
\subsection{Discrete energy: definition and dissipation}
The aim of this section is to design a discrete energy that might be preserved in the case $c=0$ and to obtain the dissipation of  the discrete energy in the case $c > 0$. To this end, let us define:
\begin{itemize}
	\item the discrete kinetic energy for $u$ as: $\displaystyle  E_{k,u}^{n} = \dfrac{1}{2} \sum_{j=1}^{N}\left(\dfrac{u_j^{n+1}-u_j^n}{\Delta t}\right)^2$
	\item the discrete potential energy for $u$ as: $\displaystyle  E_{p,u}^{n} = \dfrac{a}{2} \sum_{j=0}^{N}\left(\dfrac{u_{j+1}^{n}-u_j^n}{\Delta x}\right) \left(\dfrac{u_{j+1}^{n+1}-u_j^{n+1}}{\Delta x}\right)$
	\item the discrete kinetic energy for $y$ as: $\displaystyle  E_{k,y}^{n} = \dfrac{1}{2} \sum_{j=1}^{N}\left(\dfrac{y_j^{n+1}-y_j^n}{\Delta t}\right)^2$
	\item the discrete potential energy for $u$ as: $\displaystyle  E_{y,u}^{n} = \dfrac{1}{2} \sum_{j=0}^{N}\left(\dfrac{y_{j+1}^{n}-y_j^n}{\Delta x}\right) \left(\dfrac{y_{j+1}^{n+1}-y_j^{n+1}}{\Delta x}\right)$
\end{itemize}
The total discrete energy is then defined as
\begin{equation}\label{discrete-energy}
	\mathcal{E}^{n} = E_{k,u}^{n} + E_{p,u}^{n} + E_{k,y}^{n} + E_{p,u}^{n}.
\end{equation}

Let us prove now that this definition of the energy fulfills the two properties stated above. For this sake, we multiply the first equation of \eqref{discrete-1D} by $(u_j^{n+1}-u_j^{n-1})$ and we sum over $j = 1,\ldots,N$. We obtain: 
\begin{align}
	&\displaystyle\sum_{j=1}^{N}\dfrac{u_j^{n+1}-2u_j^n + u_j^{n-1}}{\Delta t^2}(u_j^{n+1}-u_j^{n-1})
	- a\sum_{j=1}^{N} \dfrac{u_{j+1}^{n} - 2 u_j^n + u_{j-1}^{n}}{\Delta x^2} (u_j^{n+1}-u_j^{n-1}) \nonumber\\ 
	& \displaystyle+ \sum_{j=1}^{N}b_j \dfrac{y_j^{n+1}-y_j^{n-1}}{2\Delta t}(u_j^{n+1}-u_j^{n-1}) + \sum_{j=1}^{N} c_j \dfrac{(u_j^{n+1}-u_j^{n-1})^2}{2 \Delta t}=0. \label{Discrete.energy-1D}
\end{align}
\textbf{Estimation of the first term of \eqref{Discrete.energy-1D}} We firstly have:
\begin{align}
	\displaystyle\sum_{j=1}^{N}\dfrac{u_j^{n+1}-2u_j^n + u_j^{n-1}}{\Delta t^2}(u_j^{n+1}-u_j^{n-1}) &= \displaystyle\sum_{j=1}^{N}\dfrac{u_j^{n+1}-u_j^n-(u_j^n -u_j^{n-1})}{\Delta t^2}(u_j^{n+1}-u_j^n+u_j^n-u_j^{n-1})\nonumber\\
	&= \sum_{j=1}^{N} \bigg(\dfrac{u_j^{n+1}-u_j^n}{\Delta t}\bigg)^2 -\sum_{j=1}^{N} \bigg(\dfrac{u_j^{n+1}-u_j^{n-1}}{\Delta t}\bigg)^{2} \nonumber\\
	& = 2 (E_{k,u}^{n} - E_{k,u}^{n-1}). \label{discr.EC.u}
\end{align}
\textbf{Estimation of the second term of \eqref{Discrete.energy-1D}.} Using the same trick we have: 
\begin{align*}
	- a \ \sum_{j=1}^{N} \dfrac{u_{j+1}^{n} - 2 u_j^n + u_{j-1}^{n}}{\Delta x^2} (u_j^{n+1}-u_j^{n-1}) &= - a \  \sum_{j=1}^{N} \dfrac{u_{j+1}^{n} - u_j^n -( u_j^n-u_{j-1}^{n})}{\Delta x^2} (u_j^{n+1}-u_j^{n-1}) \\
	&= - a \ \sum_{j=1}^{N}\dfrac{(u_{j+1}^n-u_j^n)(u_j^{n+1}-u_j^{n-1})}{\Delta x^2} \\
	&+ a \ \sum_{j=1}^{N+1}\dfrac{(u_{j}^n-u_{j-1}^n)(u_j^{n+1}-u_j^{n-1})}{\Delta x^2}.
\end{align*}
So, by translation of index in the second term in the previous sum, we will have: 
\begin{align}
	- a \ \sum_{j=1}^{N} \dfrac{u_{j+1}^{n} - 2 u_j^n + u_{j-1}^{n}}{\Delta x^2} (u_j^{n+1}-u_j^{n-1}) &= - a \  \sum_{j=0}^{N}\dfrac{(u_{j+1}^n-u_j^n)(u_j^{n+1}-u_j^{n-1})}{\Delta x^2}  \nonumber\\
	&+ a \ \sum_{j=0}^{N}\dfrac{(u_{j+1}^n-u_{j}^n)(u_{j+1}^{n+1}-u_{j+1}^{n-1})}{\Delta x^2} \nonumber
\end{align}
\begin{align}
	- a \ \sum_{j=1}^{N} \dfrac{u_{j+1}^{n} - 2 u_j^n + u_{j-1}^{n}}{\Delta x^2} (u_j^{n+1}-u_j^{n-1}) &=   a \ \sum_{j=0}^{N}\dfrac{(u_{j+1}^{n+1}-u_j^{n+1})(u_{j+1}^{n}-u_j^{n})}{\Delta x^2}  \nonumber\\
	&-a \ \sum_{j=0}^{N}\dfrac{(u_{j+1}^{n-1}-u_{j}^{n-1})(u_{j+1}^{n}-u_{j}^{n})}{\Delta x^2} \nonumber\\
	&=2 (E_{p,u}^{n} - E_{p,u}^{n-1}).\label{discr.Ep.u}
\end{align}
Substituting \eqref{discr.EC.u} and \eqref{discr.Ep.u} into \eqref{Discrete.energy-1D}, we get \\ 
\begin{align} \label{Energy.u.1D}
	2 \left(E_{k,u}^{n} + E_{p,u}^{n} - E_{k,u}^{n-1} - E_{p,u}^{n-1} \right)&+ 2 \Delta t \sum_{j=1}^{N} c_j \left(\dfrac{u_j^{n+1}-u_j^{n-1}}{2 \Delta t}\right)^{2} \nonumber\\
	&+ \sum_{j=1}^{N}b_j \dfrac{y_j^{n+1}-y_j^{n-1}}{2\Delta t}(u_j^{n+1}-u_j^{n-1}) = 0.
\end{align}
Similarly, by multiplying the second equation of \eqref{discrete-1D} by $(y_j^{n+1} - y_j^{n-1})$, and using the same algebraic tricks, we will get:
\begin{equation}\label{Energy.y.1D}
	2 \left(E_{k,y}^{n} + E_{p,y}^{n} - E_{k,y}^{n-1} - E_{p,y}^{n-1} \right) - \sum_{j=1}^{N}b_j \dfrac{u_j^{n+1}-u_j^{n-1}}{2\Delta t}(y_j^{n+1}-y_j^{n-1})= 0.
\end{equation}
Using the definition of the total discrete energy, \eqref{discrete-energy}, and the two equations \eqref{Energy.u.1D},  \eqref{Energy.y.1D} leads to:
\begin{equation}\label{discrete-dissipation}
	\left(\mathcal{E}^{n} - \mathcal{E}^{n-1}\right) + \Delta t  \sum_{j=1}^{N} c_j \left(\dfrac{u_j^{n+1}-u_j^{n-1}}{2 \Delta t}\right)^{2}  = 0.
\end{equation}
Consequently, the total discrete energy of system \eqref{discrete-1D} is decreasing along time.
\section{Numerical experiments: validation of the theoretical results}
In every experiment, we have chosen: 
\[ u_{0}(x) = x (x-1) \,,\, u_{1}(x) = x (x-1) \,,\, y_{0}(x) = - x (x-1) \,,\, y_{1}(x) = - x (x-1).
\]

The mesh size is chosen as $N = 100$ so that $\Delta x = 0.01$ and the time step
is chosen as $\dfrac{\Delta t}{\Delta x} = CFL$.

\bigskip
In order to validate the different theoretical results, we have chosen different functions $b$ and $c$ synthesized in the list below:
\begin{itemize}
	\item No coupling: $ b_{1}(x) = 0$ or no dissipation $c_{1}(x)=0$,
	\item Full coupling $\displaystyle b_{2}(x) = \mathds{1}_{(0,1)}(x)$ or full  dissipation $c_{2}(x)= \mathds{1}_{(0,1)}(x)$,
	\item Partial coupling $\displaystyle b_{3}(x) = \mathds{1}_{[0.1,0.2]\cup[0.8,0.9]}(x)$  or partial dissipation $c_{3}(x)= \mathds{1}_{[0.1,0.2]\cup[0.8,0.9]}(x)$,
	\item Partial coupling $\displaystyle b_{4}(x) = \mathds{1}_{[0.1,0.2]}(x)$  or partial dissipation $ c_{4}(x)= \mathds{1}_{[0.1,0.2]}(x)$,
	\item Partial coupling $\displaystyle b_{5}(x) = \mathds{1}_{[0.4,0.6]}(x)$  or partial dissipation $ c_{5}(x)= \mathds{1}_{[0.4,0.6]}(x)$.
\end{itemize}
Combining the different choices of the coupling and damping functions in order to have or not $\omega_{b} \cap \omega_{c_{+}} \neq \emptyset$  will permit us to validate the theoretical results.
\medskip
\\Let us notice that in the special case of the dimension 1, the geometric control condition GCC holds as soon as $\omega_{c_{+}} \neq \emptyset$.
\subsection{Same propagation speed: $\boldsymbol{a=1}$}
For every numerical simulation, the final time $T$ is chosen as $T = 500$.
\bigskip
\subsubsection{{No damping: conservation of the total energy}}
\begin{figure}[H]
	\centering
	\includegraphics[width=0.5\textwidth]{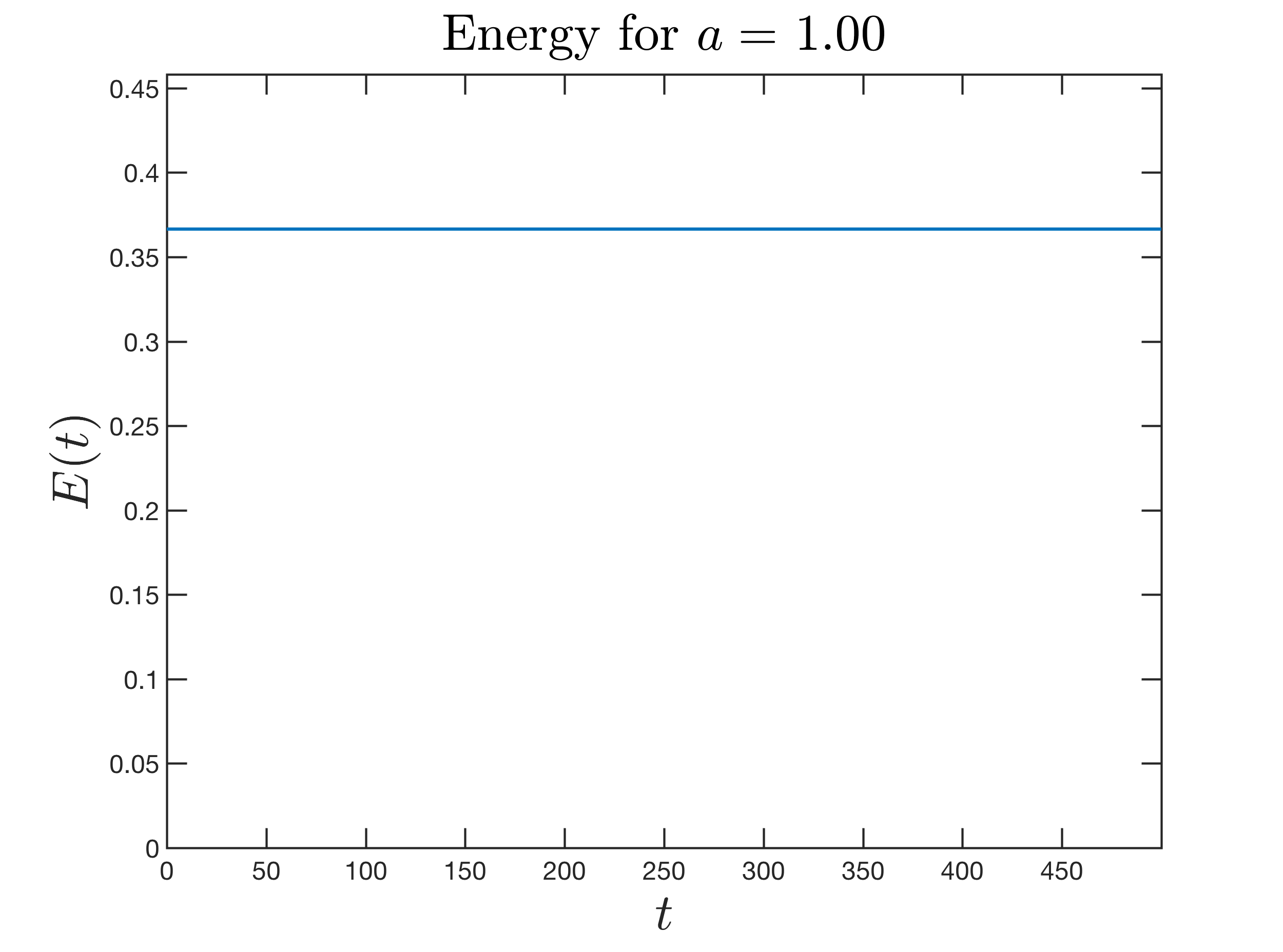}
	\captionsetup{justification=centering}
	\caption{No damping.\\$c = c_{1} = 0$, partial coupling $b = b_{3} = \mathds{1}_{[0.1,0.2]\cup[0.8,0.9]}(x)$}
	\label{b3-c1}
\end{figure}

Firstly, let us verify that when no damping are present, the discrete energy is conserved. 
We present in figure \ref{b3-c1} the numerical experiment when $c = c_{1} = 0$ and $b = b_{3} =  \mathds{1}_{[0.1,0.2]\cup[0.8,0.9]}(x)$. 
Indeed, the total energy is conserved along time.

\begin{remark}
	This numerical test where no damping is applied  shows that without a damping term, the total energy is completely conserved. This fact suggests that the numerical scheme
	does not produce numerical dissipation. So the numerical behavior observed thereafter is only due to the considered model.
\end{remark}
\bigskip
\subsubsection{{$\boldsymbol{\omega_{b} \cap \omega_{c_{+}} \neq \emptyset}$. Exponential stability}}

\begin{figure}[H]
	\centering
	\subcaptionbox{Energy.\label{b4-c3-Energy}}
	{\includegraphics[width=0.45\textwidth]{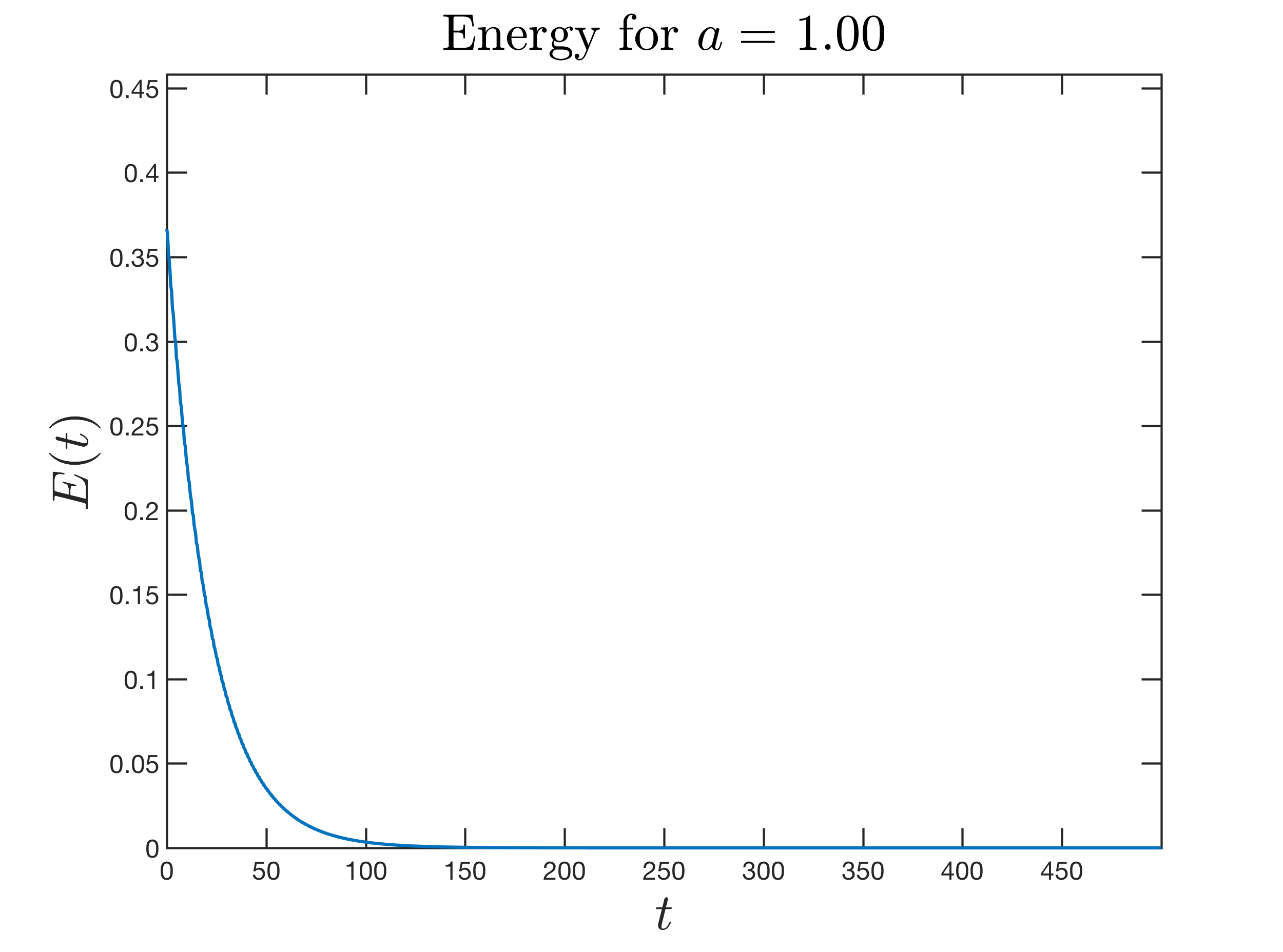}}
	\subcaptionbox{Exponential decay.\label{b4-c3-Exponential}}
	{\includegraphics[width=0.45\textwidth]{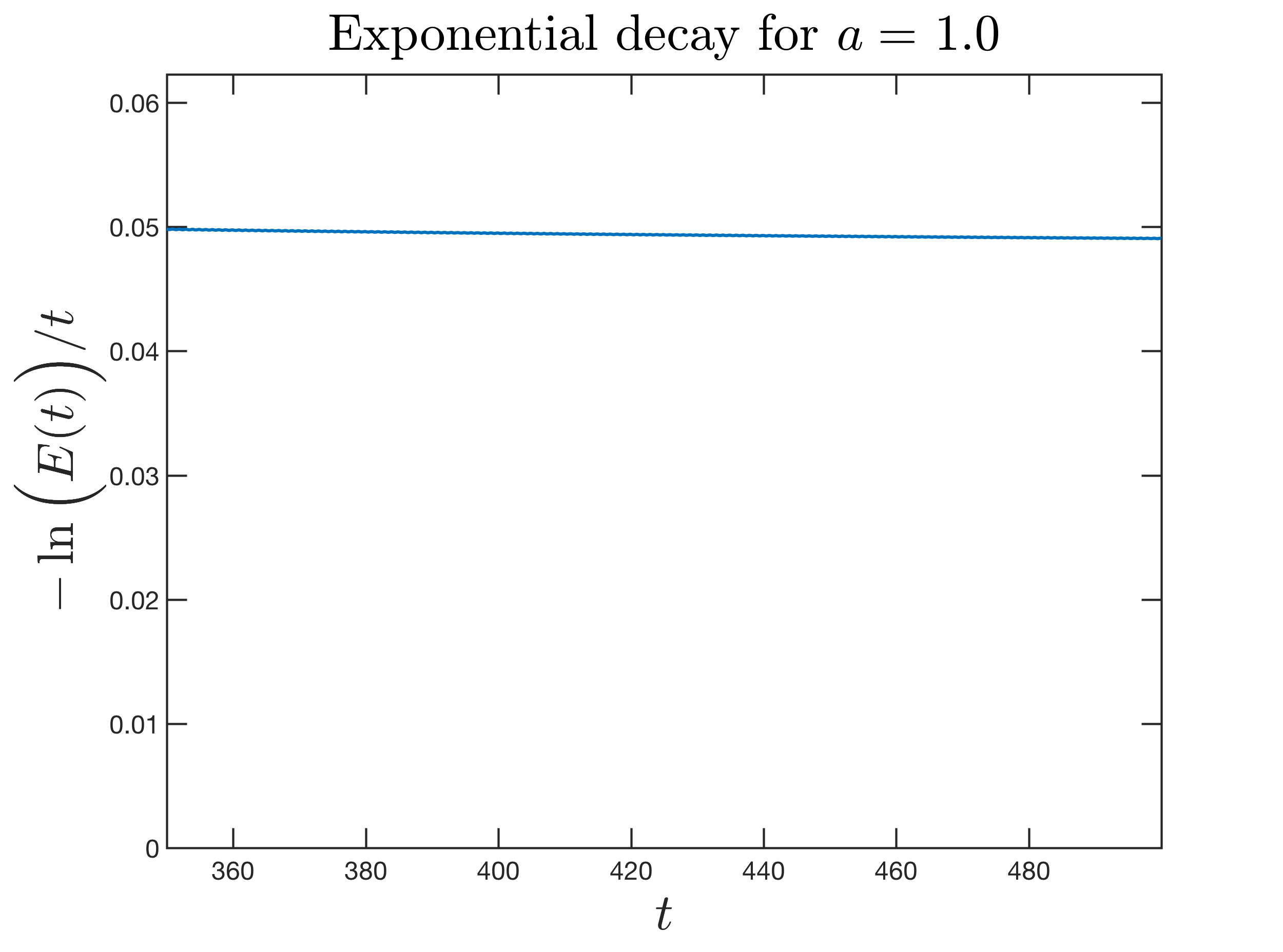}}\\[0.01\textheight]
	\subcaptionbox{Final time profile.\label{b4-c3-final}}
	{\includegraphics[width=0.45\textwidth]{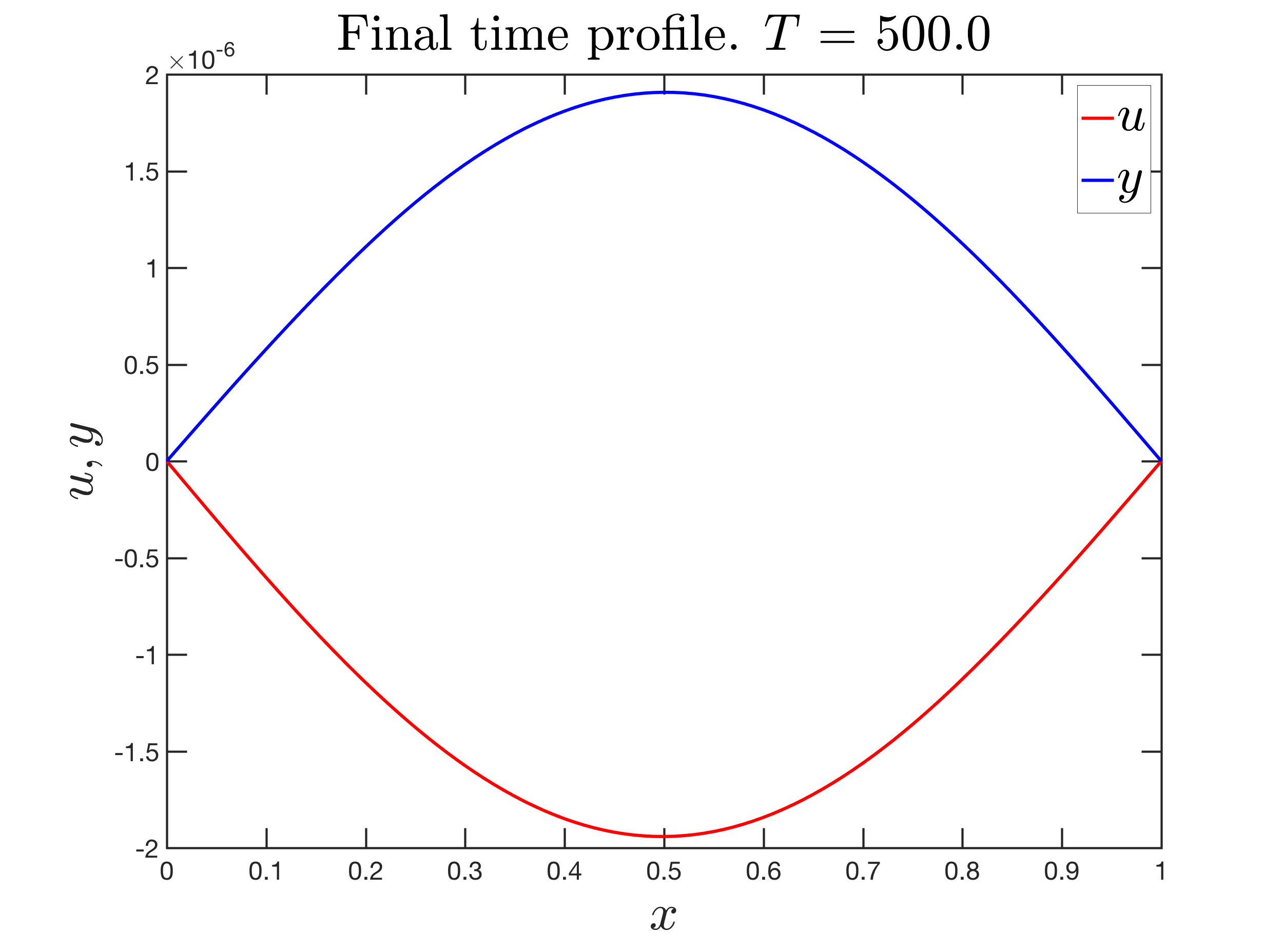}}
	\captionsetup{justification=centering}
	\caption{Long time behavior when $\omega_{b} \cap \omega_{c_{+}} \neq \emptyset$.\\
		\footnotesize{$b = b_{4}(x) = \mathds{1}_{[0.1,0.2]}(x)$ and $c = c_{3}(x)= \mathds{1}_{[0.1,0.2]\cup[0.8,0.9]}(x)$}}
	\label{b4-c3}
\end{figure}

Let us now verify the theoretical results when we suppose that  $\omega_{b} \cap \omega_{c_{+}} \neq \emptyset$.
For this sake, we present in figure \ref{b4-c3}, the total energy and the quantity $-\ln\big(E(t)\big) /t$ versus time $t$ for large time,
where we have chosen $b = b_{4}(x) = \mathds{1}_{[0.1,0.2]}(x)$ and $c = c_{3}(x)= \mathds{1}_{[0.1,0.2]\cup[0.8,0.9]}(x)$.
This choice verifies the assumption that $\omega_{b} \cap \omega_{c_{+}} \neq \emptyset$ and in figure \ref{b4-c3}, it is shown that the energy is decreasing 
and an exponential decay is observed since it seems that $-\ln\big(E(t)\big) /t$ tends to a constant as $t \rightarrow + \infty$.
The final time profile confirms that $u$ and $y$ are small and the final profiles of $u$ and $y$ are smooth as expected (high frequency oscillations are exponentially dissipated).


\bigskip
\subsubsection{$\boldsymbol{\omega_{b} \cap \omega_{c_{+}} = \emptyset}$. Unpredicted behavior}
\begin{figure}[H]
	\centering
	\subcaptionbox{Energy.\label{b4-c5-Energy}}
	{\includegraphics[width=0.45\textwidth]{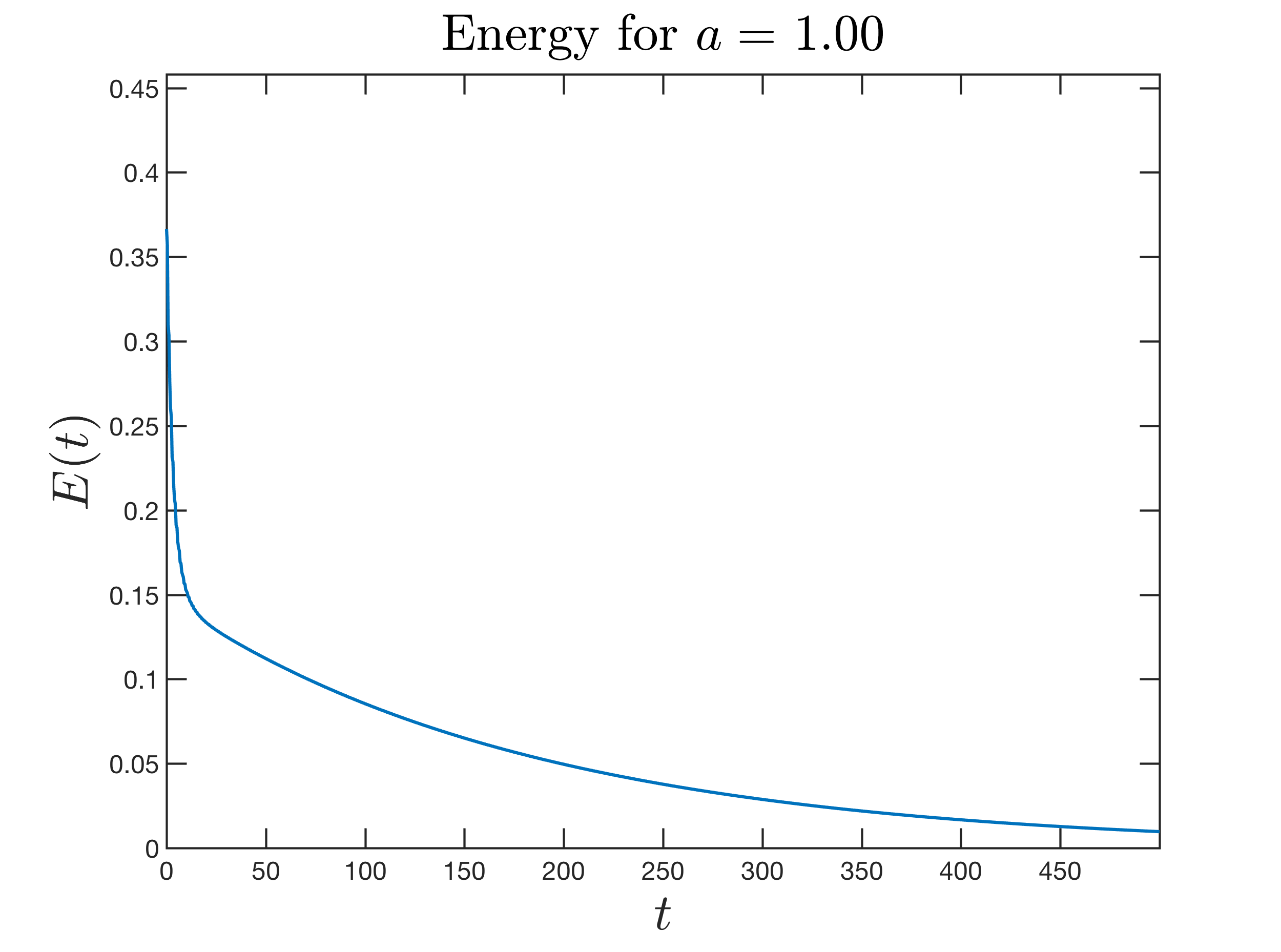}}
	\subcaptionbox{Exponential decay.\label{b4-c5-Exponential}}
	{\includegraphics[width=0.45\textwidth]{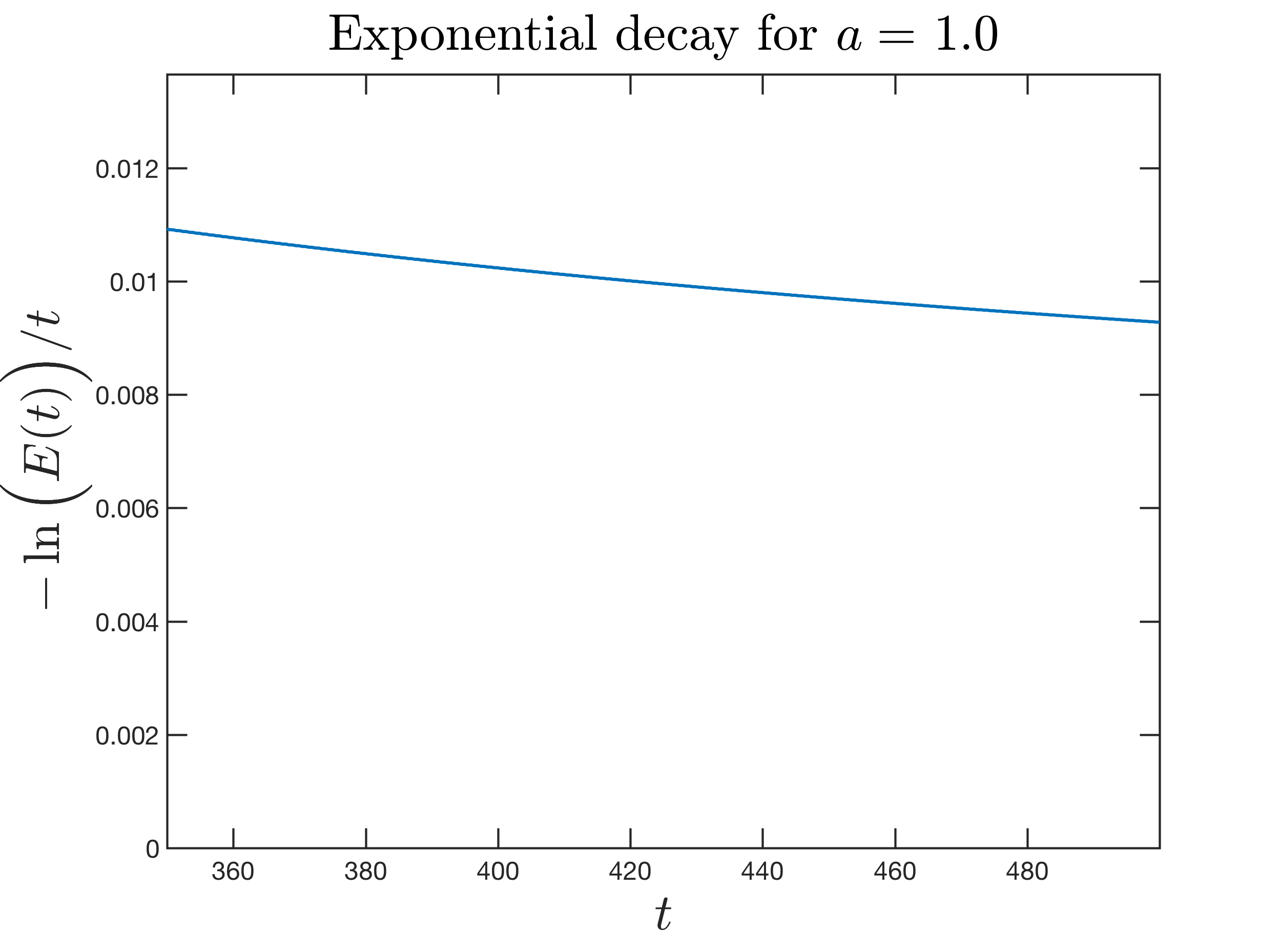}}\\[0.01\textheight]
	\subcaptionbox{Final time profile.\label{b4-c5-final}}
	{\includegraphics[width=0.45\textwidth]{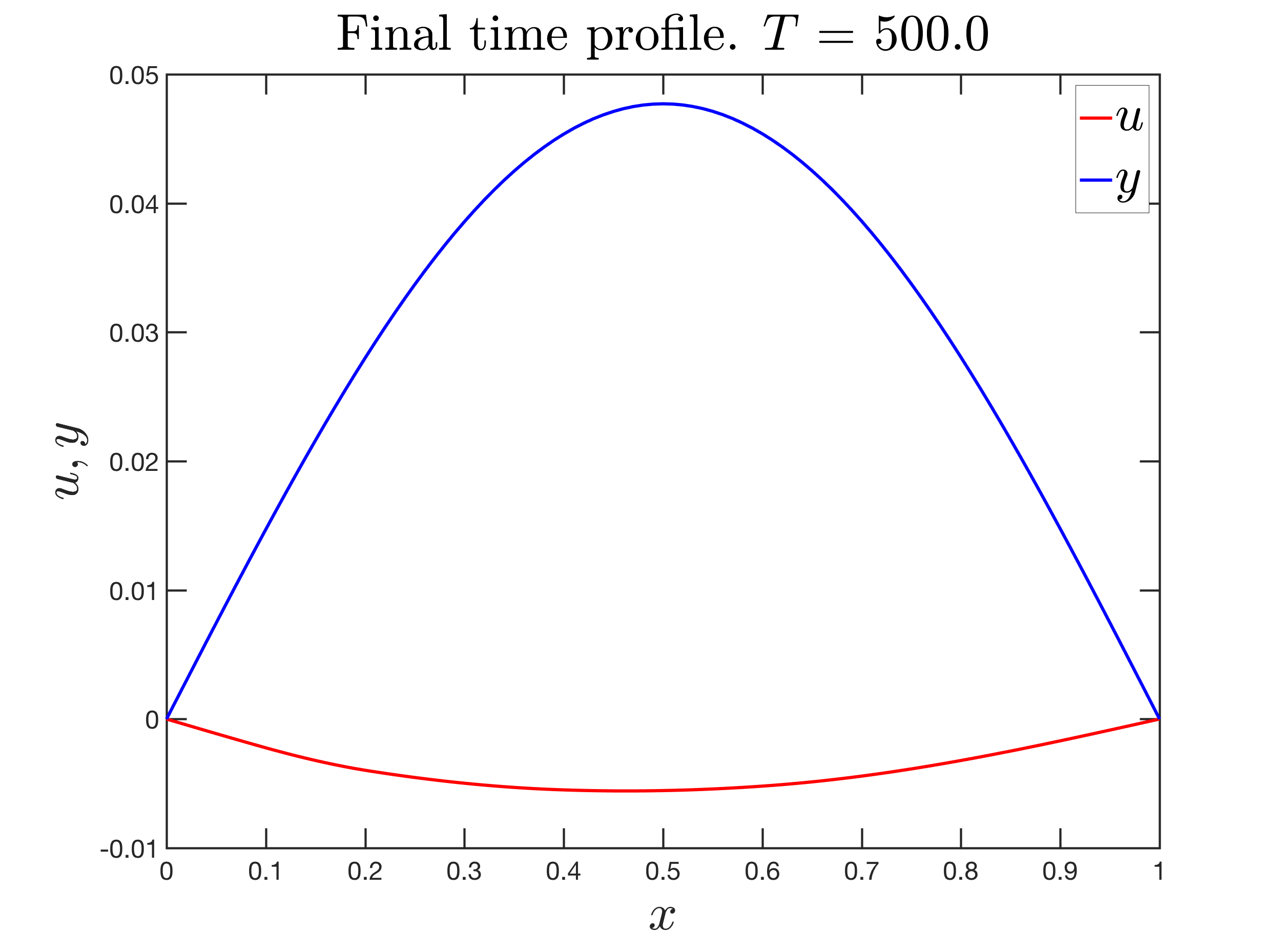}}
	\captionsetup{justification=centering}
	\caption{Long time behavior when $\omega_{b} \cap \omega_{c_{+}} = \emptyset$.\\
		\footnotesize{$b = b_{4}(x) = \mathds{1}_{[0.1,0.2]}(x)$ and $c = c_{5}(x)= \mathds{1}_{[0.4,0.6]}(x)$}}
	\label{b4-c5}
\end{figure}

At the numerical level, we are interested in the long time behavior of the solution $(u,y)$ when we suppose that  $\omega_{b} \cap \omega_{c_{+}} = \emptyset$.
For this sake, we present in figure \ref{b4-c5}, the total energy and the quantity $-\ln\big(E(t)\big) /t$ versus time $t$ for large time,
where we have chosen $b = b_{4}(x) = \mathds{1}_{[0.1,0.2]}(x)$ and $c = c_{5}(x)= \mathds{1}_{[0.4,0.6]}(x)$.
This choice verifies the assumption that $\omega_{b} \cap \omega_{c_{+}} = \emptyset$. In figure \ref{b4-c5}, it is shown that the energy is decreasing 
and an exponential decay is observed since it seems that $-\ln\big(E(t)\big) /t$ tends to a constant as $t \rightarrow + \infty$. 
The final time profile confirms that $u$ and $y$ are small and again the couple of solution $(u,y)$ is smooth .
We have not considered this case in the theoretical study and this numerical result shows a similar behavior as in the case presented before.

\begin{figure}[H]
	\centering
	\subcaptionbox{Energy.\label{b5-c4-Energy}}
	{\includegraphics[width=0.45\textwidth]{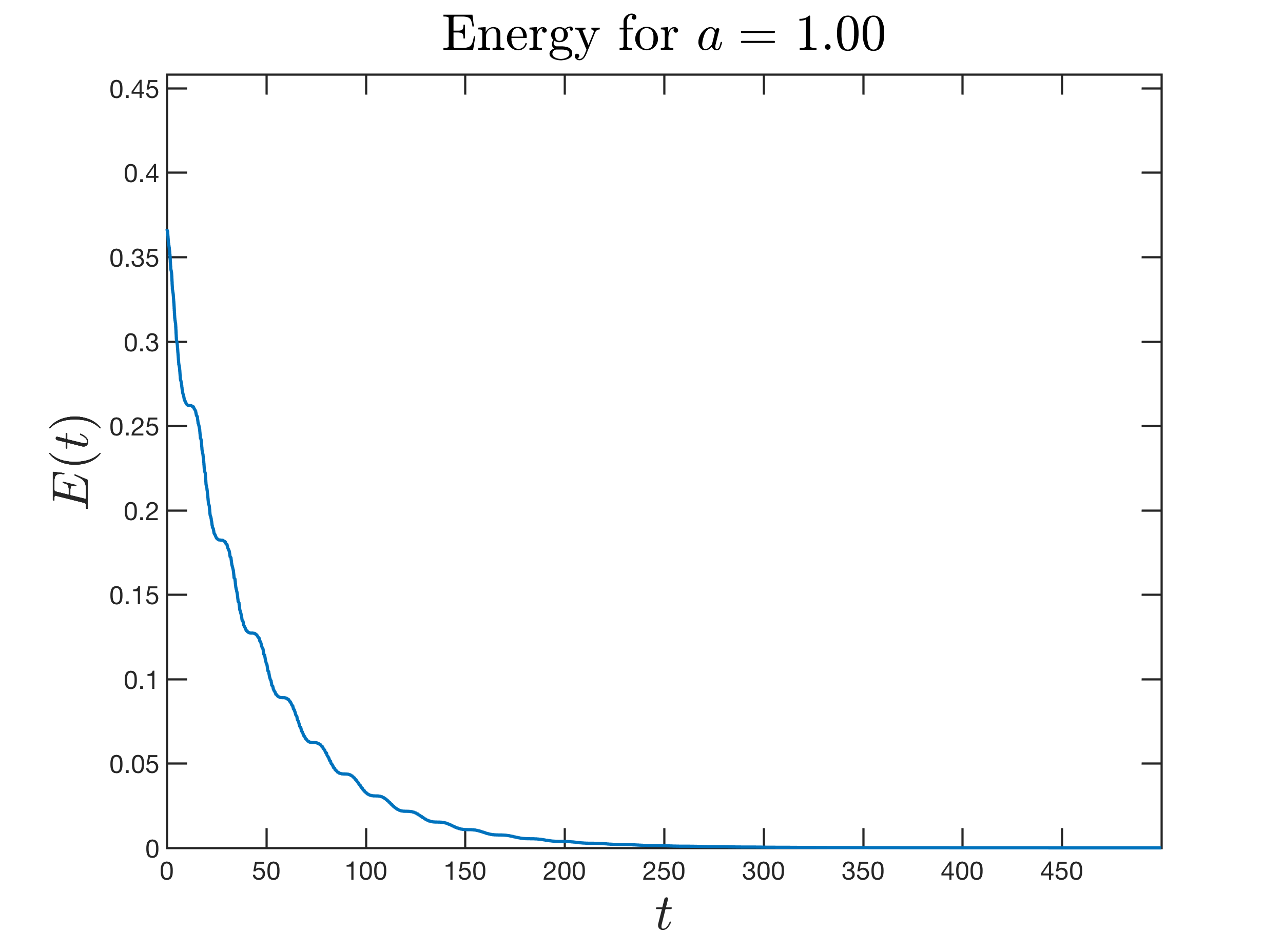}}
	\subcaptionbox{Exponential decay.\label{b5-c4-Exponential}}
	{\includegraphics[width=0.45\textwidth]{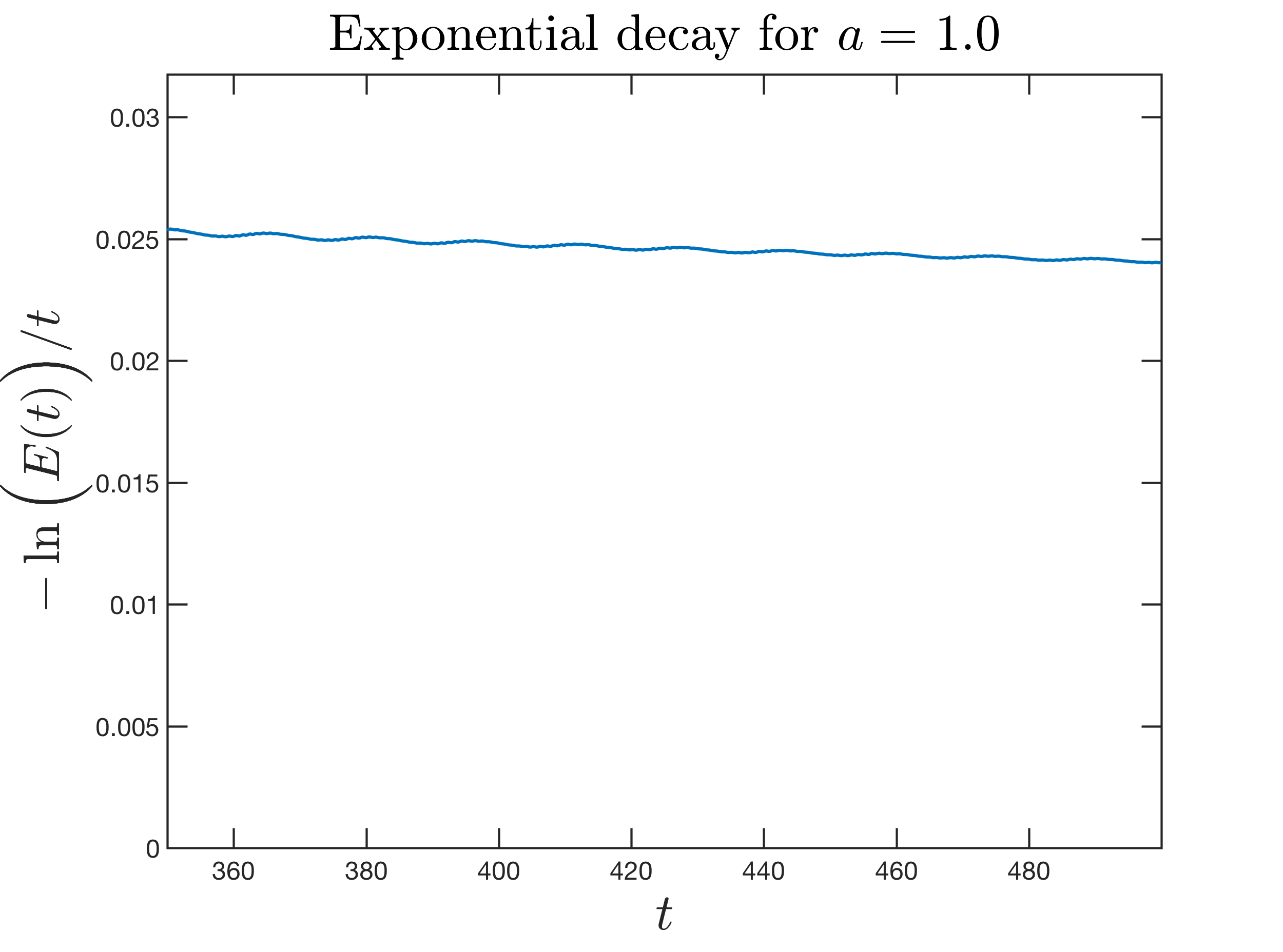}}\\[0.01\textheight]
	\subcaptionbox{Final time profile.\label{b5-c4-final}}
	{\includegraphics[width=0.45\textwidth]{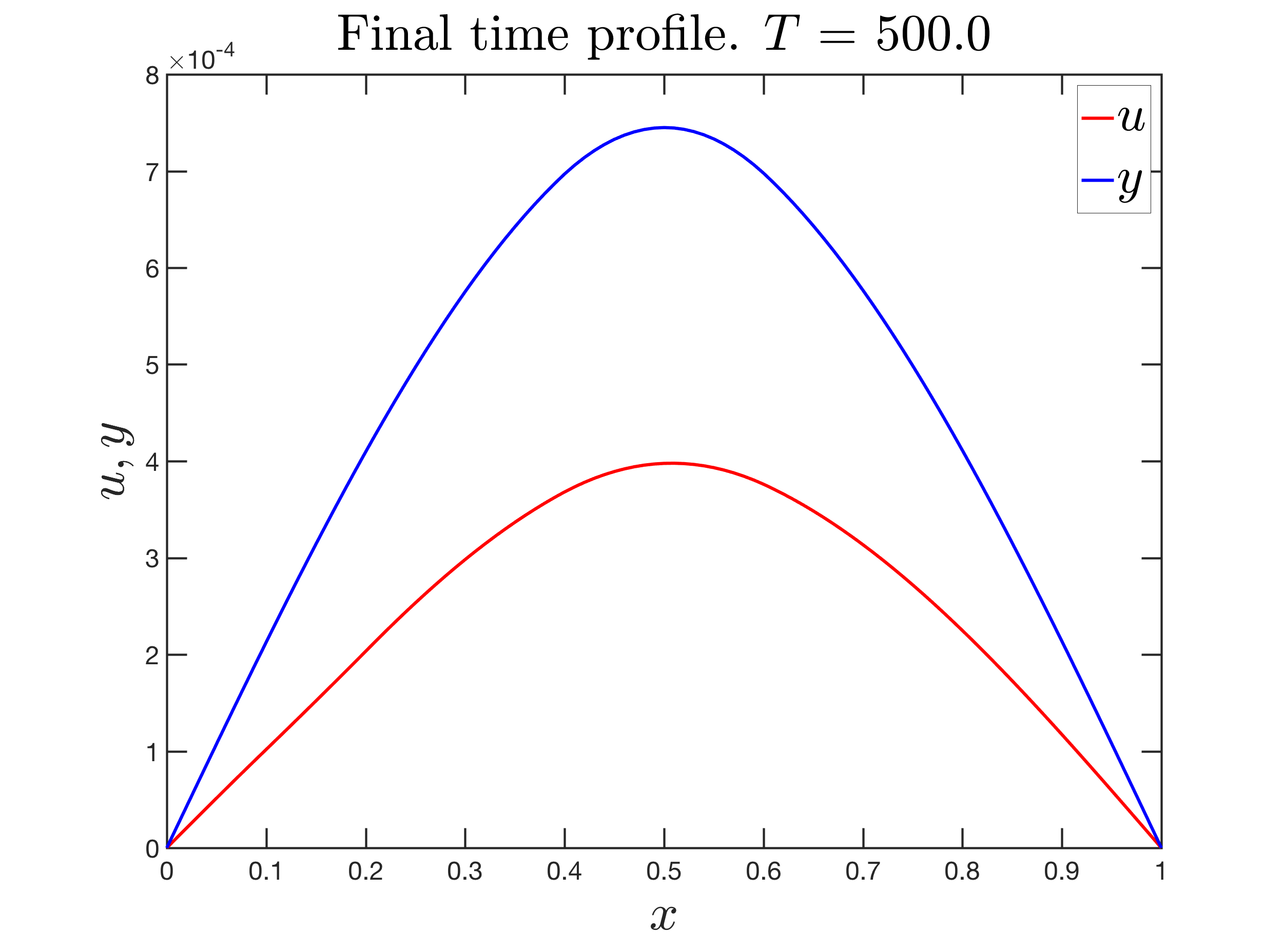}}
	\captionsetup{justification=centering}
	\caption{Long time behavior when $\omega_{b} \cap \omega_{c_{+}} = \emptyset$.\\
		\footnotesize{$b = b_{5}(x) = \mathds{1}_{[0.4,0.6]}(x)$ and $c = c_{4}(x)= \mathds{1}_{[0.1,0.2]}(x)$}}
	\label{b5-c4}
\end{figure}

So we decided to confirm this behavior by choosing $b = b_{5}(x) = \mathds{1}_{[0.4,0.6]}(x)$ and \linebreak
$c=c_{4}(x)=\mathds{1}_{[0.1,0.2]}(x)$.
This choice verifies also the assumption that $\omega_{b} \cap \omega_{c_{+}} = \emptyset$.
In figure \ref{b5-c4}, it is shown that the energy is decreasing 
and an exponential decay is observed since it seems that $-\ln\big(E(t)\big) /t$ tends to a constant as $t \rightarrow + \infty$. 
The final time profile confirms that $u$ and $y$ are small and again the couple of solution $(u,y)$ is smooth.

\begin{remark}
	Let us notice that when the propagation speeds are the same for $u$ and $y$, the final profiles of the solution $u\,,\,y$ presented in figure \ref{b4-c3-final}, figure \ref{b4-c5-final} and  in figure \ref{b5-c4-final} have the same form as the initial one, that is no spurious oscillations due to high frequency  are present.
\end{remark}
\subsection{Different propagation speed: $\boldsymbol{a > 1}$}
We investigate now the long time behavior of $(u,y)$ when the propagation speeds are different and specifically when $a> 1$. So we have chosen to take $a=2$. We firstly investigate the case when the propagation speed for $u$
is greater than the one of $y$ namely $a > 1$. We have chosen $a = 2$.
\bigskip
\subsubsection{{$\boldsymbol{\omega_{b} \cap \omega_{c_{+}} \neq \emptyset}$. Polynomial stability}}
\begin{figure}[H]
	\centering
	\subcaptionbox{Final time $T = 500$.\label{Energy-small-b4-c3-a2}}
	{\includegraphics[width=0.4\textwidth]{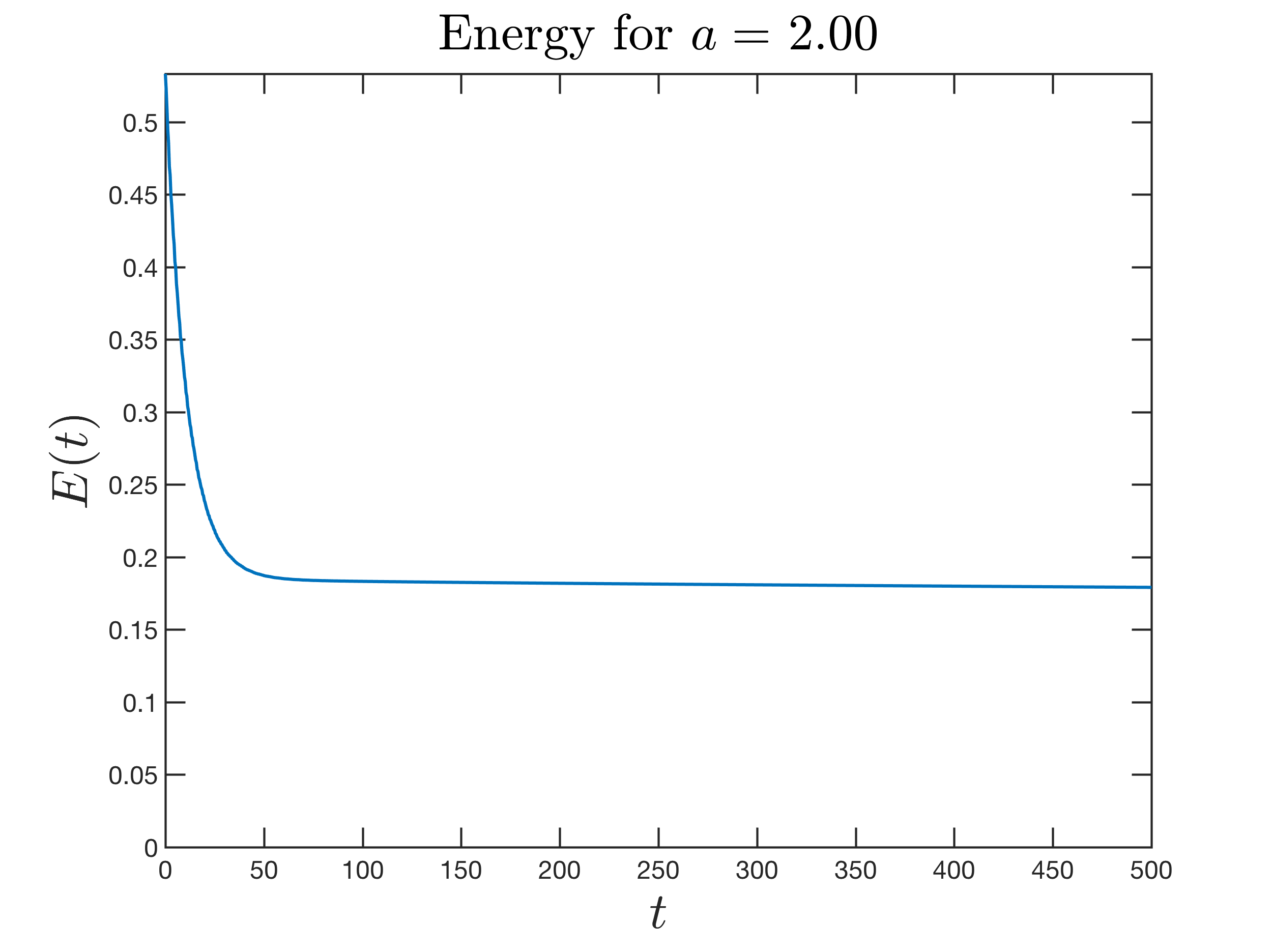}}
	\subcaptionbox{Final time $T= 500~000$.\label{Energy-b4-c3-a2}}
	{\includegraphics[width=0.4\textwidth]{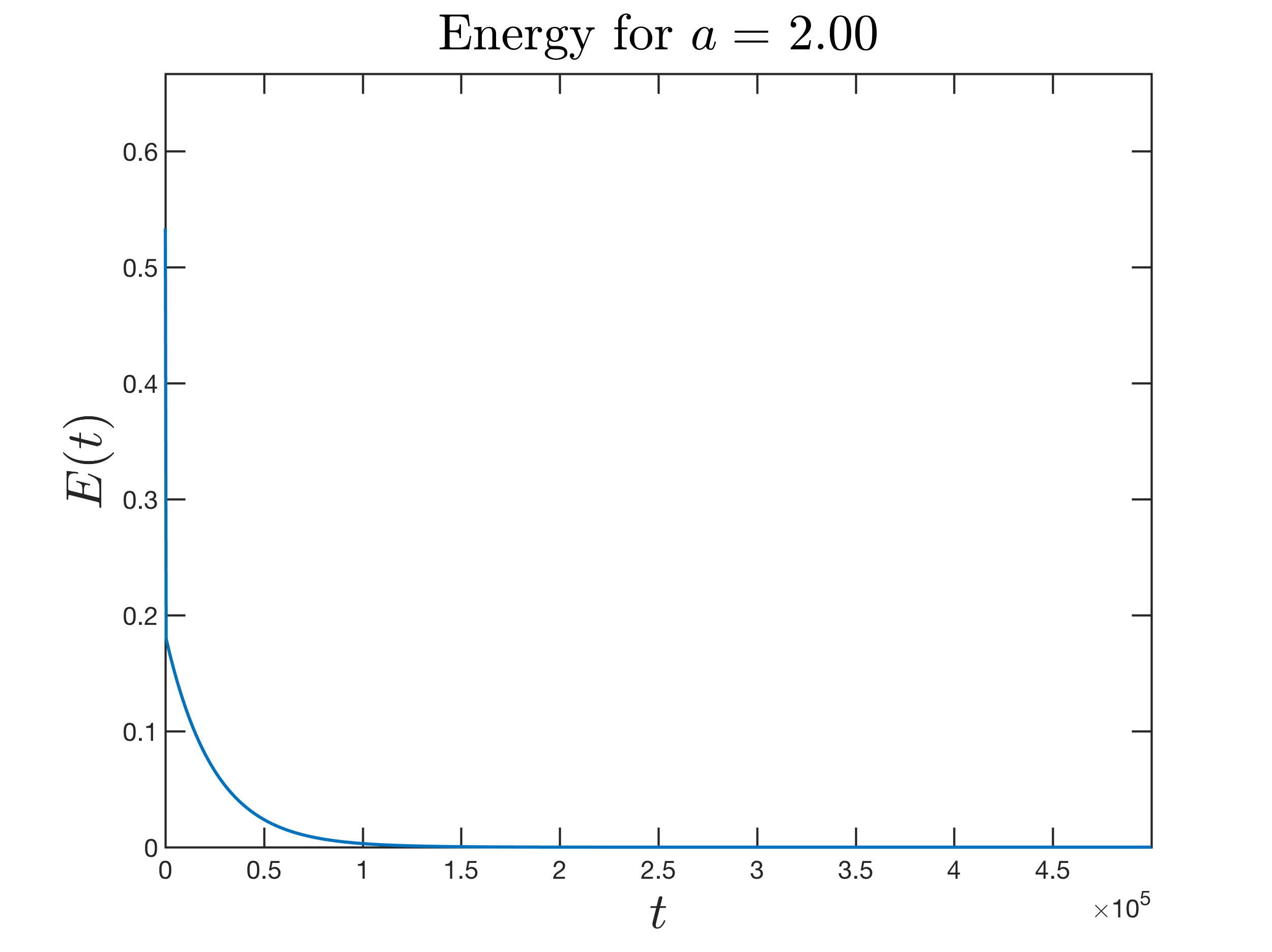}}\\[0.01\textheight]
	\captionsetup{justification=centering}
	\caption{Energy when $\omega_{b} \cap \omega_{c_{+}} \neq \emptyset$.\\
		\footnotesize{$b = b_{4}(x) = \mathds{1}_{[0.1,0.2]}(x)$ and $c = c_{3}(x)= \mathds{1}_{[0.1,0.2]\cup[0.8,0.9]}(x)$}}
	\label{Ener-b4-c3-a2}
\end{figure}

Let us now verify the theoretical results when we suppose that  $\omega_{b} \cap \omega_{c_{+}} \neq \emptyset$.
For this sake, we present in figure \ref{Ener-b4-c3-a2}, the total energy 
where we have chosen $b = b_{4}(x) = \mathds{1}_{[0.1,0.2]}(x)$ and $c = c_{3}(x)= \mathds{1}_{[0.1,0.2]\cup[0.8,0.9]}(x)$.

When taking as final time $T = 500$, it seems that the energy does not tend to zero as shown in figure \ref{Energy-small-b4-c3-a2}.
This is the reason why we have chosen for the case when $a\neq 1$  as final time $T = 500~000$ and figure \ref{Energy-b4-c3-a2}
shows that the energy finally goes to zero.

\begin{figure}[H]
	\centering
	\subcaptionbox{Exponential decay?\label{Exp-b4-c3-a2}}
	{\includegraphics[width=0.4\textwidth]{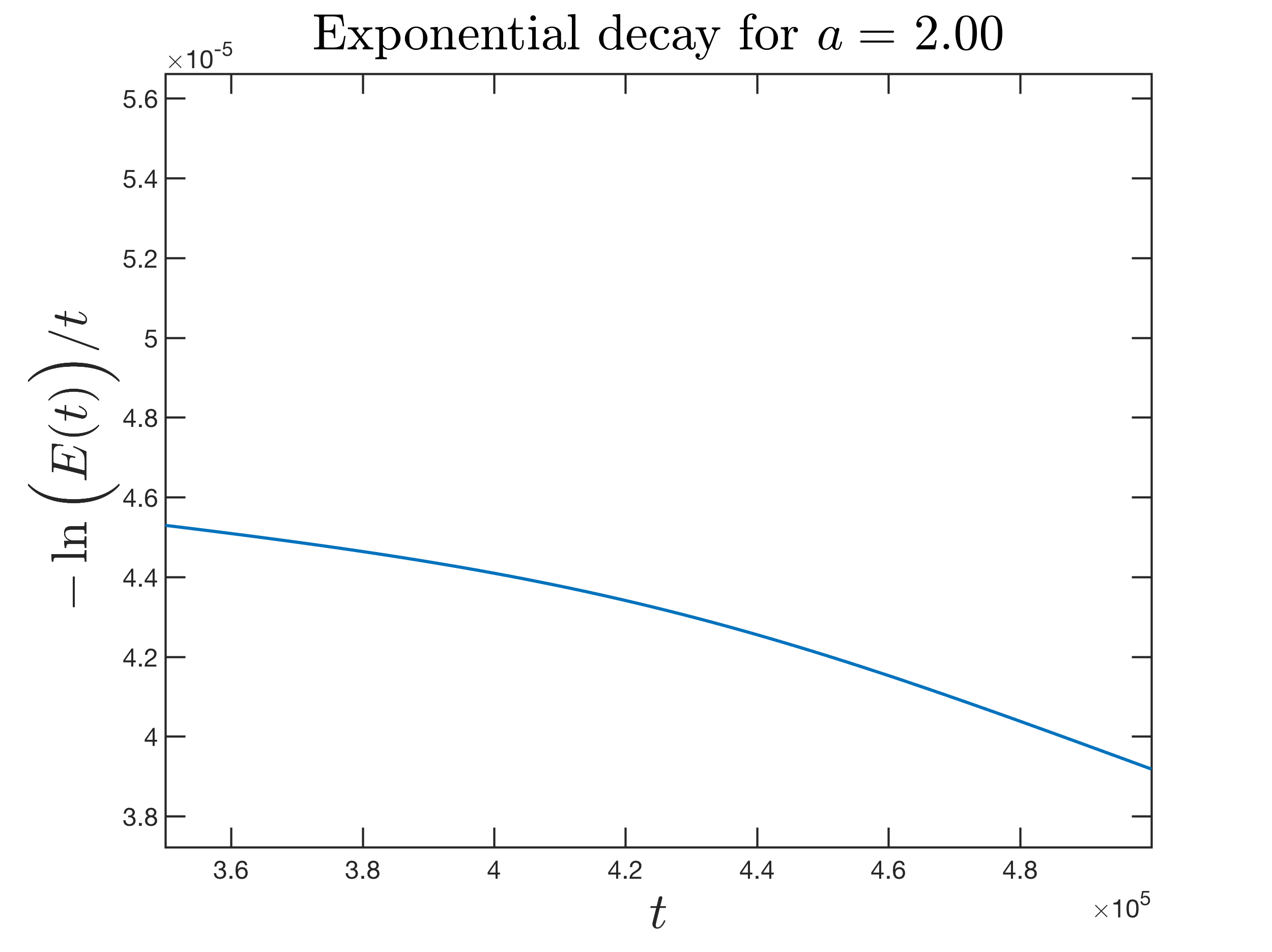}}
	\subcaptionbox{Polynomial decay in $1/t$?\label{1-b4-c3-a2}}
	{\includegraphics[width=0.4\textwidth]{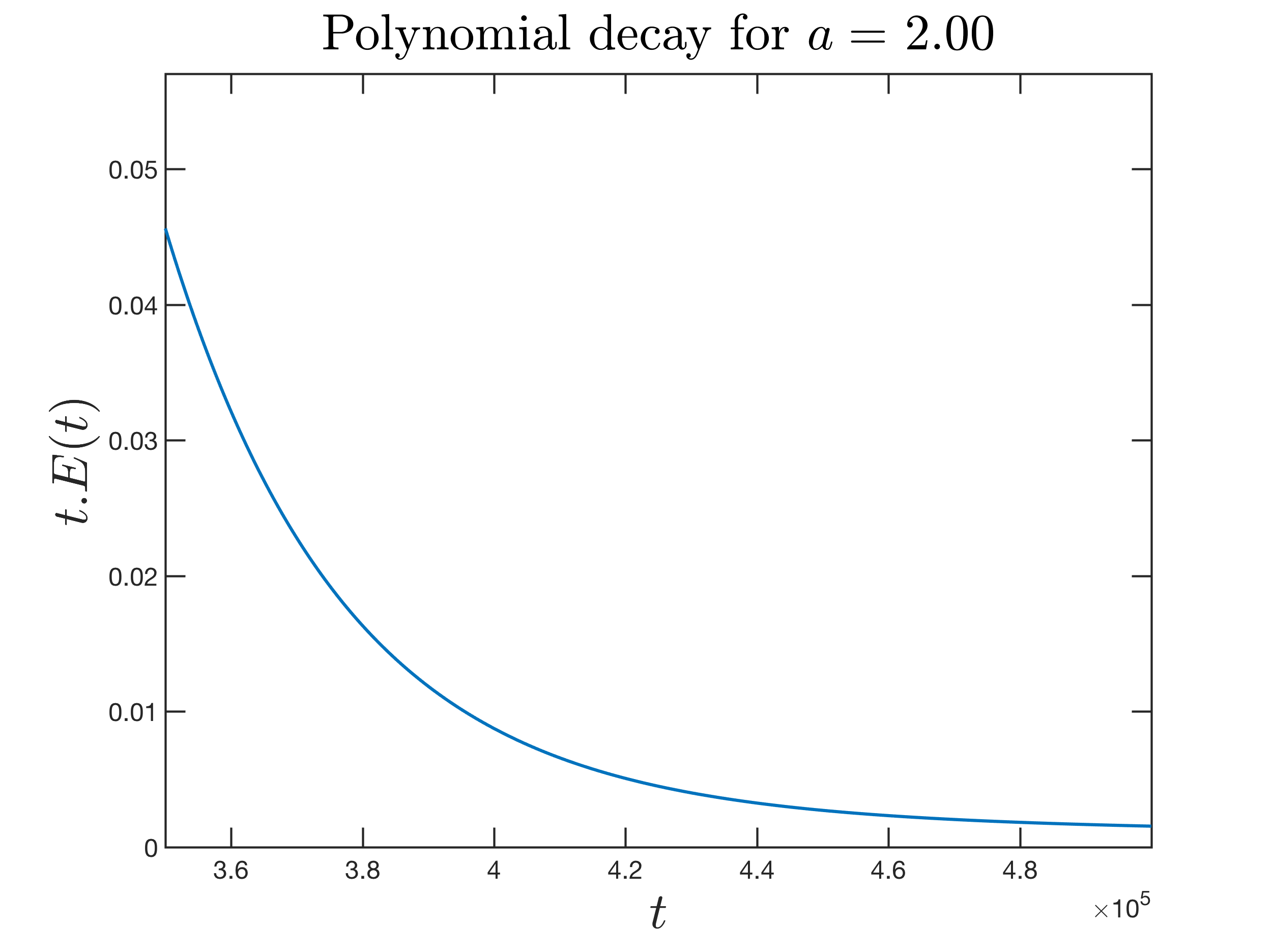}}\\[0.01\textheight]
	\subcaptionbox{Which exponent if polynomial decay?\label{Poly-b4-c3-a2}}
	{\includegraphics[width=0.4\textwidth]{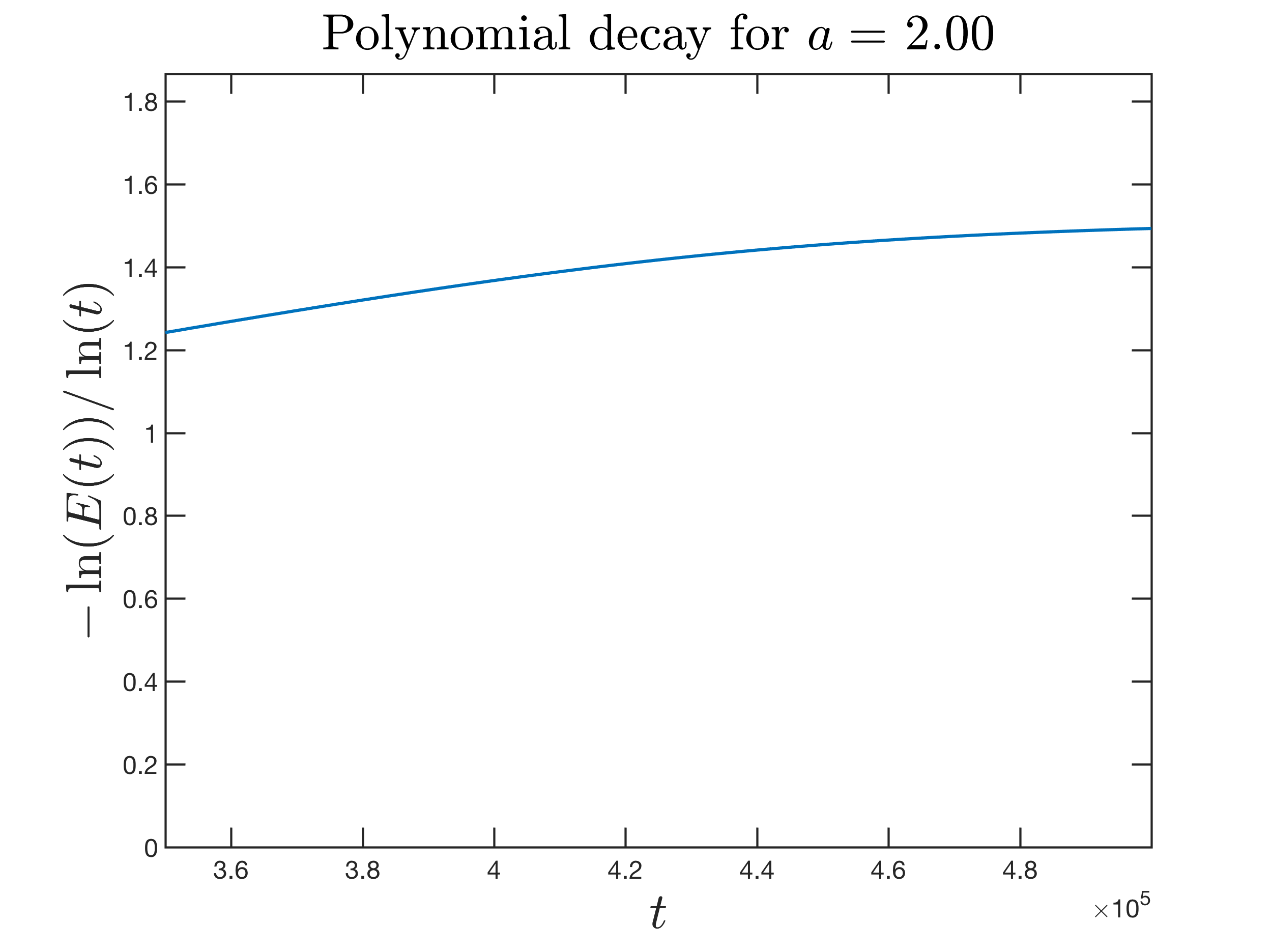}}
	\subcaptionbox{Final time profile.\label{b4-c3-a2}}
	{\includegraphics[width=0.4\textwidth]{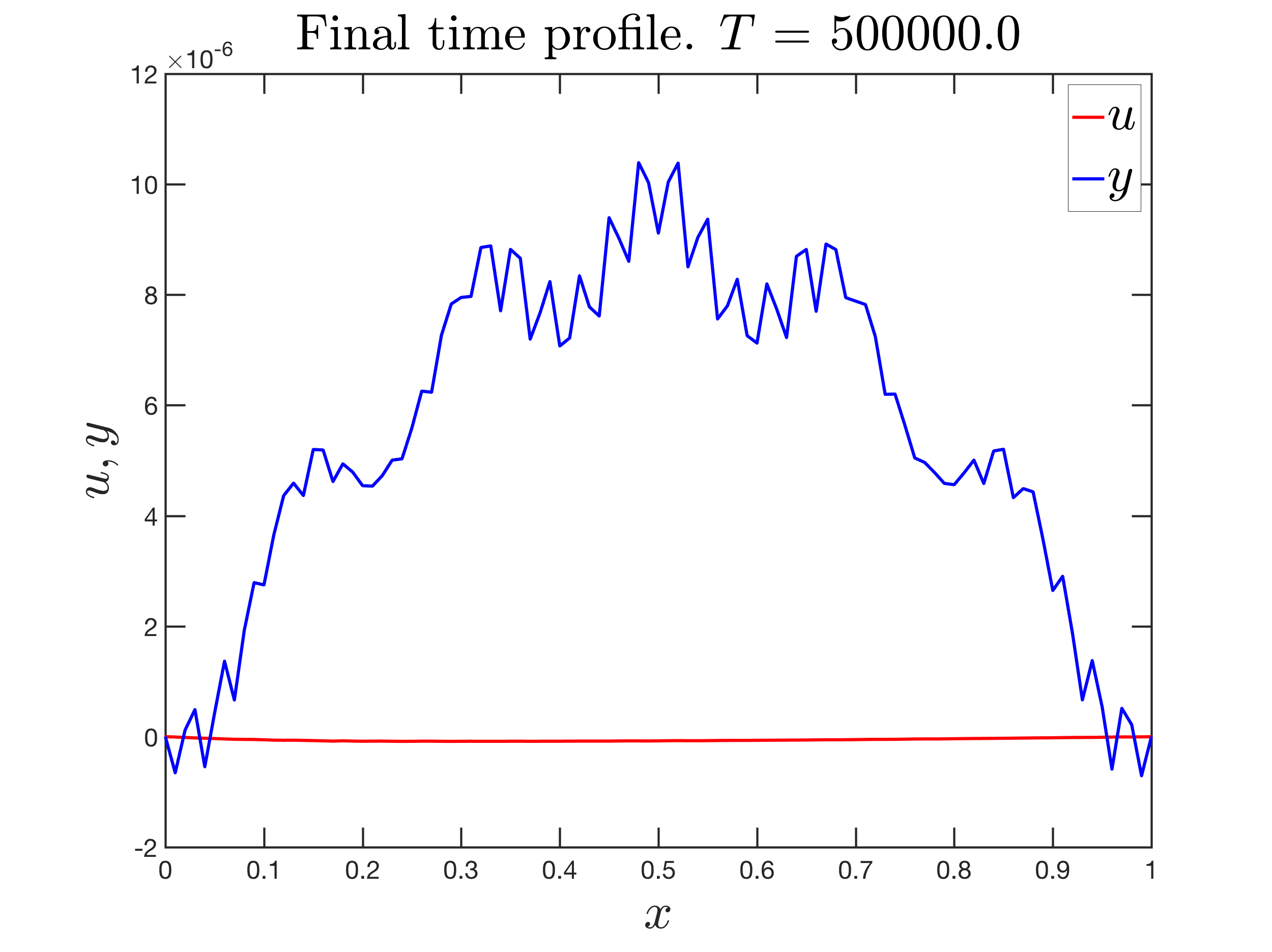}}
	\captionsetup{justification=centering}
	\caption{Long time behavior when $\omega_{b} \cap \omega_{c_{+}} \neq \emptyset$.\\
		\footnotesize{$b = b_{4}(x) = \mathds{1}_{[0.1,0.2]}(x)$ and $c = c_{3}(x)= \mathds{1}_{[0.1,0.2]\cup[0.8,0.9]}(x)$}}
	\label{Convergence-b4-c3-a2}
\end{figure}

To explore the speed of convergence to zero, we have plotted in figure \ref{Convergence-b4-c3-a2} $-\ln\big(E(t)\big) /t$ , $t \cdot E(t)$ and finally  $-\ln\big(E(t)\big) / \ln(t)$ versus $t$. Figure \ref{Exp-b4-c3-a2} shows clearly that $-\ln\big(E(t)\big) /t$ tends to zero and it permits to conclude that $E(t)$ tends to zero slower than an exponential.
Figure \ref{1-b4-c3-a2} permits to conclude that $E(t)$ tends to zero  faster than $1/t$. Finally
figure \ref{Poly-b4-c3-a2} shows that $E(t)$ tends to zero as $1/t^{\alpha}$ with $\alpha \simeq 1.4$. 

The final time profile presented in figure \ref{b4-c3-a2} confirms that $u$ and $y$ are small but it shows also that high frequencies for the unknown $y$ are not completely controlled.


\bigskip

\subsubsection{{$\boldsymbol{\omega_{b} \cap \omega_{c_{+}} = \emptyset}$. Unpredicted behavior}}
\begin{figure}[H]
	\centering
	\subcaptionbox{Final time $T = 500$.\label{Energy-small-b4-c5-a2}}
	{\includegraphics[width=0.4\textwidth]{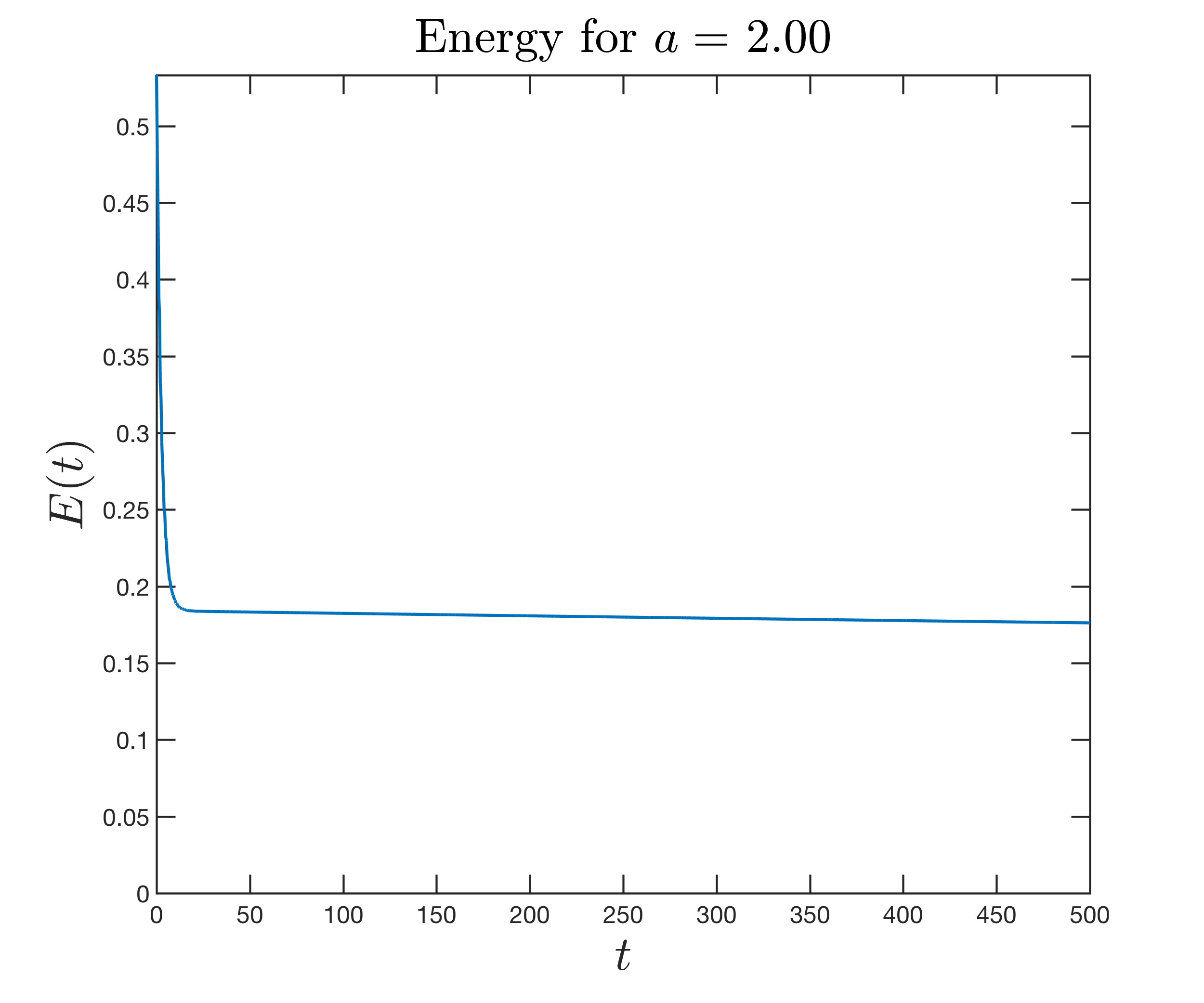}}
	\subcaptionbox{Final time $T= 500~000$.\label{Energy-b4-c5-a2}}
	{\includegraphics[width=0.4\textwidth]{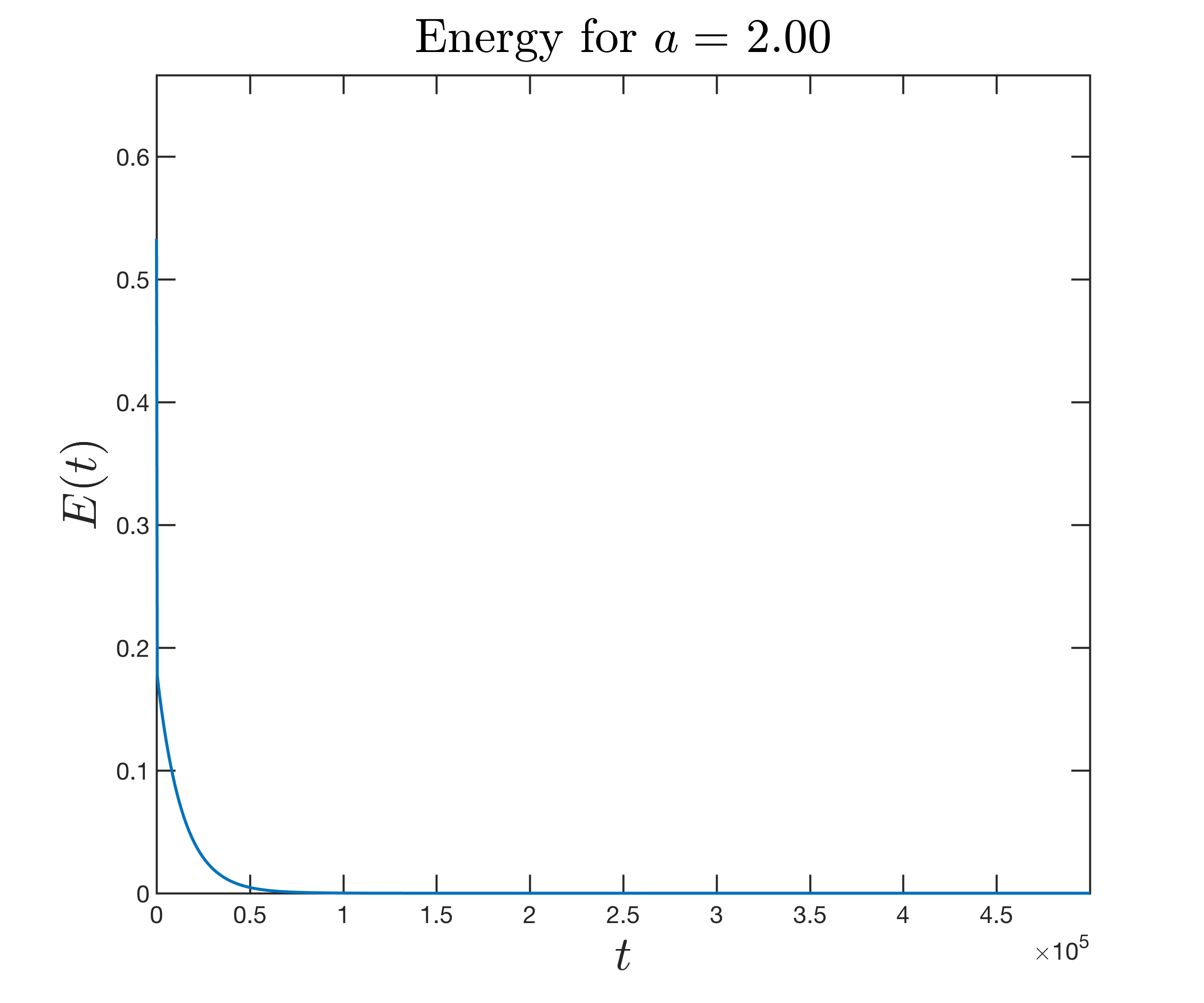}}
	\captionsetup{justification=centering}
	\caption{Energy when $\omega_{b} \cap \omega_{c_{+}} = \emptyset$.\\
		\footnotesize{$b = b_{4}(x) = \mathds{1}_{[0.1,0.2]}(x)$ and $c = c_{5}(x)= \mathds{1}_{[0.4,0.6]}(x)$}}
	\label{Ener-b4-c5-a2}
\end{figure}

At the numerical level, we are interested in the long time behavior of the solution $(u,y)$ when we suppose that  $\omega_{b} \cap \omega_{c_{+}} = \emptyset$.
For this sake, we present in figure \ref{Ener-b4-c5-a2}, the total energy  where we have chosen $b = b_{4}(x) = \mathds{1}_{[0.1,0.2]}(x)$ and $c = c_{5}(x)= \mathds{1}_{[0.4,0.6]}(x)$.

Again, when taking as final time $T = 500$, it seems that the energy does not tend to zero  as shown in figure \ref{Energy-small-b4-c5-a2}.
Taking as final time $T = 500~000$ , figure \ref{Energy-b4-c5-a2}  shows that the energy  goes finally to zero.

\begin{figure}[H]
	\centering
	\subcaptionbox{Exponential decay?\label{Exp-b4-c5-a2}}
	{\includegraphics[width=0.4\textwidth]{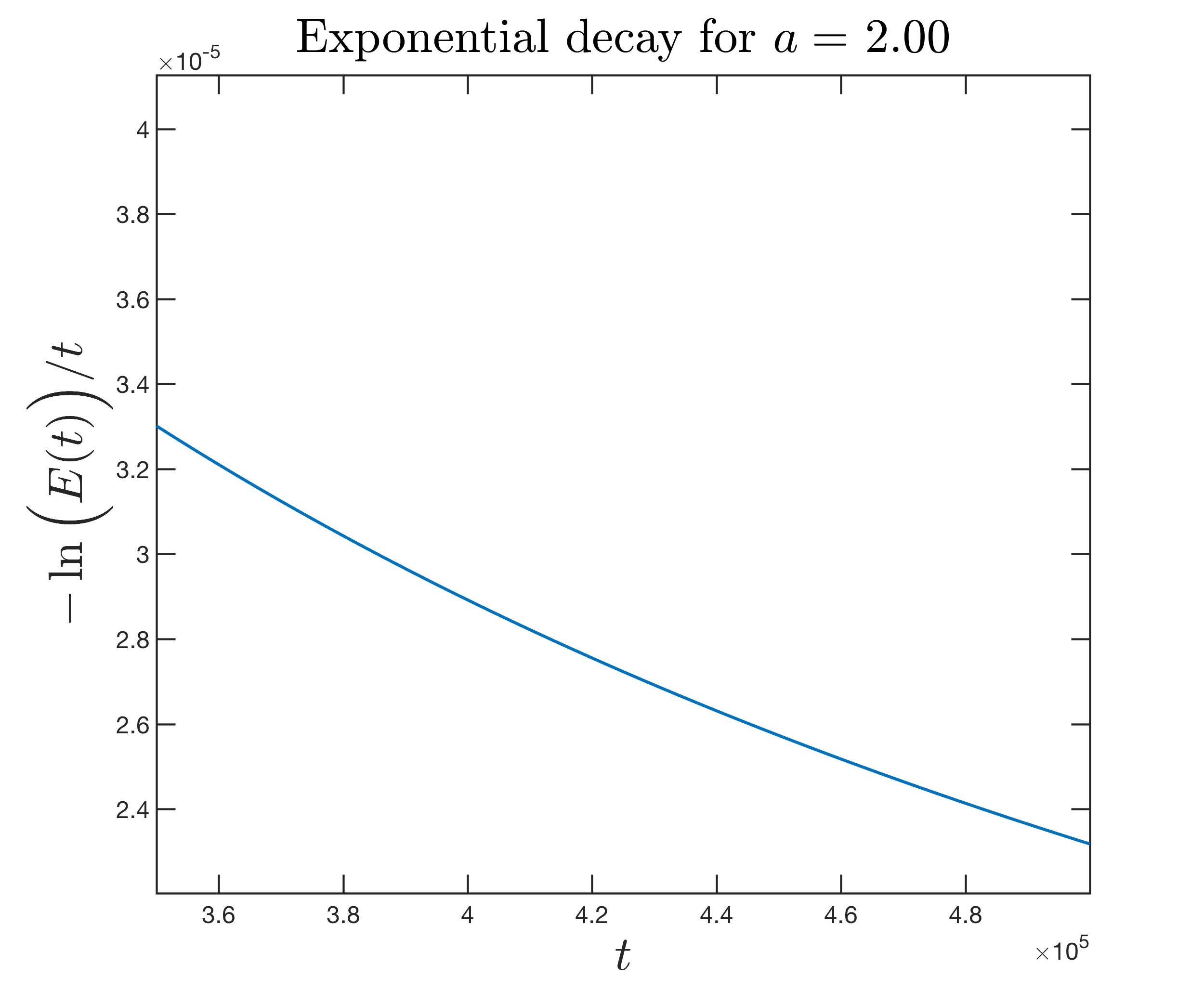}}
	\subcaptionbox{Polynomial decay in $1/t$?\label{1-b4-c5-a2}}
	{\includegraphics[width=0.4\textwidth]{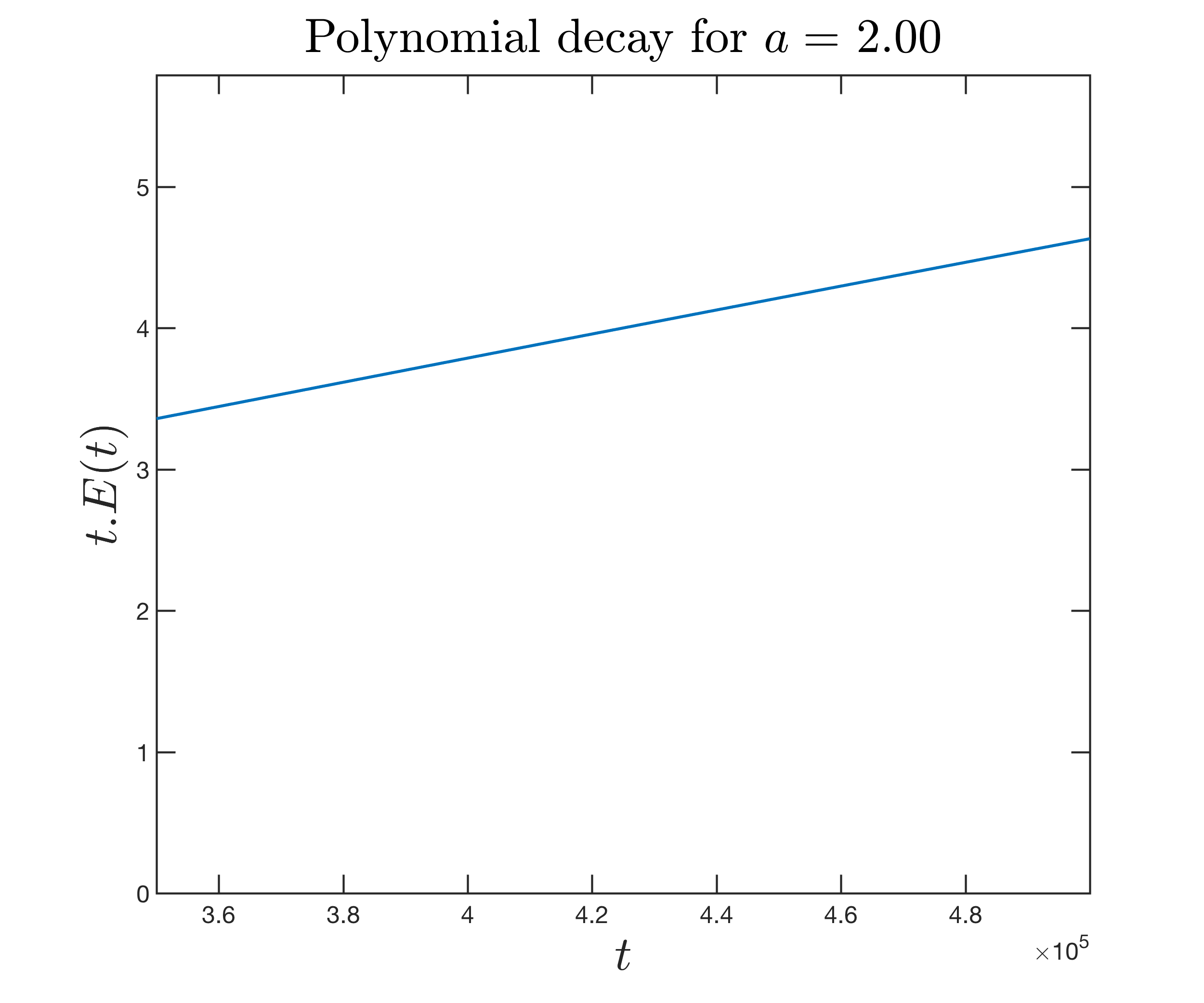}}\\[0.01\textheight]
	\subcaptionbox{Which exponent if polynomial decay?\label{Poly-b4-c5-a2}}
	{\includegraphics[width=0.4\textwidth]{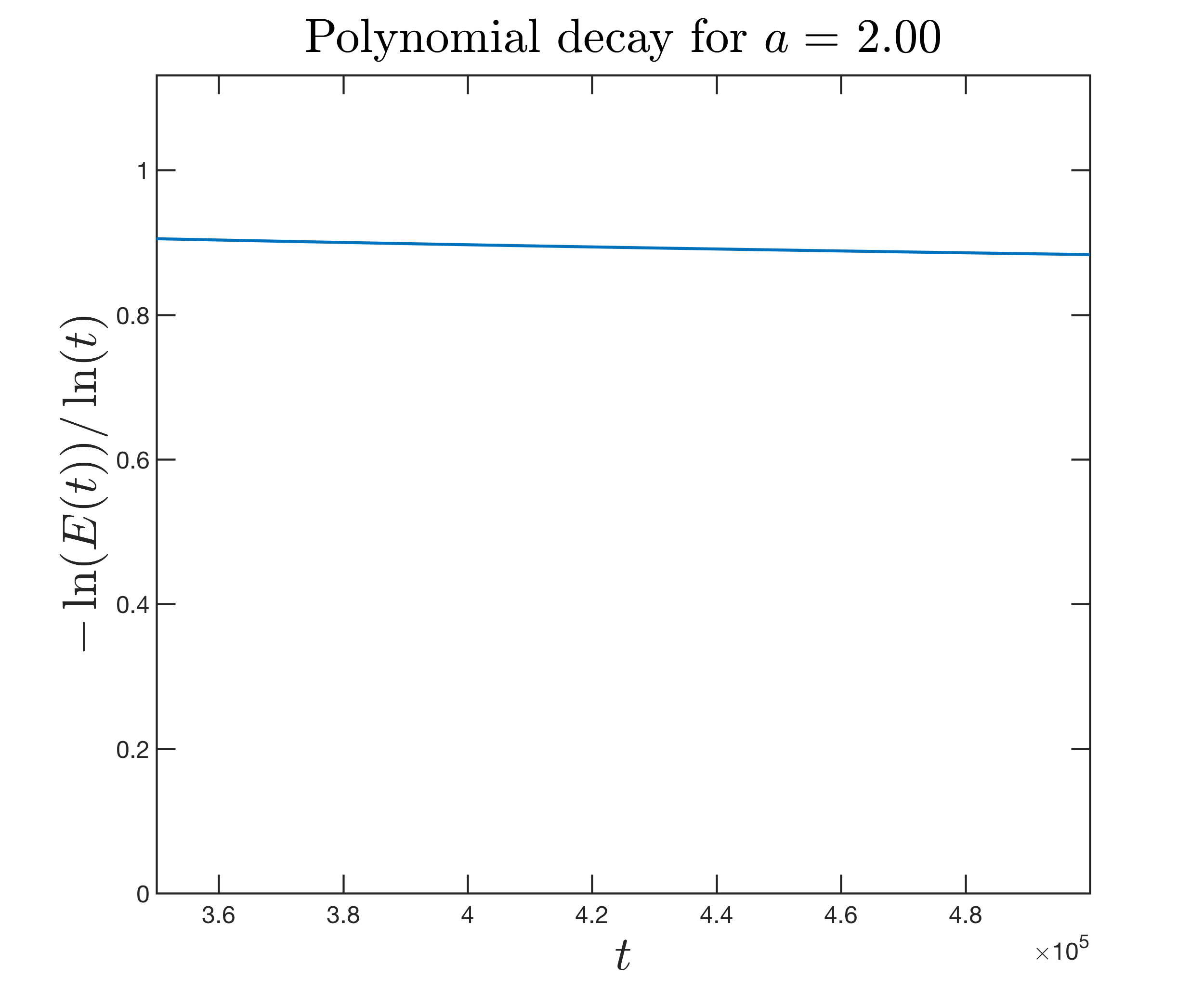}}
	\subcaptionbox{Final time profile.\label{b4-c5-a2}}
	{\includegraphics[width=0.4\textwidth]{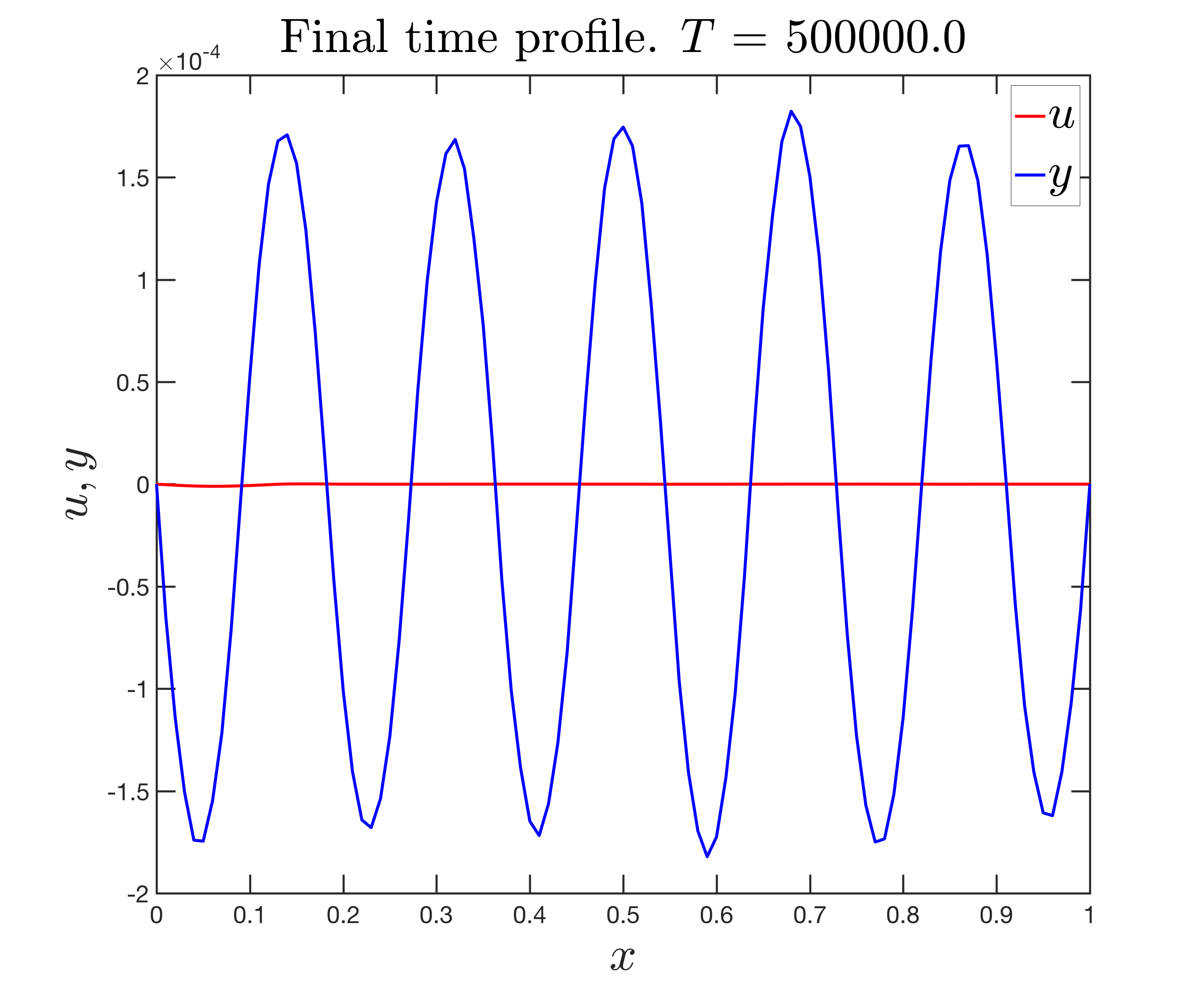}}
	\captionsetup{justification=centering}
	\caption{Long time behavior when $\omega_{b} \cap \omega_{c_{+}} = \emptyset$. \\
		\footnotesize{$b = b_{4}(x) = \mathds{1}_{[0.1,0.2]}(x)$ and $c = c_{5}(x)= \mathds{1}_{[0.4,0.6]}(x)$}}
	\label{Convergence-b4-c5-a2}
\end{figure}

To explore the speed of convergence to zero, we have plotted in figure \ref{Convergence-b4-c5-a2} $-\ln\big(E(t)\big) /t$ , $t \cdot E(t)$ and finally  $-\ln\big(E(t)\big) / \ln(t)$ versus $t$. Figure \ref{Exp-b4-c5-a2} shows clearly that $-\ln\big(E(t)\big) /t$ tends to zero and it permits to conclude that $E(t)$ tends to zero slower than an exponential   but figure \ref{1-b4-c5-a2} shows that 
$E(t)$ tends to zero slower than $1/t$. This fact is confirmed by figure \ref{Poly-b4-c5-a2} which shows that $E(t)$ tends to zero as $1/t^{\alpha}$ with $\alpha \simeq 0.9$. Eventually, taking a larger time
could conclude that the convergence is like $1/t$.
\medskip
\\Again, the final time profile presented in figure \ref{b4-c5-a2} confirms that $u$ and $y$ are small but it shows also that high frequencies for the unknown $y$ are not completely controlled.

As for the case when the two propagation speeds were identical this results was not predicted by the theoretical results.

\begin{figure}[H]
	\centering
	\subcaptionbox{Final time $T = 500$.\label{Energy-small-b5-c4-a2}}
	{\includegraphics[width=0.4\textwidth]{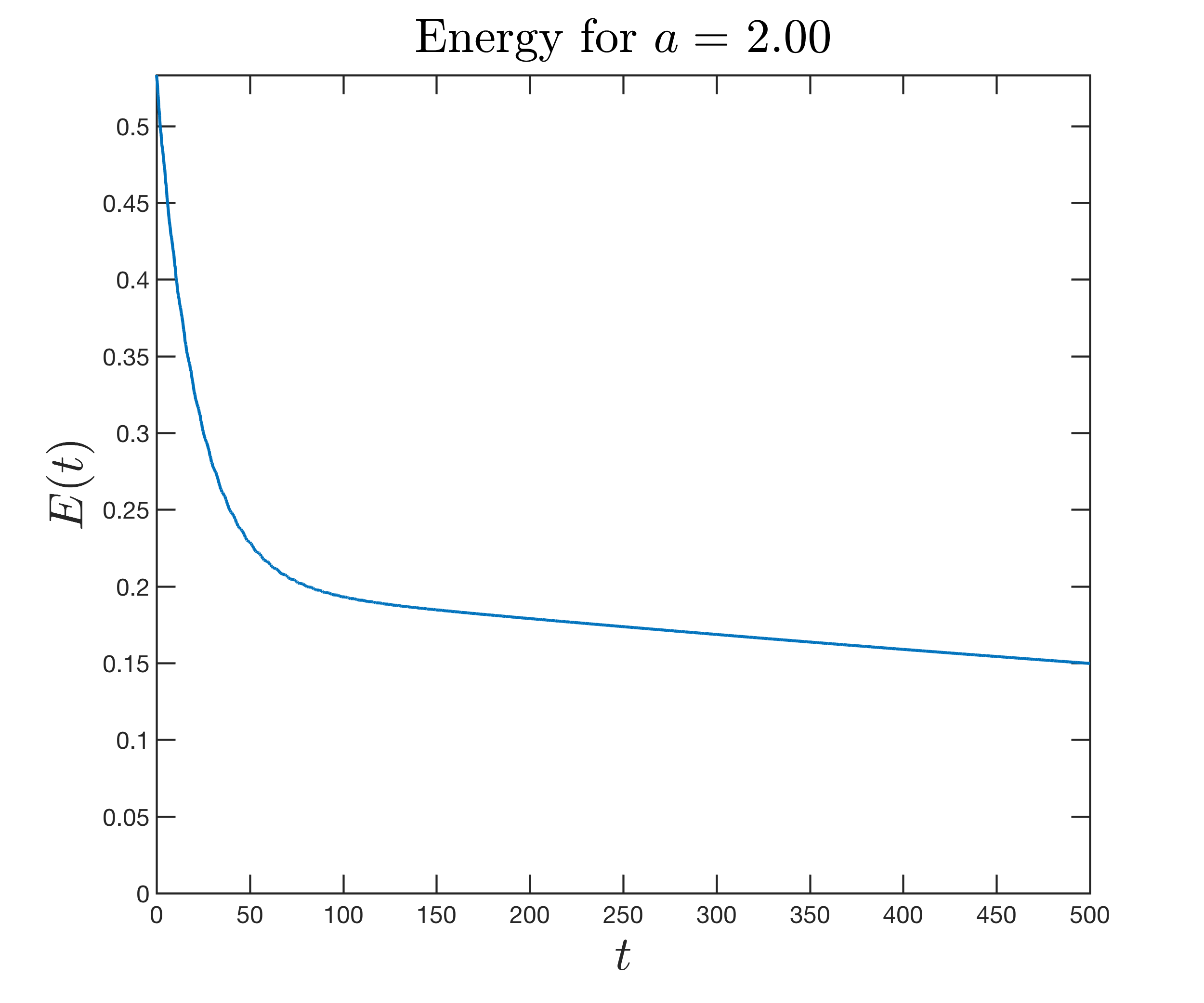}}
	\subcaptionbox{Final time $T= 500~000$.\label{Energy-b5-c4-a2}}
	{\includegraphics[width=0.4\textwidth]{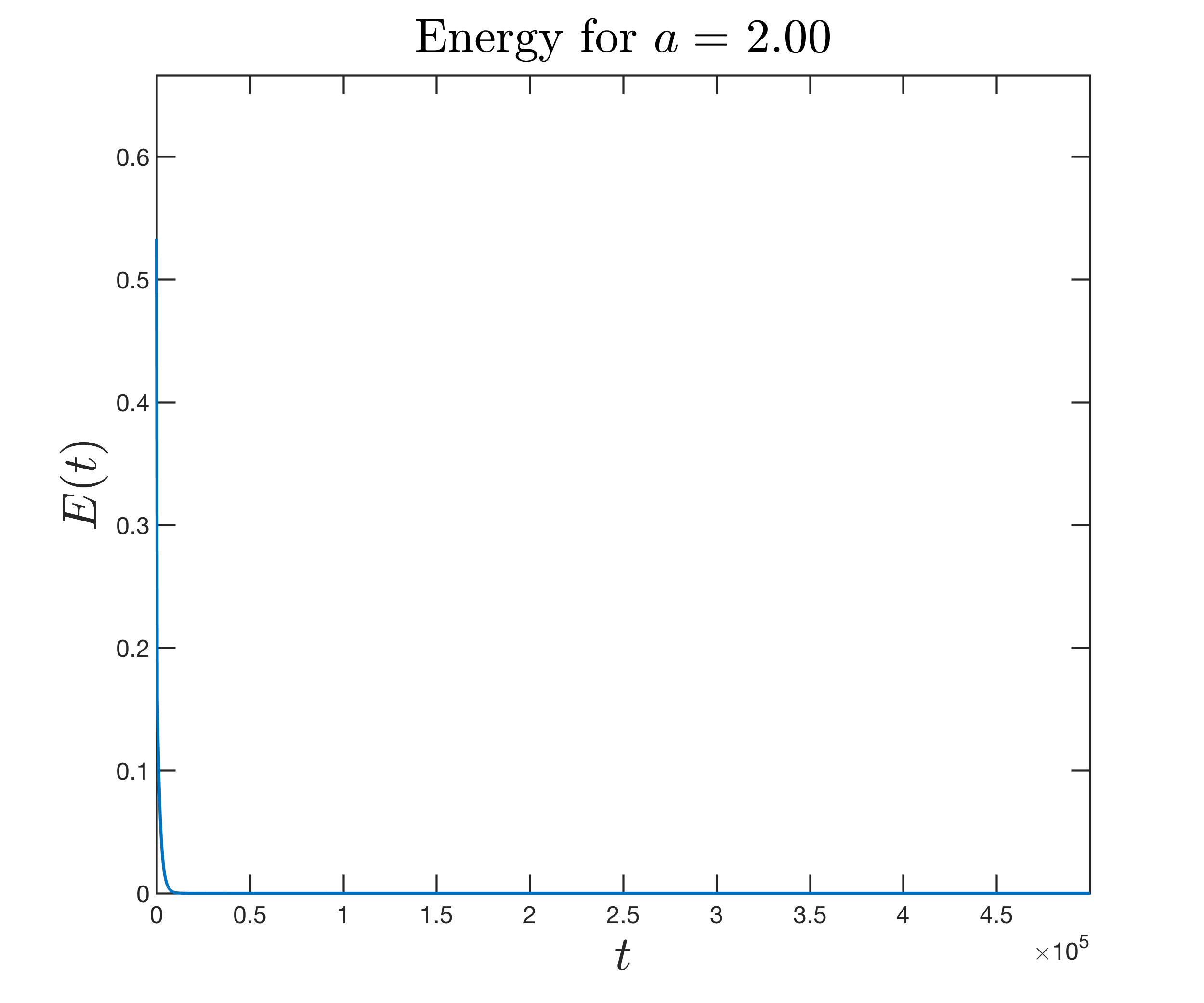}}
	\captionsetup{justification=centering}
	\caption{Energy when $\omega_{b} \cap \omega_{c_{+}} = \emptyset$.\\
		\footnotesize{$b = b_{4}(x) = \mathds{1}_{[0.1,0.2]}(x)$ and $c = c_{5}(x)= \mathds{1}_{[0.4,0.6]}(x)$}}
	\label{Ener-b5-c4-a2}
\end{figure}

So we decided to confirm this  behavior by choosing $b = b_{5}(x) = \mathds{1}_{[0.4,0.6]}(x)$ and \linebreak
$c = c_{4}(x)= \mathds{1}_{[0.1,0.2]}(x)$.
Again, when taking as final time $T = 500$, it seems that the energy does not tends to zero  as shown in figure \ref{Energy-small-b5-c4-a2}.
Taking as final time $T = 500~000$, figure \ref{Energy-b4-c5-a2}  shows that the energy goes finally to zero.

\begin{figure}[H]
	\centering
	\subcaptionbox{Exponential decay?\label{Exp-b5-c4-a2}}
	{\includegraphics[width=0.4\textwidth]{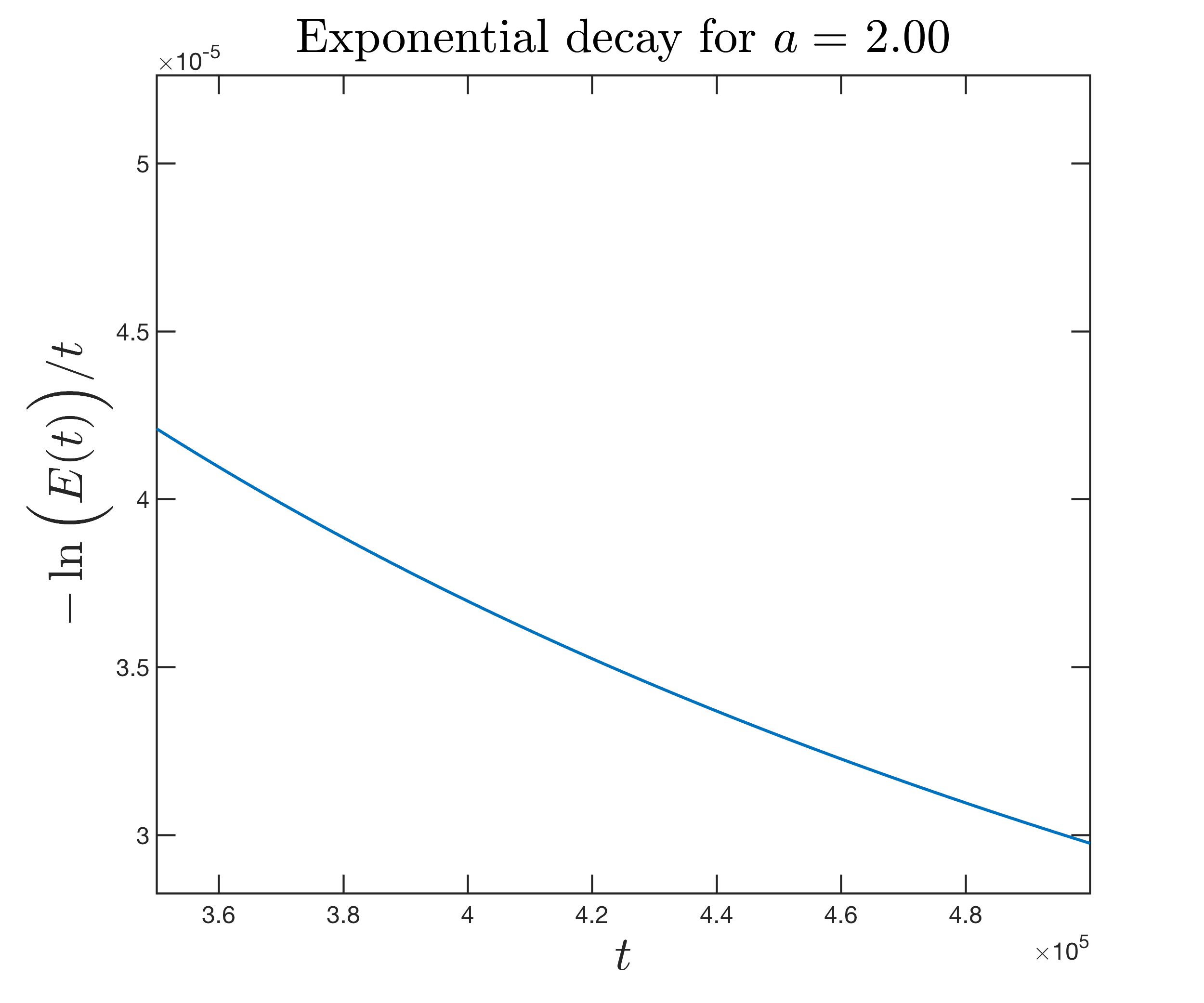}}
	\subcaptionbox{Polynomial decay in $1/t$?\label{1-b5-c4-a2}}
	{\includegraphics[width=0.4\textwidth]{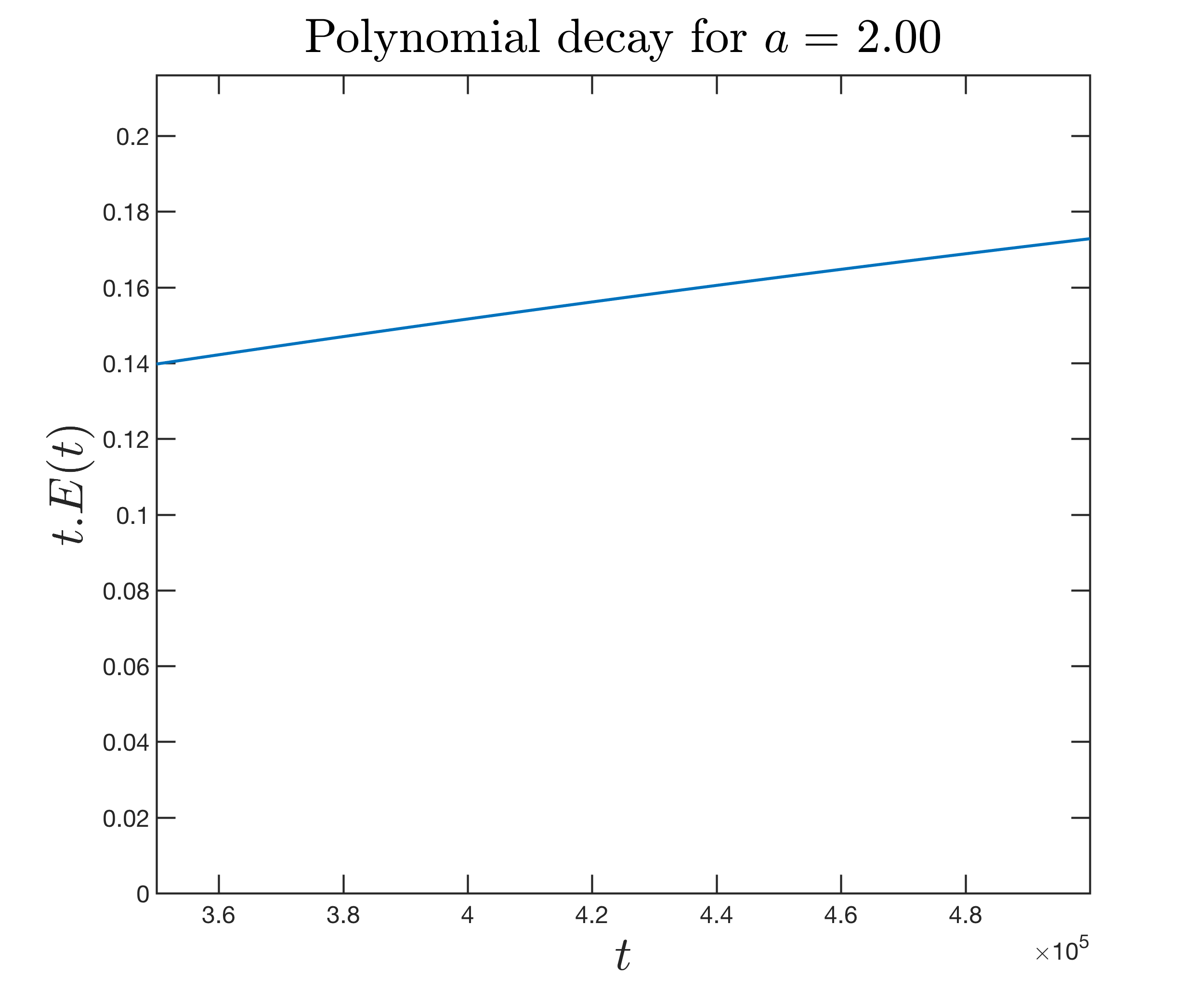}}\\[0.01\textheight]
	\subcaptionbox{Which exponent if polynomial decay?\label{Poly-b5-c4-a2}}
	{\includegraphics[width=0.4\textwidth]{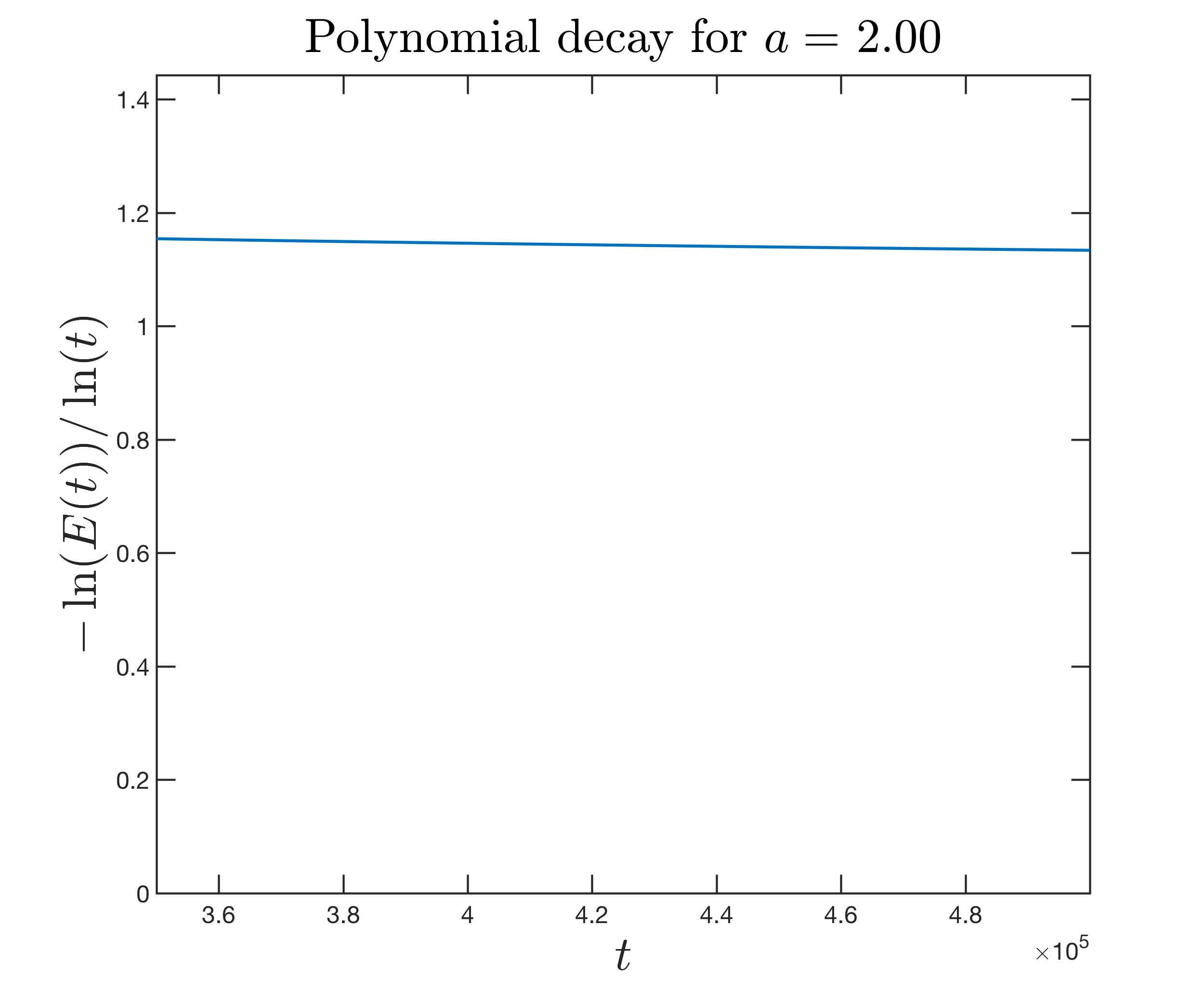}}
	\subcaptionbox{Final time profile.\label{b5-c4-a2}}
	{\includegraphics[width=0.4\textwidth]{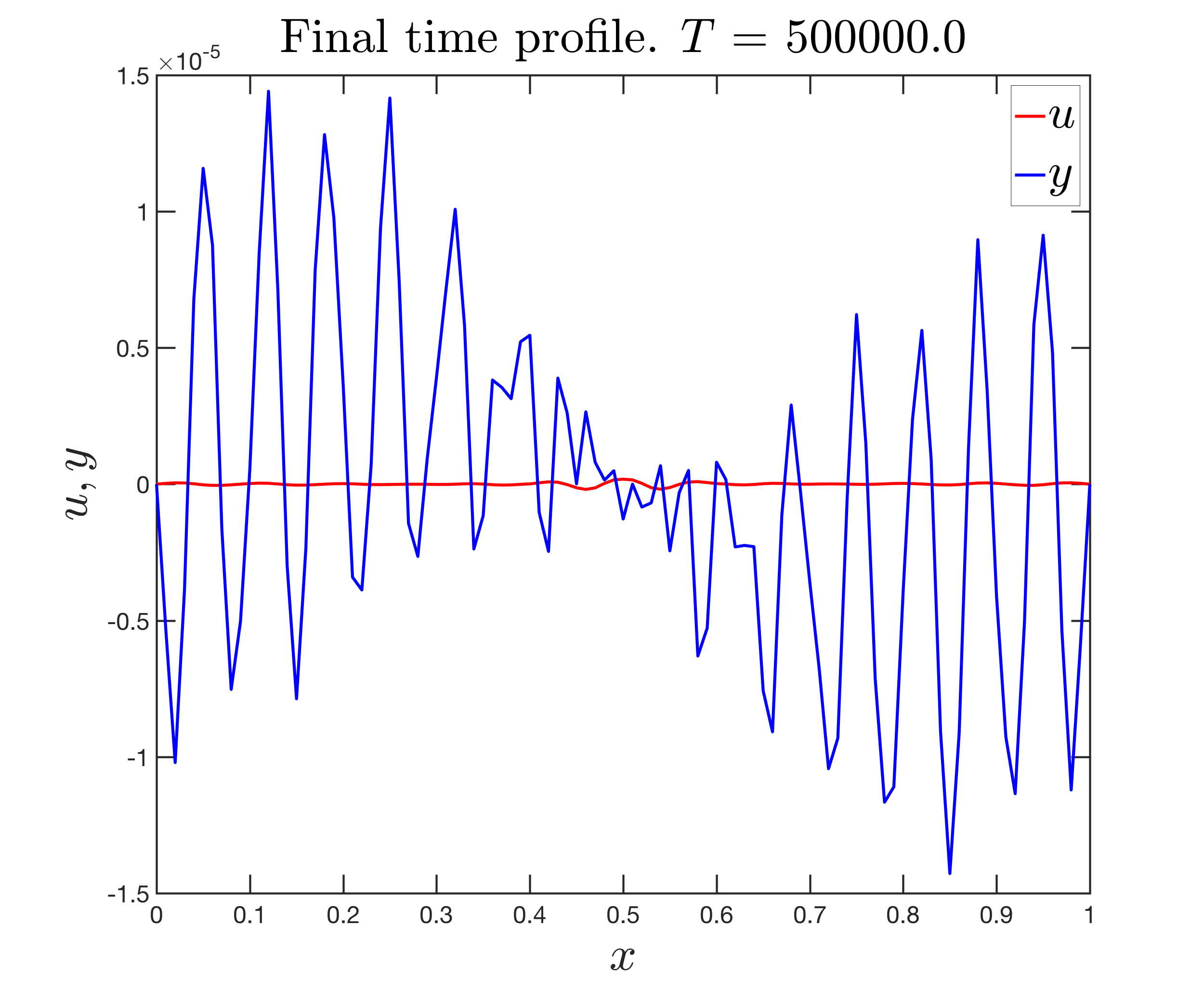}}
	\captionsetup{justification=centering}
	\caption{Long time behavior when $\omega_{b} \cap \omega_{c_{+}} = \emptyset$. \\
		\footnotesize{$b = b_{5}(x) = \mathds{1}_{[0.4,0.6]}(x)$ and $c = c_{4}(x)= \mathds{1}_{[0.1,0.2]}(x)$}}
	\label{Convergence-b5-c4-a2}
\end{figure}

To explore the speed of convergence to zero, we have plotted in figure \ref{Convergence-b5-c4-a2} $-\ln\big(E(t)\big) /t$ , $t \cdot E(t)$ and finally  $-\ln\big(E(t)\big) / \ln(t)$ versus $t$. 
Figure \ref{Exp-b5-c4-a2} shows clearly that $-\ln\big(E(t)\big) /t$ tends to zero and it permits to conclude that $E(t)$  tends to zero slower than an exponential and 
figure \ref{1-b5-c4-a2} permits to conclude that the convergence is faster than $1/t$.
Finally figure \ref{Poly-b5-c4-a2} shows that $E(t)$ tends to zero as $1/t^{\alpha}$ with $\alpha \simeq 1.19$. 

\begin{remark}
	The final time profile presented in figure \ref{b4-c3-a2} , figure \ref{b4-c5-a2} and figure \ref{b5-c4-a2} confirms that $u$ and $y$ are small 
	but it shows also that high frequencies for the unknown $y$ are not completely controlled.
\end{remark}
\subsection{Different propagation speed:  $\boldsymbol{a < 1}$}
When $a \neq 1$, in order to see if the same behavior occurs no matter if $a$ is greater or less than 1, we investigate now the long time behavior of $(u,y)$ when the propagation speeds  is less than the one of $y$ namely $a < 1$. We have chosen $a = 0.5$.

\subsubsection{{$\boldsymbol{\omega_{b} \cap \omega_{c_{+}} \neq \emptyset}$. Polynomial stability}}
\begin{figure}[H]
	\centering
	\subcaptionbox{Final time $T = 500$.\label{Energy-small-b4-c3-a1/2}}
	{\includegraphics[width=0.4\textwidth]{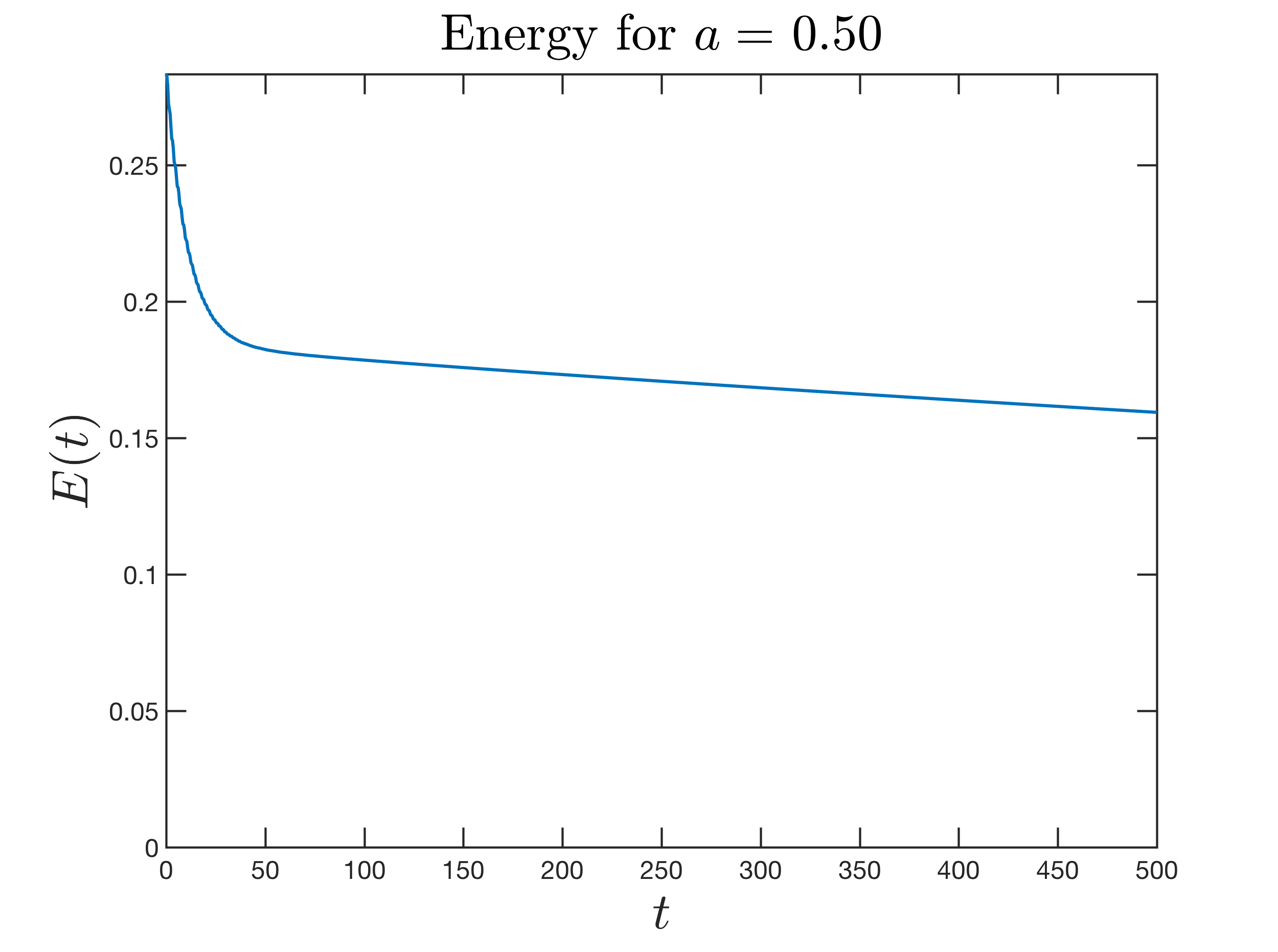}}
	\subcaptionbox{Final time $T= 500~000$.\label{Energy-b4-c3-a1/2}}
	{\includegraphics[width=0.4\textwidth]{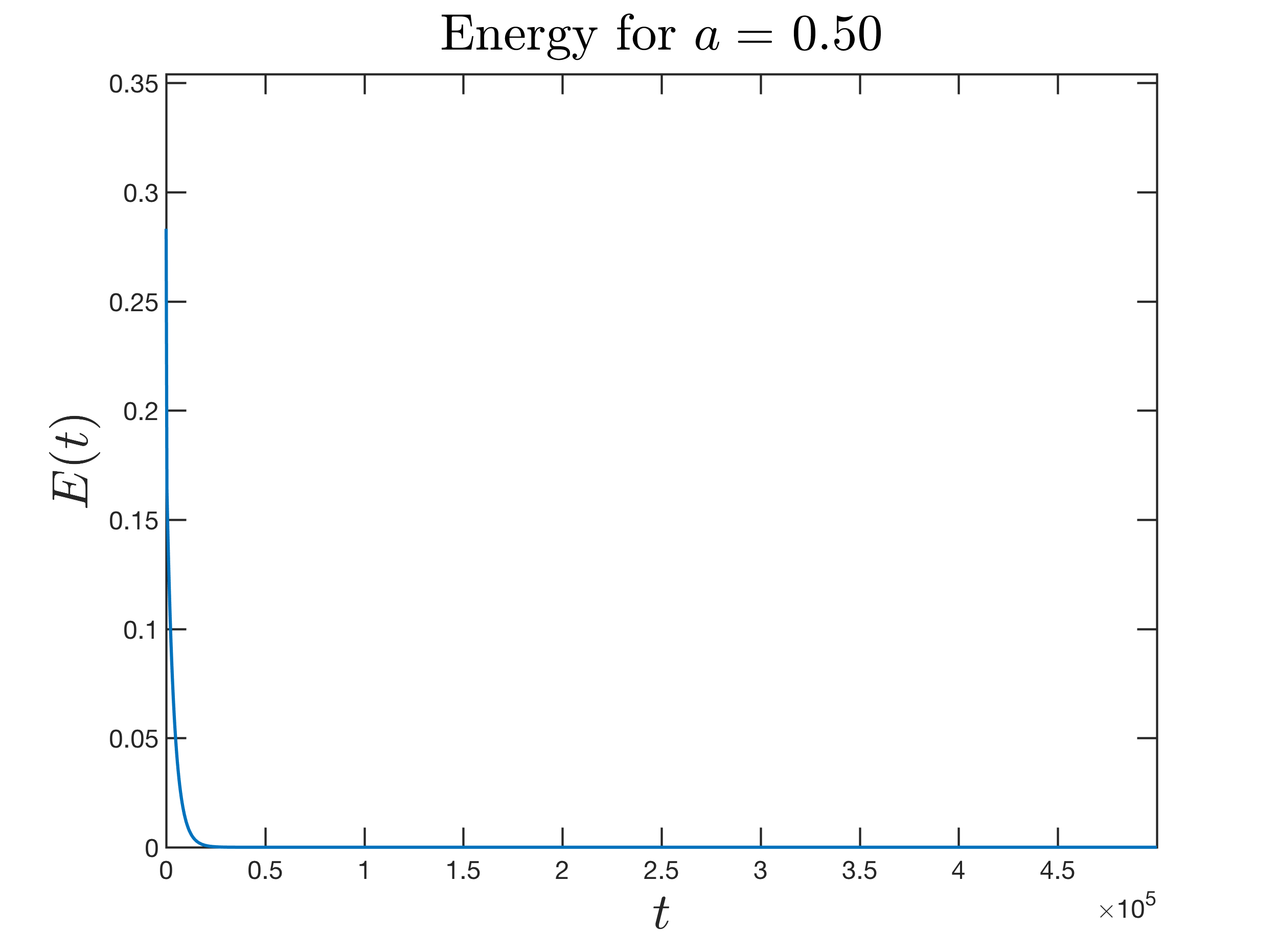}}\\[0.01\textheight]
	\captionsetup{justification=centering}
	\caption{Energy when $\omega_{b} \cap \omega_{c_{+}} \neq \emptyset$.\\
		\footnotesize{$b = b_{4}(x) = \mathds{1}_{[0.1,0.2]}(x)$ and $c = c_{3}(x)= \mathds{1}_{[0.1,0.2]\cup[0.8,0.9]}(x)$}}
	\label{Ener-b4-c3-a1/2}
\end{figure}

Let us now verify the theoretical results when we suppose that  $\omega_{b} \cap \omega_{c_{+}} \neq \emptyset$.
For this sake, we present in figure \ref{Energy-small-b4-c3-a1/2}, the total energy 
where we have chosen $b = b_{4}(x) = \mathds{1}_{[0.1,0.2]}(x)$ and $c = c_{3}(x)= \mathds{1}_{[0.1,0.2]\cup[0.8,0.9]}(x)$.

When taking as final time $T = 500$, it seems that the energy does not tend to zero as shown in figure \ref{Energy-small-b4-c3-a1/2}.
Taking as final time $T = 500~000$, figure \ref{Energy-b4-c3-a1/2}  shows that the energy  goes finally to zero.

\begin{figure}[H]
	\centering
	\subcaptionbox{Exponential decay?\label{Exp-b4-c3-a1/2}}
	{\includegraphics[width=0.4\textwidth]{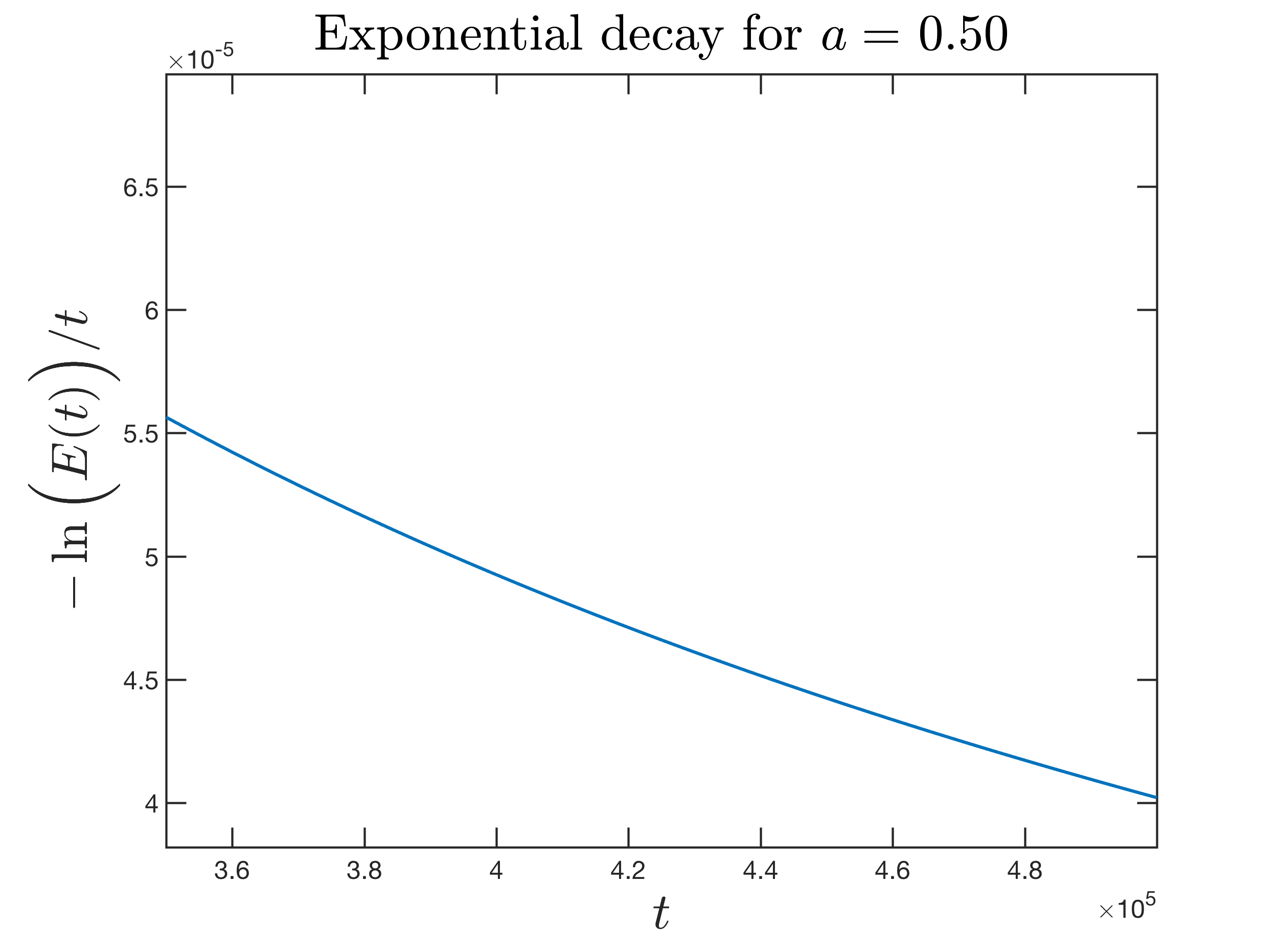}}
	\subcaptionbox{Polynomial decay in $1/t$?\label{1-b4-c3-a1/2}}
	{\includegraphics[width=0.4\textwidth]{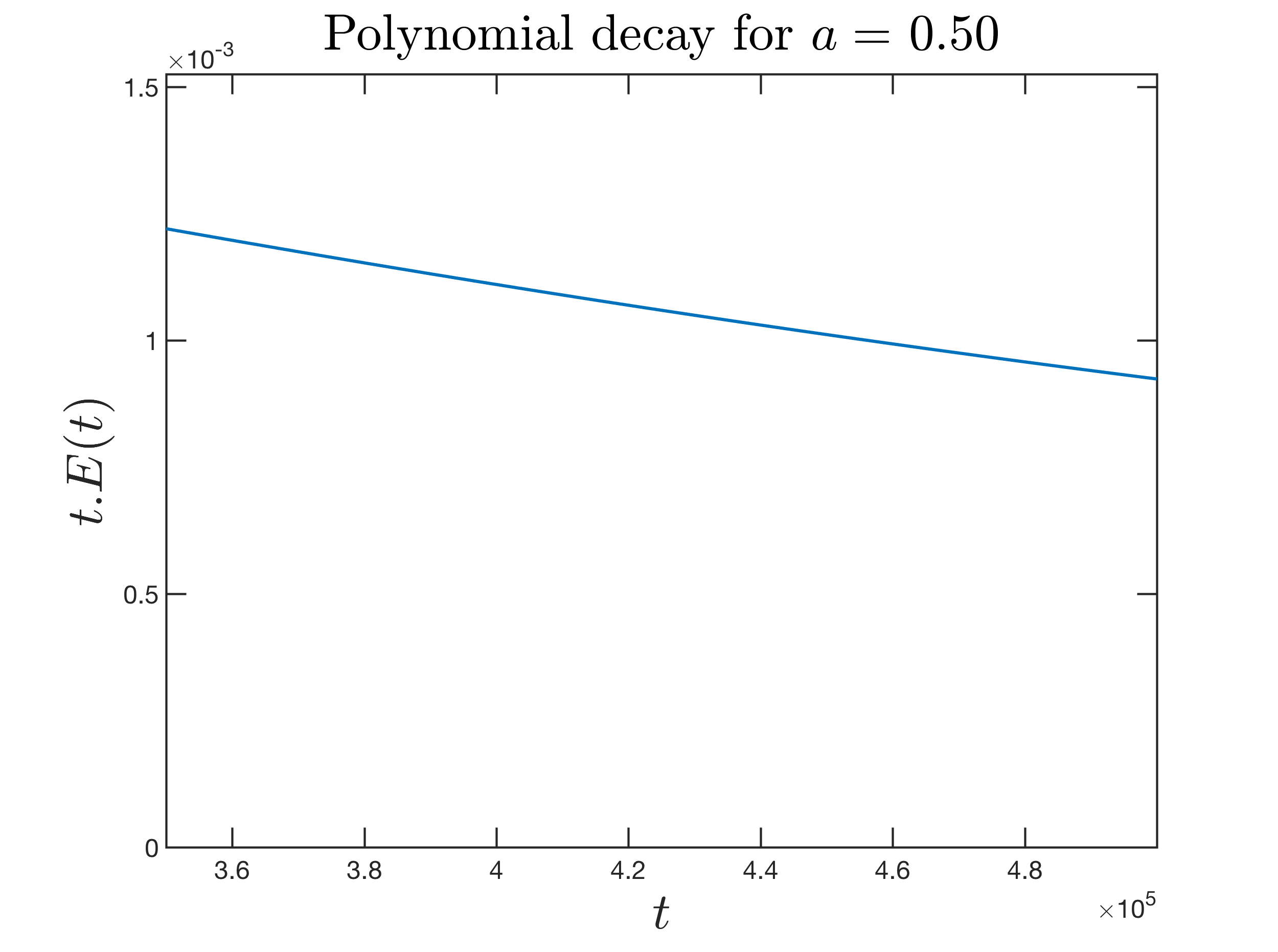}}\\[0.01\textheight]
	\subcaptionbox{Which exponent if polynomial decay?\label{Poly-b4-c3-a1/2}}
	{\includegraphics[width=0.4\textwidth]{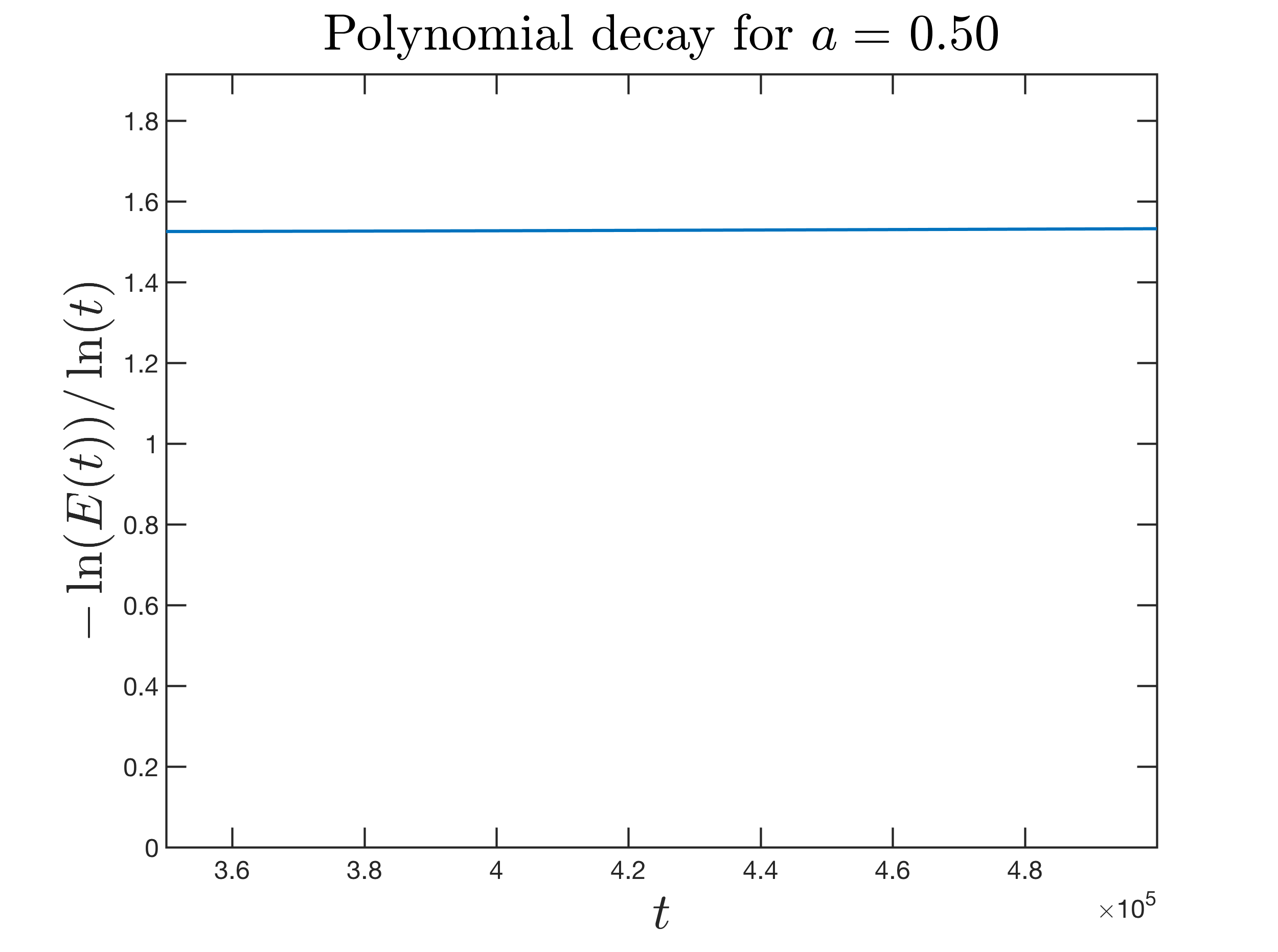}}
	\subcaptionbox{Final time profile.\label{b4-c3-a1/2}}
	{\includegraphics[width=0.4\textwidth]{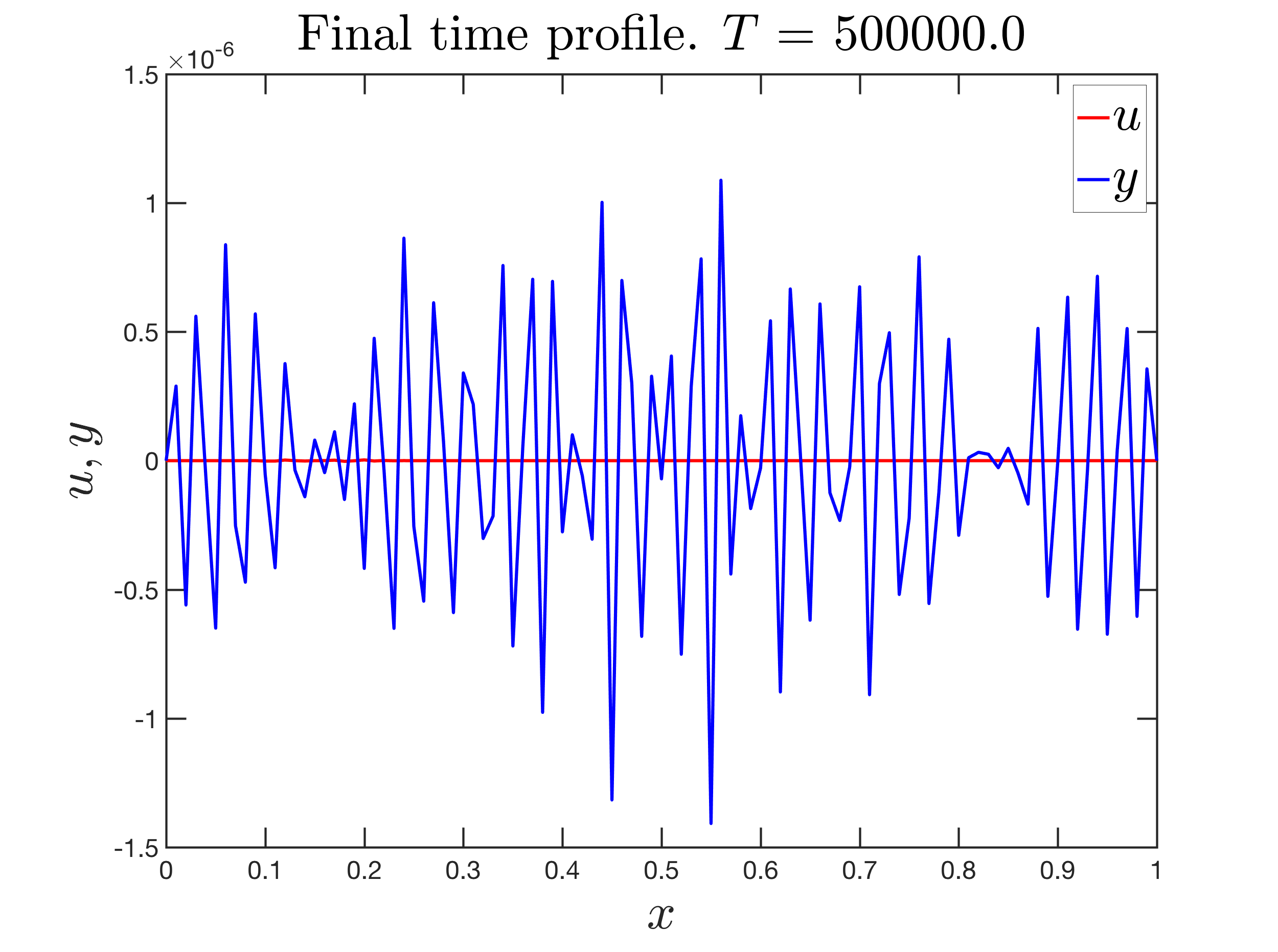}}
	\captionsetup{justification=centering}
	\caption{Long time behavior when $\omega_{b} \cap \omega_{c_{+}} \neq \emptyset$.\\
		\footnotesize{$b = b_{4}(x) = \mathds{1}_{[0.1,0.2]}(x)$ and $c = c_{3}(x)= \mathds{1}_{[0.1,0.2]\cup[0.8,0.9]}(x)$}}
	\label{Convergence-b4-c3-a1/2}
\end{figure}

To explore the speed of convergence to zero, we have plotted in figure \ref{Convergence-b4-c3-a1/2} $-\ln\big(E(t)\big) /t$ , $t \cdot E(t)$ and finally  $-\ln\big(E(t)\big) / \ln(t)$ versus $t$. Figure \ref{Exp-b4-c3-a1/2} shows clearly that $-\ln\big(E(t)\big) /t$ tends to zero slower than an exponential.
Figure \ref{1-b4-c3-a1/2} permits to conclude that $E(t)$ tends to zero  faster than $1/t$. Finally
figure \ref{Poly-b4-c3-a1/2} shows that $E(t)$ tends to zero as $1/t^{\alpha}$ with $\alpha \simeq 1.5$. 

The final time profile confirms that $u$ and $y$ are small but
it shows also that high frequencies for the unknown $y$ are not completely controlled.

\bigskip
\subsubsection{{$\boldsymbol{\omega_{b} \cap \omega_{c_{+}} = \emptyset}$: Unpredicted behavior}}
\begin{figure}[H]
	\centering
	\subcaptionbox{Final time $T = 500$.\label{Energy-small-b4-c5-a1/2}}
	{\includegraphics[width=0.4\textwidth]{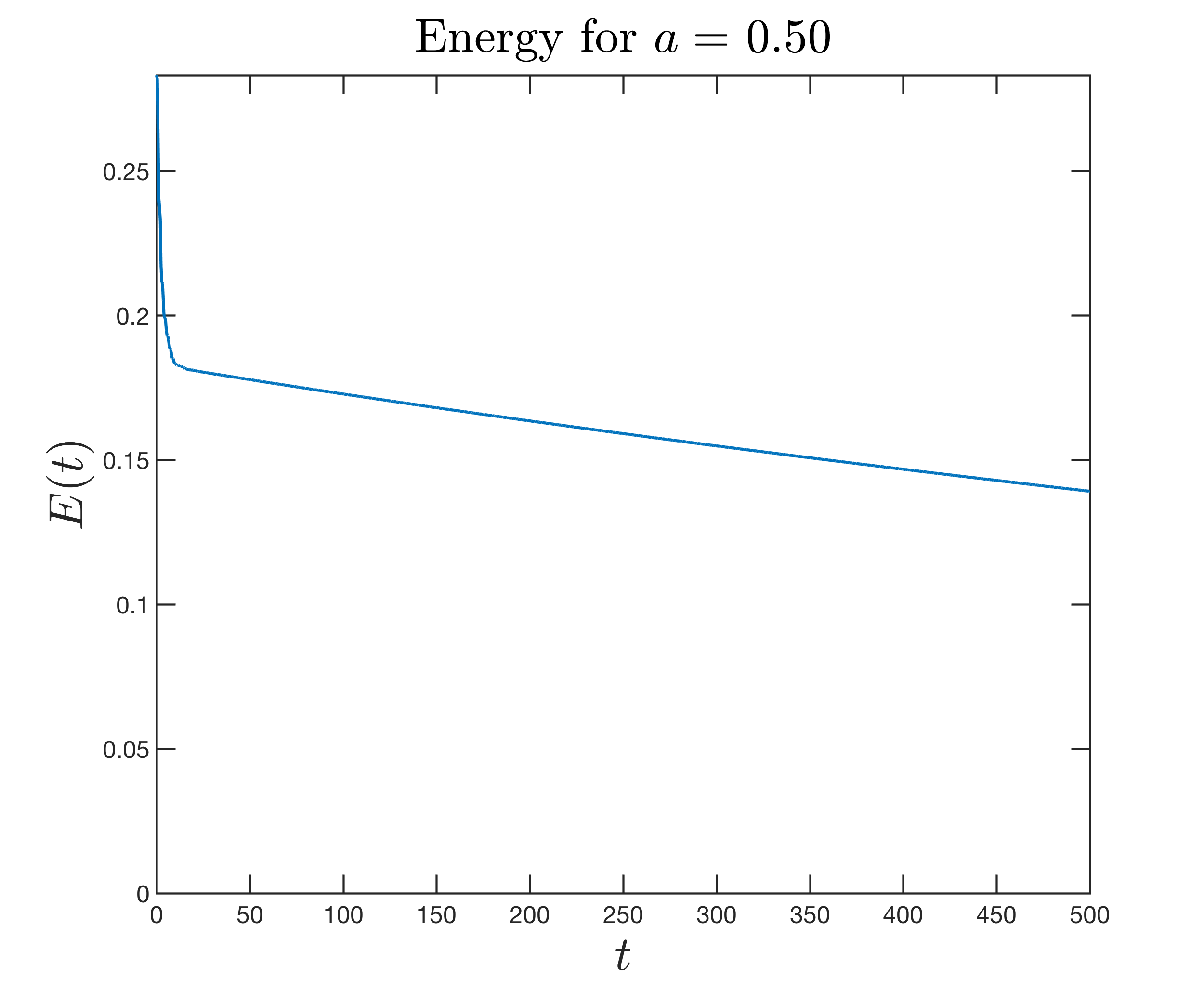}}
	\subcaptionbox{Final time $T= 500~000$.\label{Energy-b4-c5-a1/2}}
	{\includegraphics[width=0.4\textwidth]{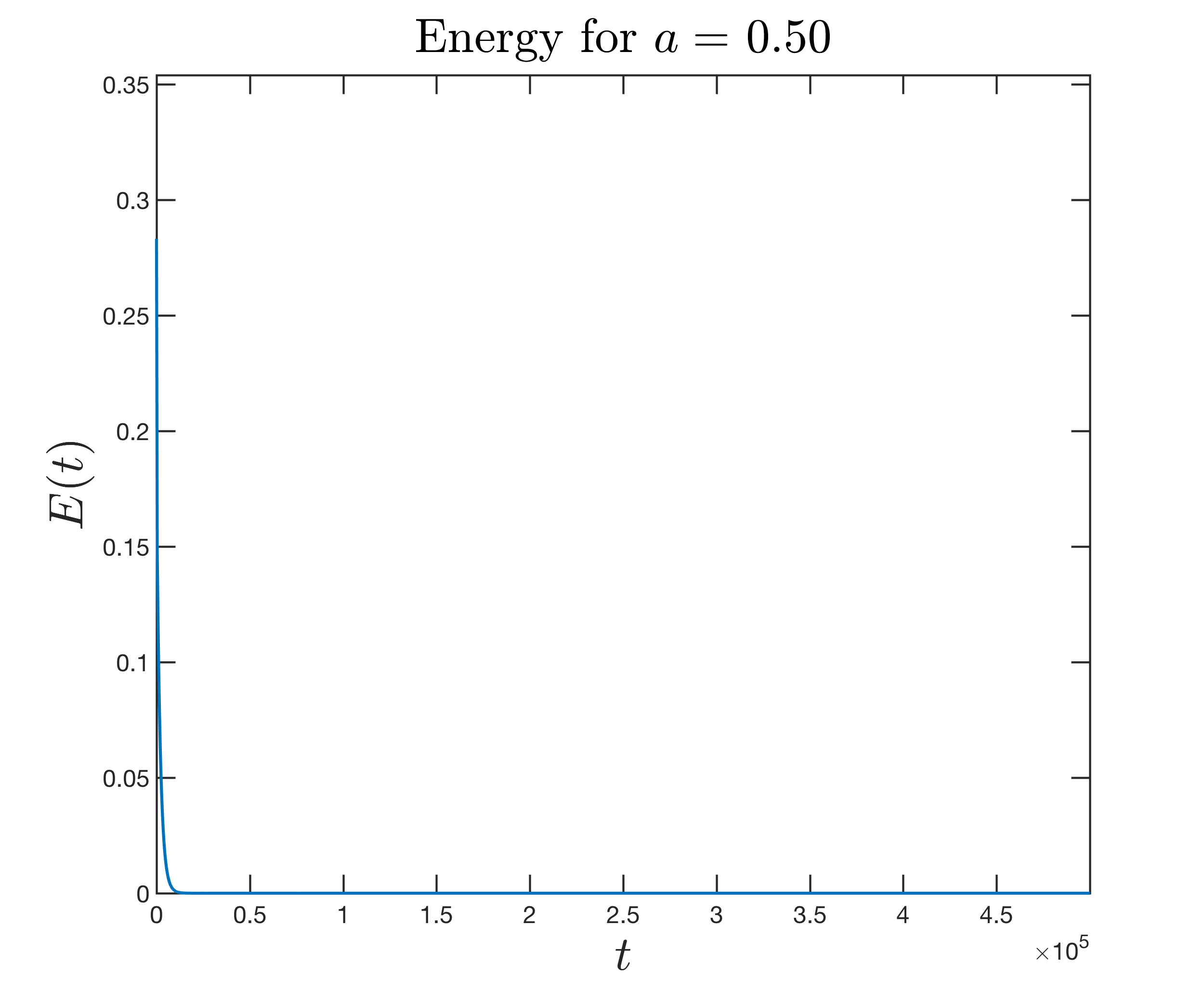}}
	\captionsetup{justification=centering}
	\caption{Energy when $\omega_{b} \cap \omega_{c_{+}} = \emptyset$.\\
		\footnotesize{$b = b_{4}(x) = \mathds{1}_{[0.1,0.2]}(x)$ and $c = c_{5}(x)= \mathds{1}_{[0.4,0.6]}(x)$}}
	\label{Ener-b4-c5-a1/2}
\end{figure}

Again, the numerical level, we are interested in the long time behavior of the solution $(u,y)$ when we suppose that  $\omega_{b} \cap \omega_{c_{+}} = \emptyset$.
For this sake, we present in figure \ref{Energy-small-b4-c5-a1/2}, the total energy  where we have chosen $b = b_{4}(x) = \mathds{1}_{[0.1,0.2]}(x)$ and $c = c_{5}(x)= \mathds{1}_{[0.4,0.6]}(x)$.

Again, when taking as final time $T = 500$, it seems that the energy does not tend to zero  as shown in figure \ref{Energy-small-b4-c5-a1/2}.
Taking as final time $T = 500~000$, figure \ref{Energy-b4-c5-a1/2}  shows that the energy finally goes to zero.

\begin{figure}[H]
	\centering
	\subcaptionbox{Exponential decay?\label{Exp-b4-c5-a1/2}}
	{\includegraphics[width=0.4\textwidth]{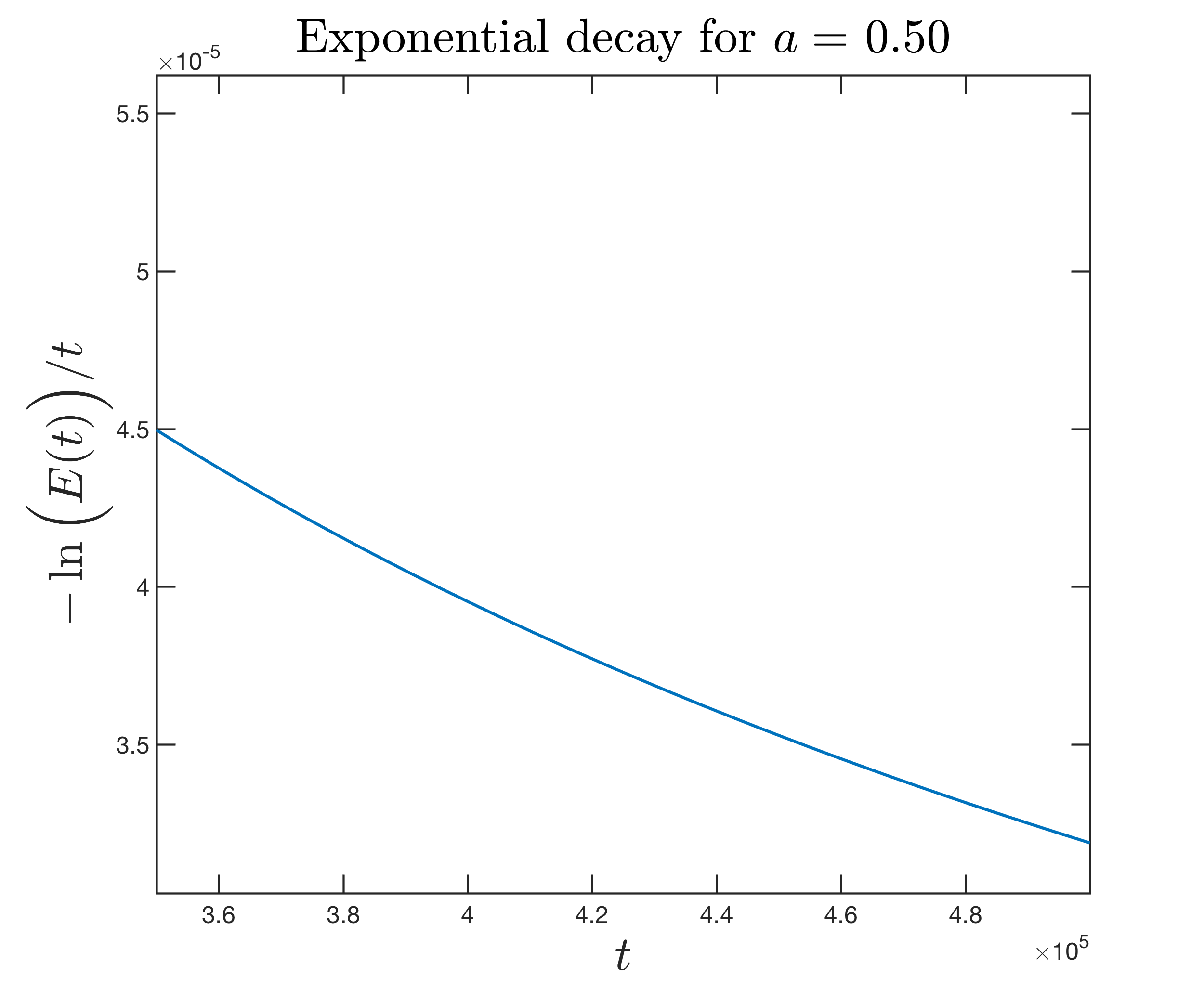}}
	\subcaptionbox{Polynomial decay in $1/t$?\label{1-b4-c5-a1/2}}
	{\includegraphics[width=0.4\textwidth]{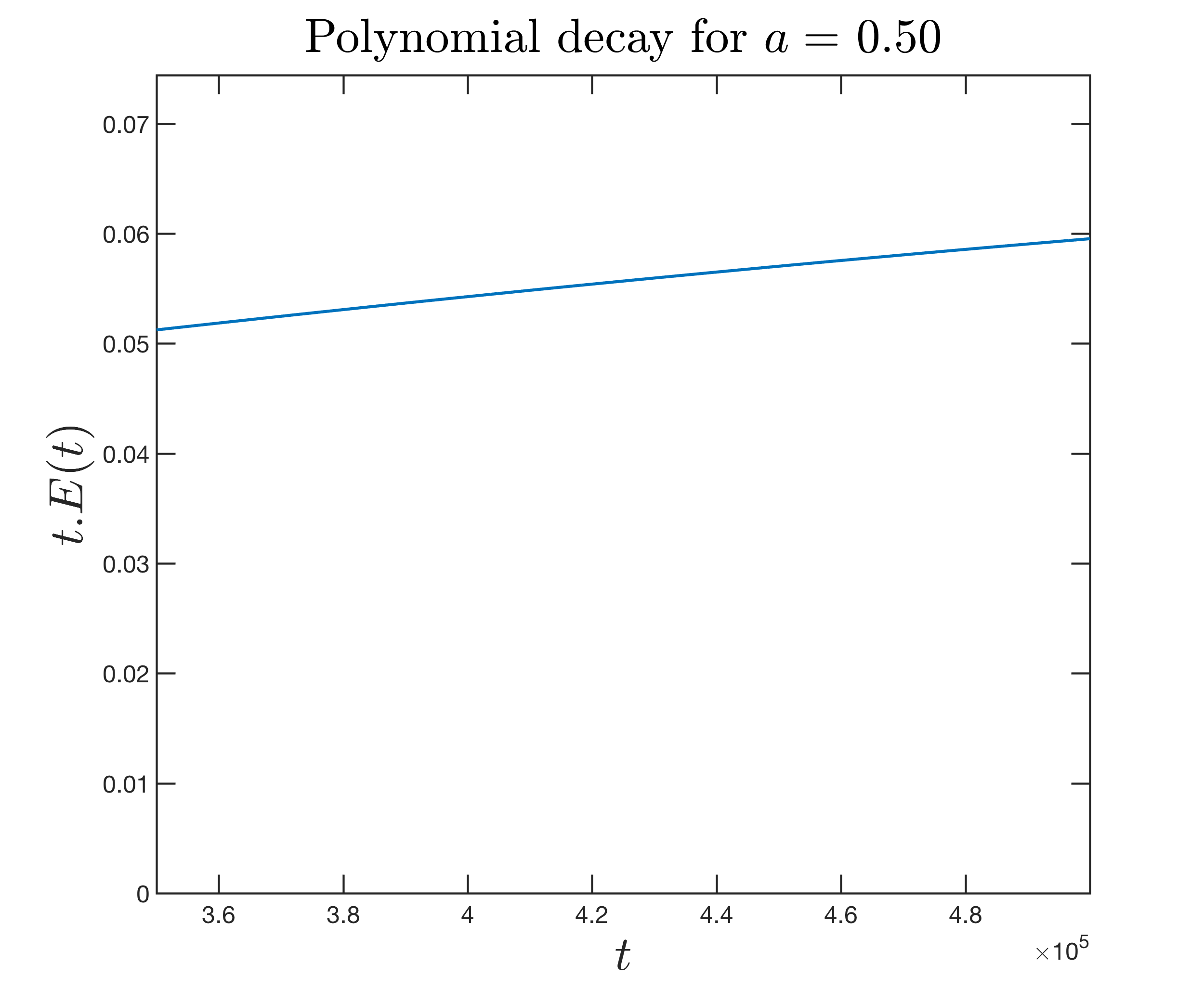}}\\[0.01\textheight]
	\subcaptionbox{Which exponent if polynomial decay?\label{Poly-b4-c5-a1/2}}
	{\includegraphics[width=0.4\textwidth]{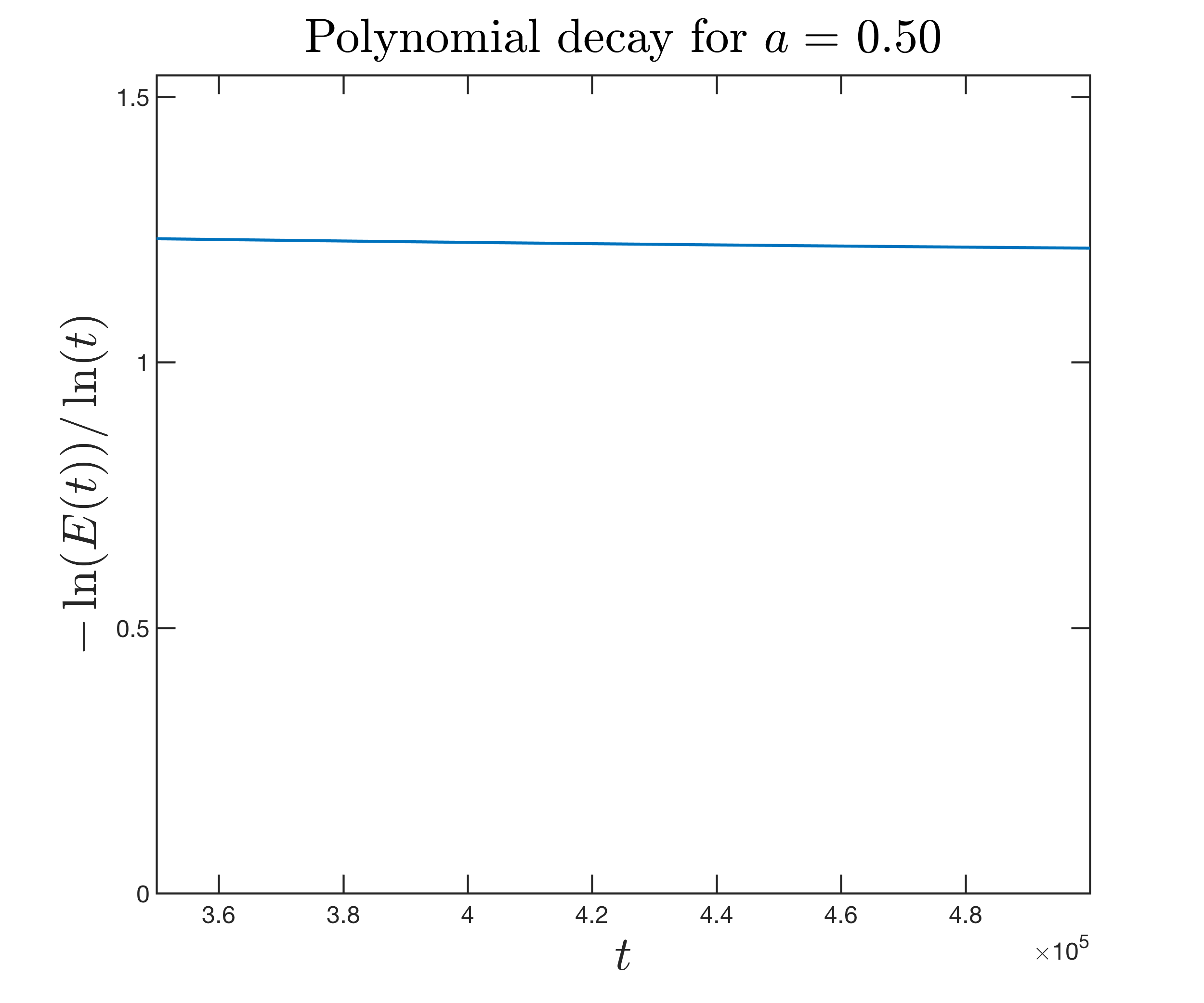}}
	\subcaptionbox{Final time profile.\label{b4-c5-a1/2}}
	{\includegraphics[width=0.4\textwidth]{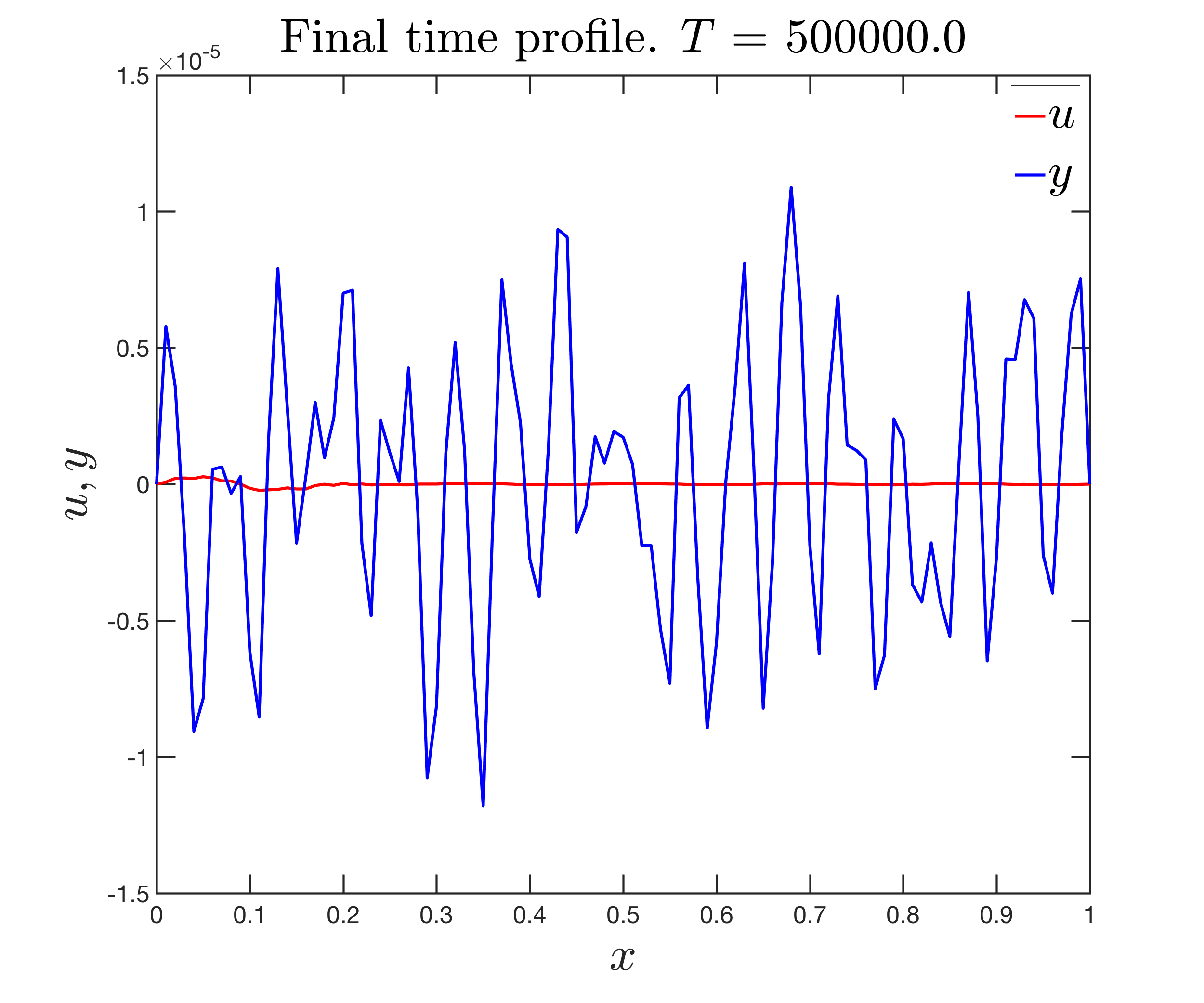}}
	\captionsetup{justification=centering}
	\caption{Long time behavior when $\omega_{b} \cap \omega_{c_{+}} = \emptyset$. \\
		\footnotesize{$b = b_{4}(x) = \mathds{1}_{[0.1,0.2]}(x)$ and $c = c_{5}(x)= \mathds{1}_{[0.4,0.6]}(x)$}}
	\label{Convergence-b4-c5-a1/2}
\end{figure}

To explore the speed of convergence to zero, we have plotted in figure \ref{Convergence-b4-c5-a1/2} $-\ln\big(E(t)\big) /t$ , $t \cdot E(t)$ and finally  $-\ln\big(E(t)\big) / \ln(t)$ versus $t$. Figure \ref{Exp-b4-c5-a1/2} shows clearly that $-\ln\big(E(t)\big) /t$ tends to zero slower than an exponential. But figure \ref{1-b4-c5-a1/2} shows that 
$E(t)$ tends to zero { faster} than $1/t$.  Finally figure \ref{Poly-b4-c5-a1/2} shows that $E(t)$ tends to zero as $1/t^{\alpha}$ with $\alpha \simeq 1.25$. 

Again, the final time profile presented in figure \ref{b4-c5-a1/2} confirms that $u$ and $y$ are small but it shows also that high frequencies for the unknown $y$ are not completely controlled.
This result was  not predicted by the theoretical results.

\begin{figure}[H]
	\centering
	\subcaptionbox{Final time $T = 500$.\label{Energy-small-b5-c4-a1/2}}
	{\includegraphics[width=0.4\textwidth]{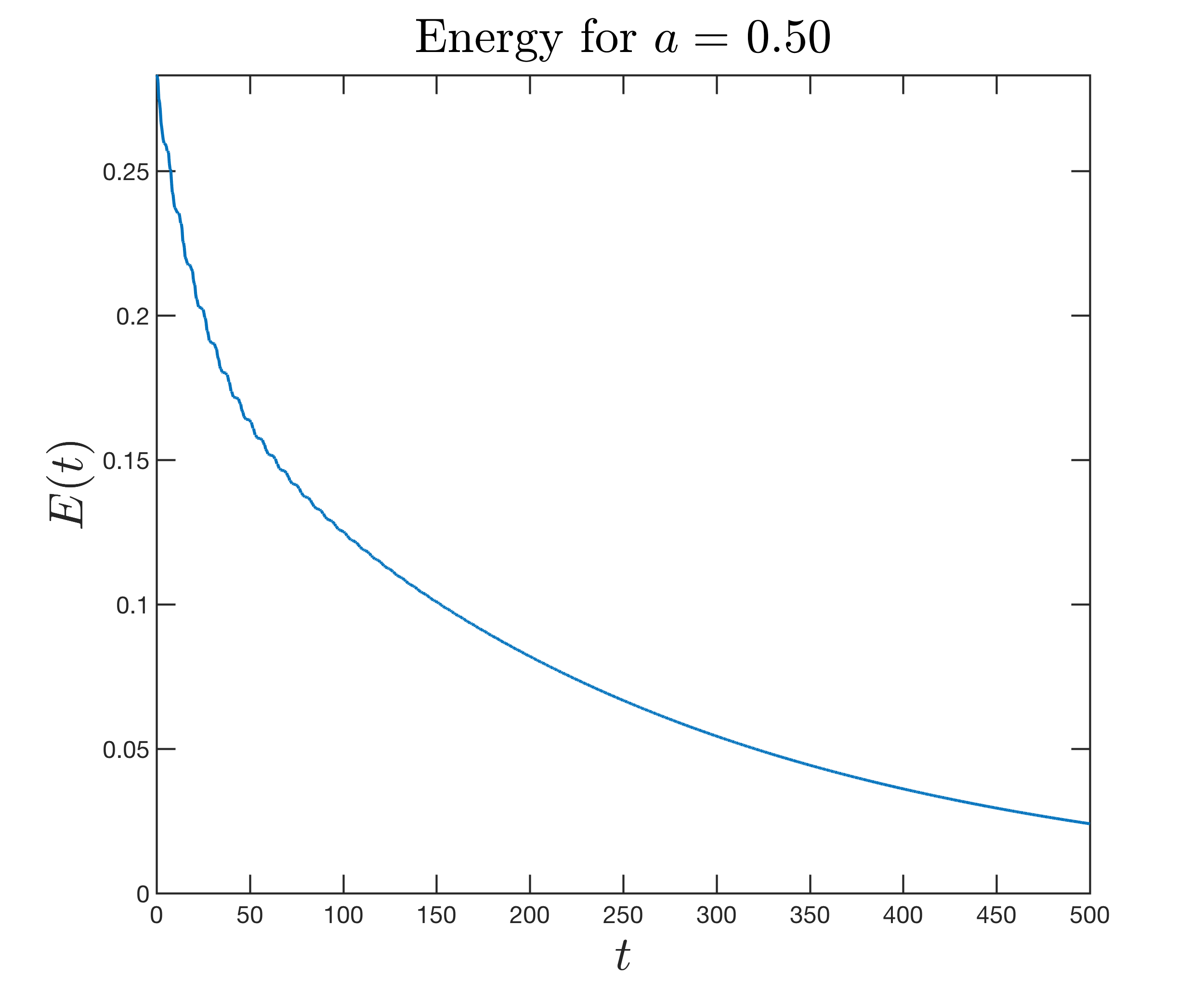}}
	\subcaptionbox{Final time $T= 500~000$.\label{Energy-b5-c4-a1/2}}
	{\includegraphics[width=0.4\textwidth]{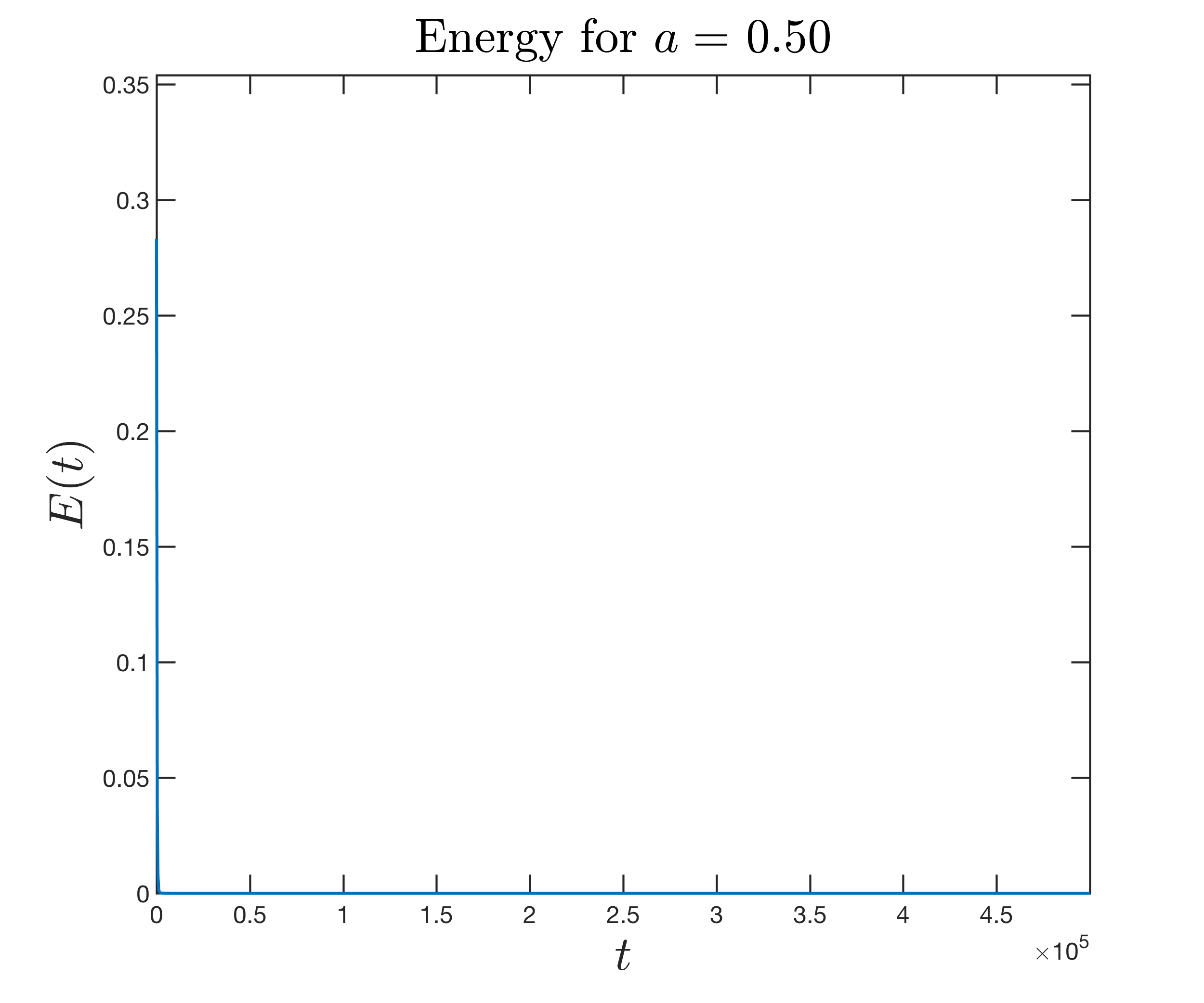}}
	\captionsetup{justification=centering}
	\caption{Energy when $\omega_{b} \cap \omega_{c_{+}} = \emptyset$.\\
		\footnotesize{$b = b_{4}(x) = \mathds{1}_{[0.1,0.2]}(x)$ and $c = c_{5}(x)= \mathds{1}_{[0.4,0.6]}(x)$}}
	\label{Ener-b5-c4-a1/2}
\end{figure}

So we decided to confirm this  behavior by choosing $b = b_{5}(x) = \mathds{1}_{[0.4,0.6]}(x)$ and $c = c_{4}(x)= \mathds{1}_{[0.1,0.2]}(x)$.
Again, when taking as final time $T = 500$, it seems that the energy does not tend to zero  as shown in figure \ref{Energy-small-b5-c4-a1/2}.
Taking as final time $T = 500~000$, figure \ref{Energy-b5-c4-a1/2}  shows that the energy goes finally to zero.

\begin{figure}[H]
	\centering
	\subcaptionbox{Exponential decay?\label{Exp-b5-c4-a1/2}}
	{\includegraphics[width=0.4\textwidth]{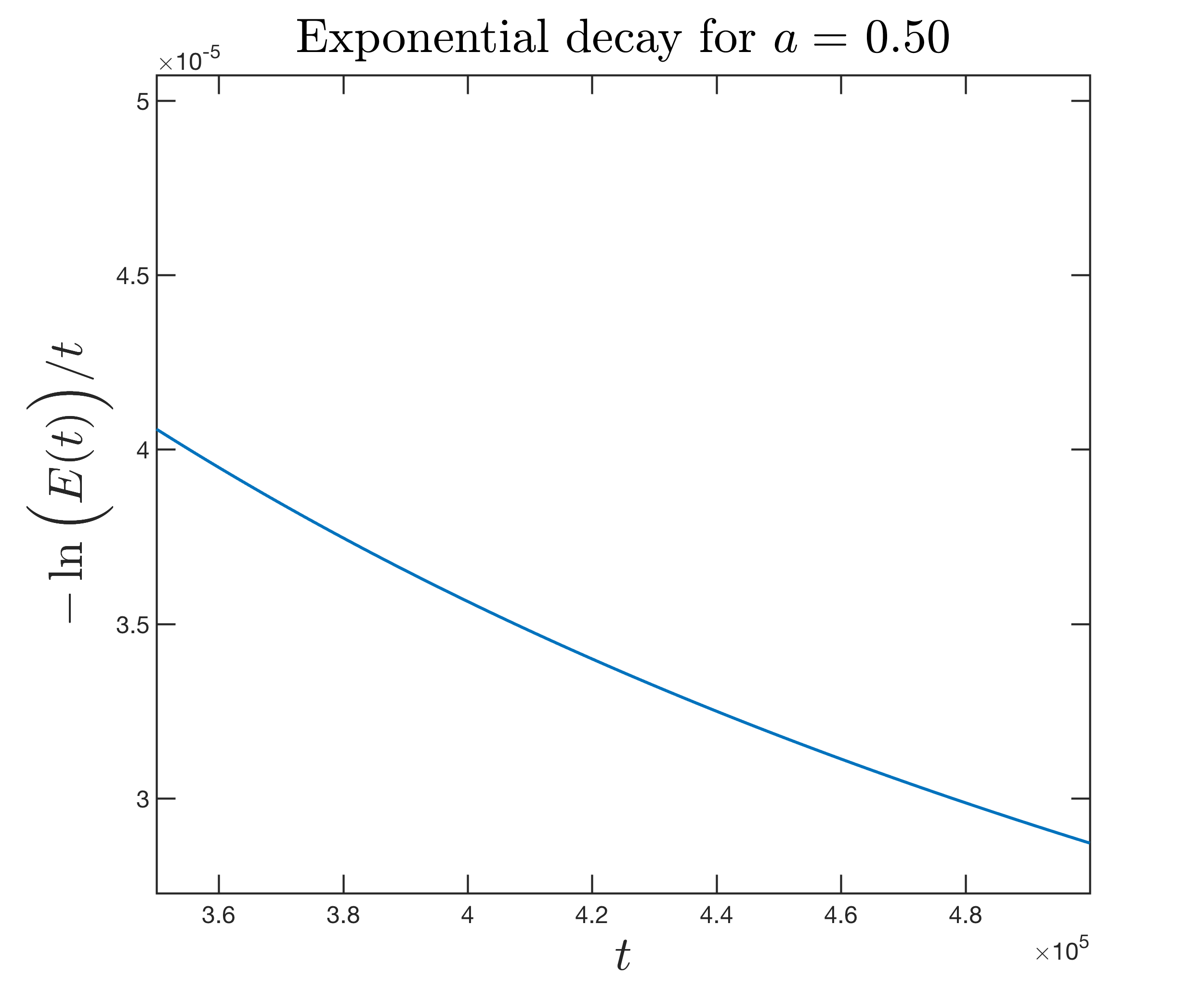}}
	\subcaptionbox{Polynomial decay in $1/t$?\label{1-b5-c4-a1/2}}
	{\includegraphics[width=0.4\textwidth]{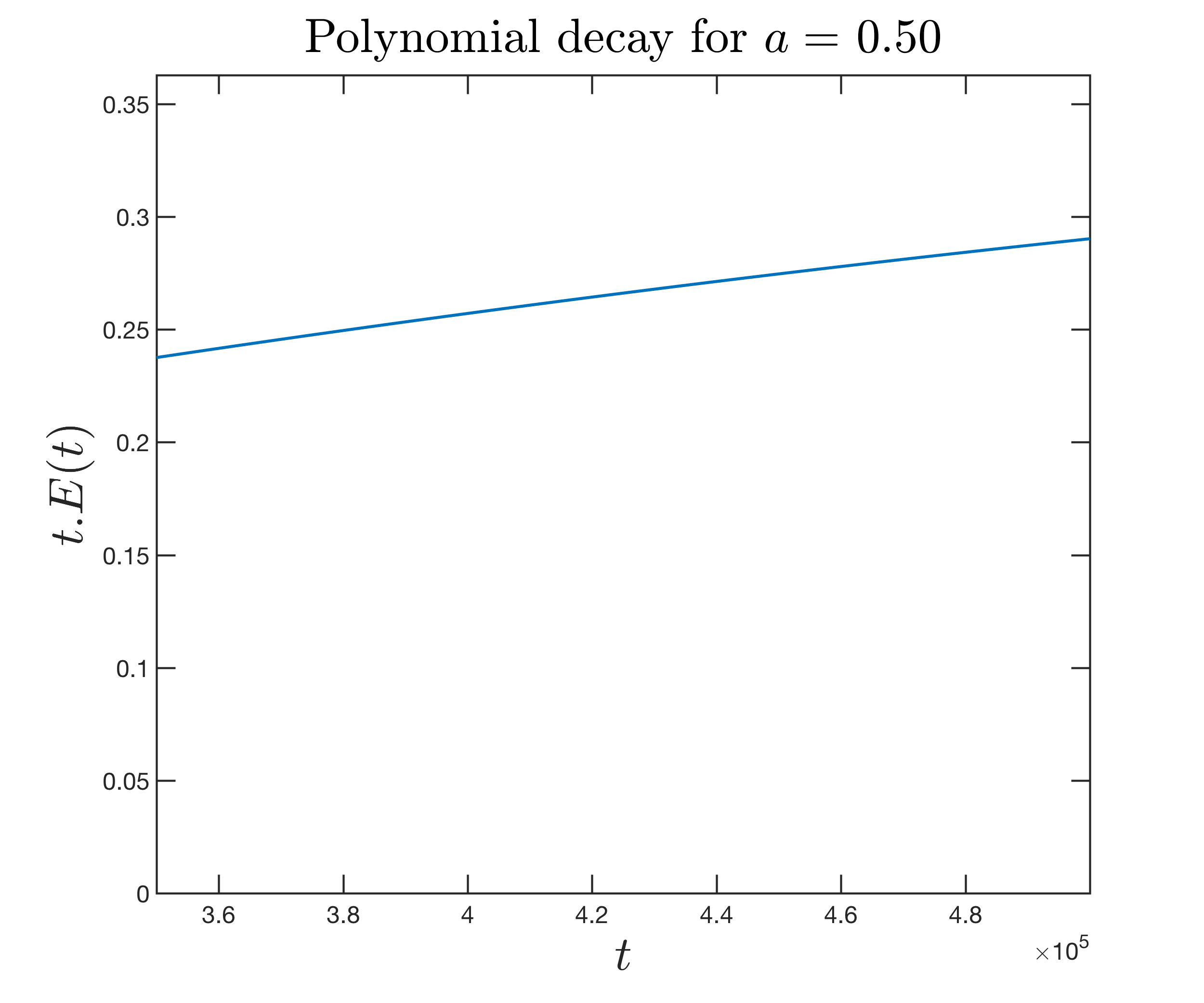}}\\[0.01\textheight]
	\subcaptionbox{Which exponent if polynomial decay?\label{Poly-b5-c4-a1/2}}
	{\includegraphics[width=0.4\textwidth]{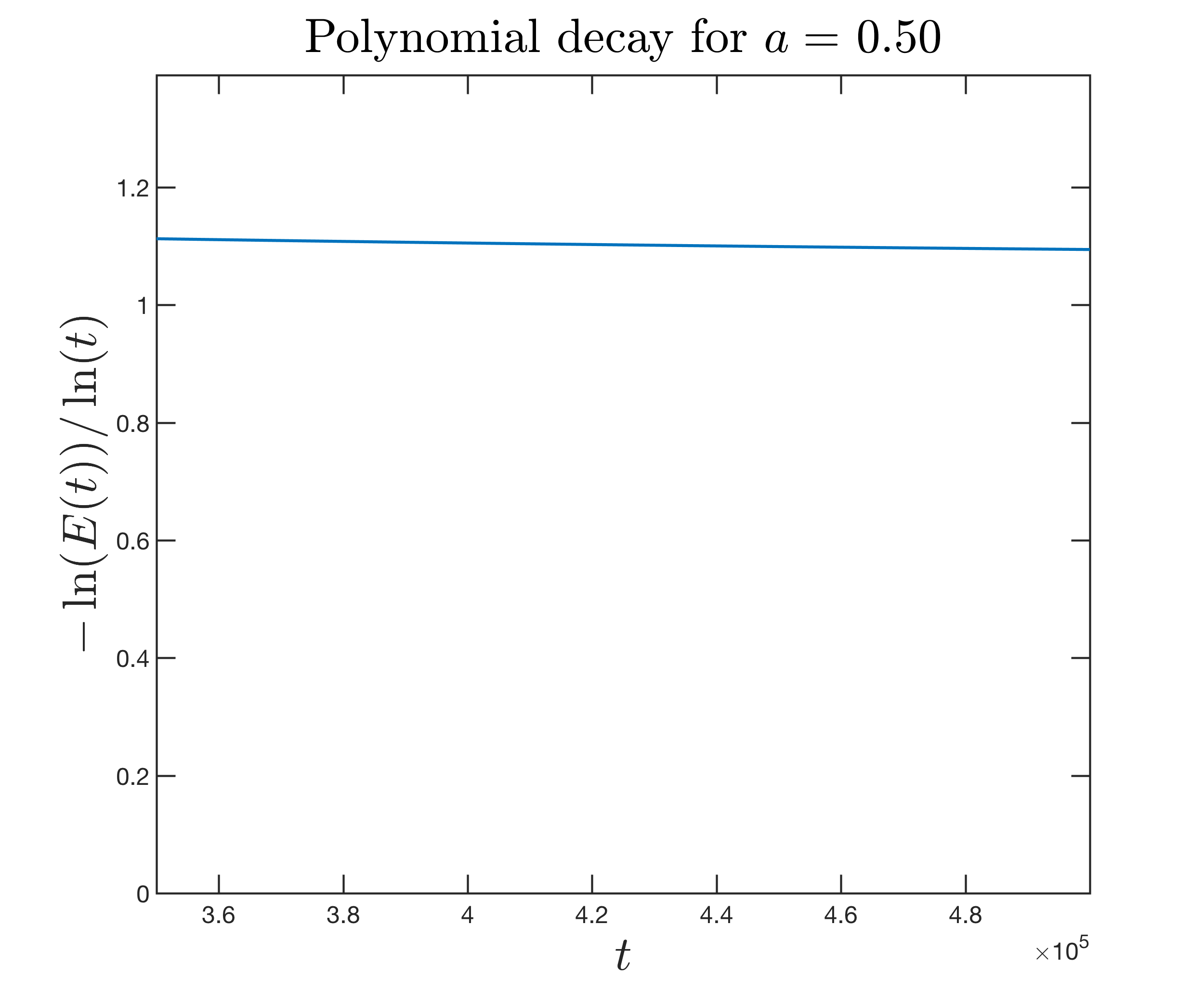}}
	\subcaptionbox{Final time profile.\label{b5-c4-a1/2}}
	{\includegraphics[width=0.4\textwidth]{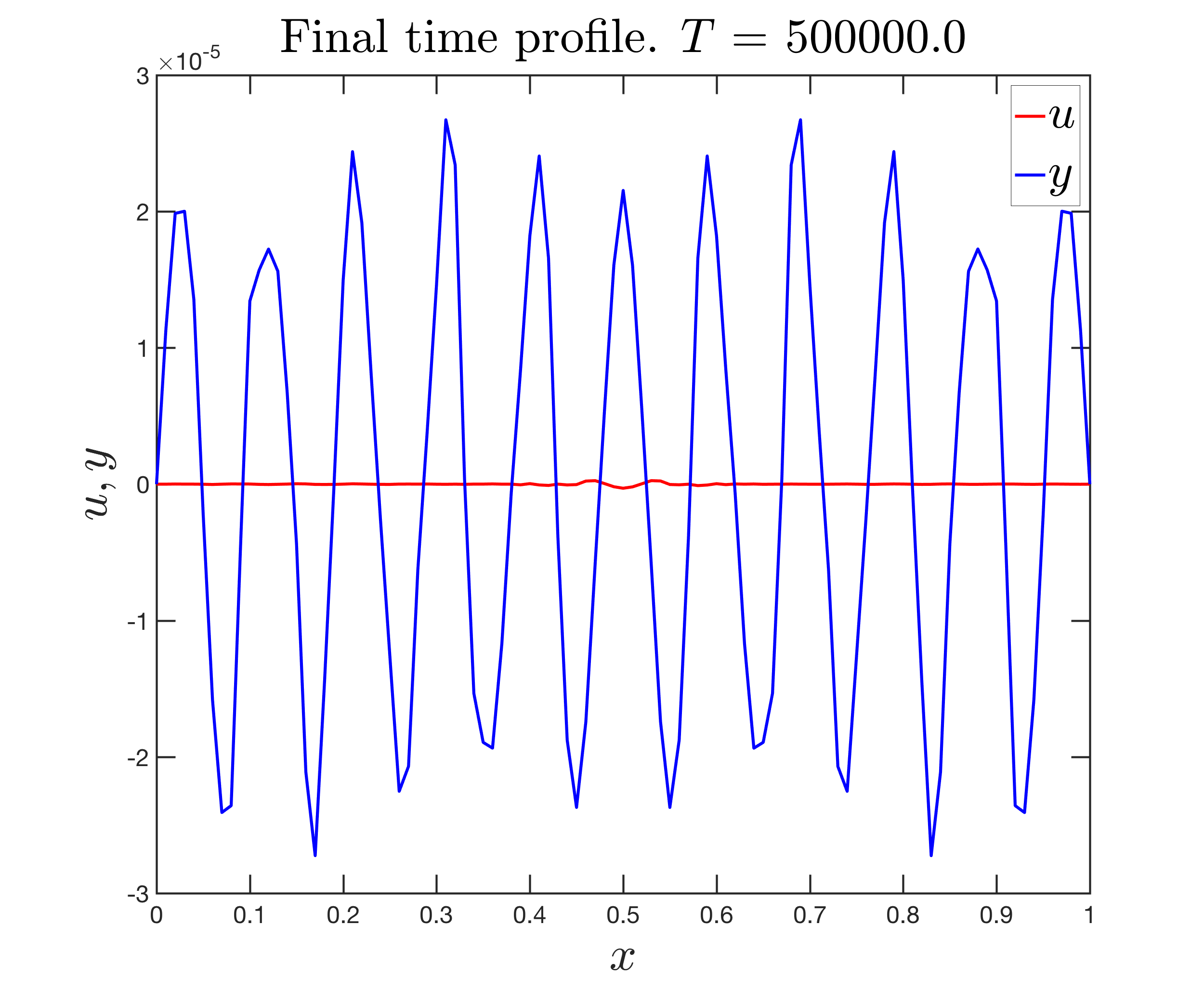}}
	\captionsetup{justification=centering}
	\caption{Long time behavior when $\omega_{b} \cap \omega_{c_{+}} = \emptyset$. \\
		\footnotesize{$b = b_{5}(x) = \mathds{1}_{[0.4,0.6]}(x)$ and $c = c_{4}(x)= \mathds{1}_{[0.1,0.2]}(x)$}}
	\label{Convergence-b5-c4-a1/2}
\end{figure}

To explore the speed of convergence to zero, we have plotted in figure \ref{Convergence-b5-c4-a1/2} $-\ln\big(E(t)\big) /t$ , $t \cdot E(t)$ and finally  $-\ln\big(E(t)\big) / \ln(t)$ versus $t$. 
Figure \ref{Exp-b5-c4-a1/2} shows clearly that $-\ln\big(E(t)\big) /t$ tends to zero and it permits to conclude that $E(t)$ tends to zero slower than an exponential but figure \ref{1-b5-c4-a1/2} shows that 
$E(t)$ tends to zero { faster} than $1/t$. 
Finally figure \ref{Poly-b5-c4-a1/2} shows that $E(t)$ tends to zero as $1/t^{\alpha}$ with $\alpha \simeq 1.15$. 

Again, the final time profile presented in figure \ref{b5-c4-a1/2} confirms that $u$ and $y$ are small but it shows also that high frequencies for the unknown $y$ are not completely controlled.

\begin{remark}
	The final time profile presented in figure \ref{b4-c3-a1/2} , figure \ref{b4-c5-a1/2} and figure \ref{b5-c4-a1/2} confirms that $u$ and $y$ are small 
	but it shows also that high frequencies for the unknown $y$ are not completely controlled.
\end{remark}

\begin{remark}
	When the propagation speeds are not equal, the solution $(u,y)$ has the same behavior no matter if $a > 1$ or $a < 1$. The polynomial convergence is numerically better than $1/t$ but
	it will be probably be $1/t$ for greater time. For reason of computation time, we did not perform very long simulation to confirm.
\end{remark}

\subsubsection*{Acknowledgments}  
The authors are grateful to the anonymous referees and the editor  for their valuable comments and useful suggestions.

The authors thanks professor Kais Ammari for their valuable discussions and comments.

Amina Mortada and Chiraz Kassem would like to thank the AUF agency for its support in the framework of the PCSI project untitled {\it Theoretical and Numerical Study of Some Mathematical Problems and Applications}

Ali Wehbe would like to thank the CNRS and the LAMA laboratory of Mathematics of the Université Savoie Mont Blanc for their supports. 

\begin{figure}
	\centering
	{\includegraphics[width=0.45\textwidth]{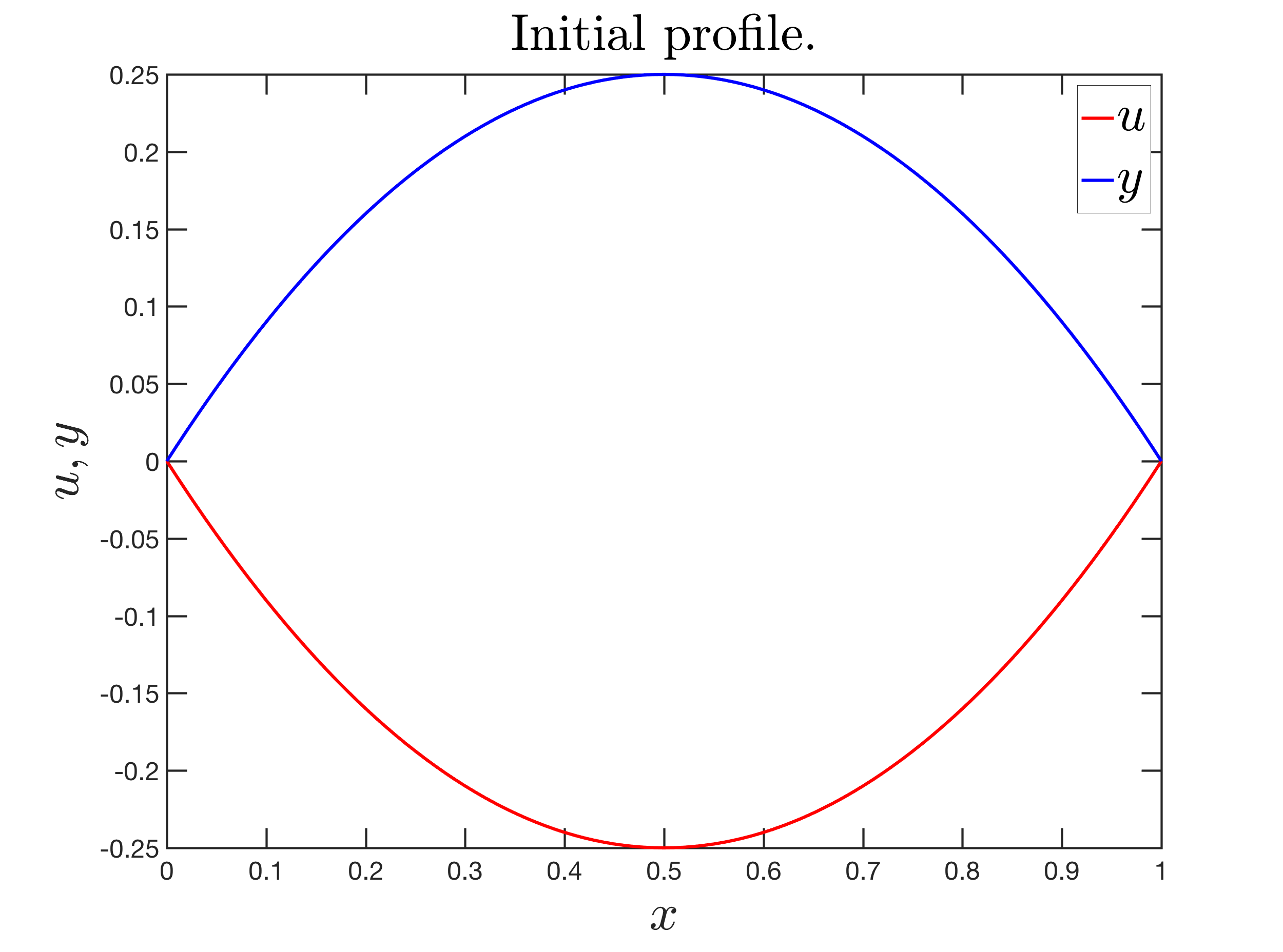}}
	{\includegraphics[width=0.45\textwidth]{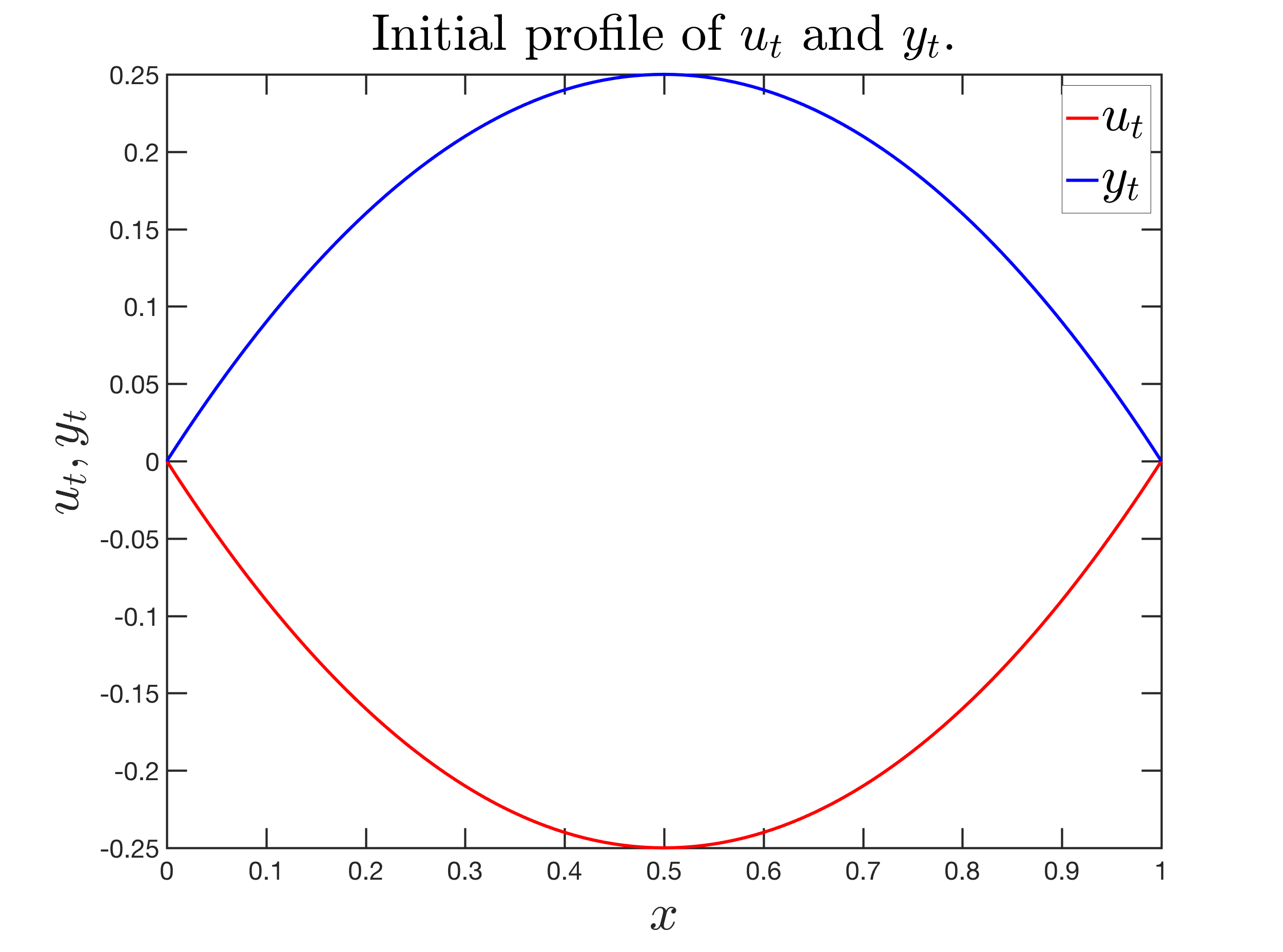}}\\
	\caption{Initial profiles}
\end{figure}

\vspace{-2.5cm}

\pagebreak
\pagebreak
\pagebreak
\pagebreak

\vspace*{-0.5cm}

\pagebreak

\vspace*{-0.5cm}

\pagebreak

\vspace*{-0.5cm}

\pagebreak

\vspace*{-0.5cm}

\pagebreak

\vspace*{-0.5cm}

\pagebreak

\vspace*{-0.5cm}

\bibliographystyle{plain}

\end{document}